\documentclass[11pt,reqno,a4paper]{amsart}
\usepackage[fleqn,tbtags]{mathtools}
\usepackage{amssymb}
\usepackage{amsthm}
\usepackage[T1]{fontenc}
\usepackage{lmodern}
\usepackage{graphicx}
\usepackage[shortlabels]{enumitem}
\setlist{topsep=4pt plus 1pt minus 1pt,itemsep=2pt plus 1pt,
         parsep=0pt,partopsep=0pt}
\usepackage{microtype}
\usepackage[nodisplayskipstretch]{setspace}
\usepackage{xcolor}
\usepackage{tikz}
\usetikzlibrary{arrows.meta,positioning,calc}
\usepackage[a4paper,margin=25mm]{geometry}
\AtBeginDocument{%
  \setlength{\abovedisplayskip}{8pt plus 2pt minus 2pt}%
  \setlength{\belowdisplayskip}{8pt plus 2pt minus 2pt}%
  \setlength{\abovedisplayshortskip}{3pt plus 2pt}%
  \setlength{\belowdisplayshortskip}{5pt plus 2pt minus 1pt}%
}
\usepackage[colorlinks=true,linkcolor=red,citecolor=blue,urlcolor=blue,
            hypertexnames=false]{hyperref}
\usepackage[nameinlink,capitalize,noabbrev]{cleveref}

\theoremstyle{plain}
\newtheorem{theorem}{Theorem}[section]
\newtheorem{proposition}[theorem]{Proposition}
\newtheorem{lemma}[theorem]{Lemma}
\newtheorem{corollary}[theorem]{Corollary}
\theoremstyle{definition}
\newtheorem{definition}[theorem]{Definition}
\theoremstyle{remark}
\newtheorem{remark}[theorem]{Remark}

\theoremstyle{plain}

\newcommand{\ddc}{dd^c}
\newcommand{\Supp}{\operatorname{Supp}}

\newcommand{\D}{{\mathbb{D}}}
\newcommand{\PP}{{\mathbb{P}}}

\newcommand{\dA}{dA}
\DeclareMathOperator{\Ric}{Ric}
\DeclareMathOperator{\tr}{tr}
\newcommand{\dd}{\,d}
\allowdisplaybreaks[1]
\numberwithin{equation}{section}
\hypersetup{pdftitle={A K\"ahler Miyaoka--Mori Theorem}}
\title[Rational curves on compact K\"ahler manifolds]{Disc functionals and rational curves}
\date{}

\author[Y. Li]{Yinji Li}
\address{Yinji Li: Institute of Mathematics\\Academy of Mathematics and Systems Science\\Chinese Academy of
		Sciences\\Beijing\\100190\\P. R. China}
\email{liyinji@amss.ac.cn}
\author[H. Sun]{Haoyuan Sun}
\address{Haoyuan Sun: School of Mathematical Sciences\\ Beijing Normal University\\ Beijing 100875\\ P. R. China}
\email{202531130037@mail.bnu.edu.cn}
\author[Z. Wang]{Zhiwei Wang}
\address{Zhiwei Wang: Laboratory of Mathematics and Complex Systems (Ministry of Education)\\ School of Mathematical Sciences\\ Beijing Normal University\\ Beijing 100875\\ P. R. China}
\email{zhiwei@bnu.edu.cn}

\author[X. Zhou]{Xiangyu Zhou}
\address{Xiangyu Zhou: Institute of Mathematics\\ Academy of Mathematics and Systems Science\\
	and Hua Loo--Keng Key Laboratory of Mathematics\\ Chinese Academy of Sciences\\ Beijing 100190\\ P. R. China}
\email{xyzhou@math.ac.cn}

\subjclass[2020]{32Q15; 32U35}
\keywords{Discs functionals, rational curves, K\"ahler manifolds}

\begin{document}
\raggedbottom

\begin{abstract}
We prove a criterion for the existence of rational curves on compact K\"ahler manifolds using disc functionals. As an application, we obtain the existence of rational curves with an optimal degree bound on compact K\"ahler pair with non-nef canonical class. 
\end{abstract}

\maketitle
\tableofcontents

\section{Introduction}
\label{sec:introduction}

Our goal is to construct rational curves on normal compact K\"ahler pairs with non-nef canonical class from a point of view of several complex variables.
Construction of rational curves is an important problem in algebraic geometry and analytic geometry: rational curves of negative canonical degree generate
the extremal rays to be contracted, see for example,
\cite{KM98,BCHM10,HP16intro,HP24intro,HLX26intro,HX26intro}.
Mori introduced bend-and-break technique to construct rational curves in his proof of Hartshorne's conjecture
on ample tangent bundles \cite{Mor79} and later used it to study
extremal rays \cite{Mor82}. This method is a combination of deformation theory, reduction to positive characteristic and Frobenius morphism.
Using bend-and-break, Miyaoka and Mori \cite{MM86} showed
that if $K_X\cdot C<0$ for a curve $C$ on a smooth projective
$n$-fold, then for any ample divisor $H$ and any $x\in C$ there
is a rational curve $R$ through $x$ with
$H\cdot R\le 2n(H\cdot C)/(-K_X\cdot C)$; see also
\cite{Kol96} when $H$ is a nef real divisor. Recently,
Jovinelly, Lehmann and Riedl \cite{JLR} sharpened
the coefficient to the optimal $n+1$ by breaking off rational
components at several marked points at once, which they described as bend-and-shatter. 

On the analytic side, Siu and Yau \cite{SY80} constructed rational
curves from harmonic spheres in their proof of the Frankel conjecture.
Siu later proposed using multiplier ideals and complex Monge--Amp\`ere
equations to construct rational curves on Fano manifolds \cite{Siu09}.
If $X$ is a projective manifold, Boucksom,
Demailly, P\u{a}un and Peternell
\cite{BDPP13} proved the uniruled criterion: $X$ is uniruled if and only if $K_X$ is not pseudoeffective.
They obtained this
characterization by combining their cone duality theorem with Mori's reduction to positive characteristic. 
Ou \cite{Ou} proved the characterization in full generality when $X$ is compact K\"ahler by foliation theory and geometric property of optimal $L^2$ extension \cite{GZ15a,GZ15b,GZ17,ZZ16}.


\vspace{0.4cm}
\subsection{Disc functional and criterion of existence of rational curves on compact K\"ahler manifold}

Disc functional plays an important role in several complex variables, see for example \cite{Pol91,Pol93,Ros03,DF12}.
We refer the readers to \cite[Section~16]{Kis00} for a nice survey. 
The construction of plurisubharmonic functions by taking infima links disc envelope functions to minimum principles. The extended future tube conjecture was proved in~\cite{Zhou98}, by discovering a deep connection between the conjecture and the Kiselman--Loeb minimum principle,
as explained in~\cite[Section~4]{Zhou02}.

\vspace{0.15cm}
Let $(X,\omega)$ be a connected compact K\"ahler manifold of dimension $n$ and $\theta:=\operatorname{Ric}(\omega)$, so that $[-\theta]=2\pi c_1(K_X)$.
For a holomorphic disc $f:\D\to X$, its area and weighted area functionals are defined as follows.
\[
 A(f):=\int_\D f^*\omega,\qquad
 W(f):=\int_\D(1-|z|^2)\,f^*\omega.
\]
Let $\rho$ be a smooth closed real $(1,1)$-form on $X$, we define the disc functional $Q_\rho$ by
\[
Q_\rho(f):=\frac{1}{2\pi}\int_\D \log(1/|z|)f^*\rho\,.
\]
At a given point $x\in X$, we set
\[
\mathcal A(X,x):=\Big\{f:\overline\D\to X\text{ holomorphic with }f(0)=x\Big\}.
\]
Using the disc functional $Q$, we obtain a criterion of the existence of rational curves on compact K\"ahler manifold.

\begin{theorem}
\label{thm:pointwise-criterion}
Let $\delta>0$. 
If for some $c>0$ and $q(f):=Q_{\theta-\delta\omega}(f)-cW(f)$ one has
\begin{equation}\label{eq:pointwise-escape}
 \sup_{f\in\mathcal A(X,x)}q(f)=+\infty,
\end{equation}
then there is a rational curve $R\subset X$ through $x$ with
\[
 [\omega]\cdot R\le\frac{2\pi(n+1)}{\delta}.
\]
\end{theorem}

On one hand, the logarithmic term $\log(1/|z|)$ in disc functional $Q$ "squeezes out" a rational curves at $z=0$, as explained in Section~\ref{subsec:stategy}.
On the other hand, the disc functional $Q_\rho$ arises naturally when studying the positivity of $(1,1)$-class $[\rho]$.

\begin{theorem}
\label{thm:disc functional}
    Let $X$ be a connected complex manifold and $\rho$ be a smooth closed real $(1,1)$-form. The disc envelope function
    \[
    x\longrightarrow\inf_{f\in \mathcal A(X,x)}Q_\rho(f)
    \]
    is either $\rho$-plurisubharmonic or identically $-\infty$.
\end{theorem}
For a similar result, see
Magn\'usson's paper \cite{Magnusson}. We provide a complete and detailed proof of Theorem~\ref{thm:disc functional} in Section~\ref{sec:nonpsef}.

Now we suppose $X$ is a compact K\"ahler manifold and $K_X$ is not pseudoeffective. Since the
pseudoeffective cone is closed, for some $t>0$ the class $ [-\theta+2\pi t\omega]=2\pi[K_X+t\omega]$
is still not pseudoeffective, namely $\operatorname{PSH}(X,-\theta+2\pi t\omega)=\emptyset$. Theorem~\ref{thm:disc functional} yields that
\[
 \inf_{f\in\mathcal A(X,x)}Q_{-\theta+2\pi t\omega}(f)=-\infty
 \qquad\text{for every }x\in X.
\]
$2\log(1/|z|)\ge 1-|z|^2$ gives $4\pi Q_\omega(f)\ge W(f)$.
For any $0<\delta<2\pi t$, we set $c=(2\pi t-\delta)/(4\pi)$ and obtain
\[
 Q_{\theta-\delta\omega}(f)-cW(f)
 \ge Q_{\theta-2\pi t\omega}(f).
\]
Therefore \eqref{eq:pointwise-escape} holds for every
$x\in X$ and
Theorem~\ref{thm:pointwise-criterion} gives
a rational curve through every point of $X$. 
This recovers the
implication from non-pseudoeffectivity to uniruledness in Ou's characterization of uniruled compact K\"ahler manifolds \cite{Ou}. 

We call $t_*:=\inf\{t>0:[K_X+t\omega]\text{ is pseudoeffective}\}$ the pseudoeffective threshold of $K_X$. In Appendix~\ref{app:psef threshold} we construct a compact K\"ahler manifold $(X,\omega)$ satisfying
\[
\sup_{f\in\mathcal A(X,x)}
Q_{\theta-\delta\omega}(f)
=
\begin{cases}
+\infty,&0<\delta<2\pi t_*,\\
0,&\delta\ge 2\pi t_*.
\end{cases}
\qquad \forall\,x\in X.
\]

With a little more work, this argument extends to compact K\"ahler pairs. By a pair $(X,B)$, we always mean that $X$ is normal, $B\geq0$ is an effective $\mathbb{R}$-divisor and $K_X+B$ is $\mathbb{R}$-Cartier.
\begin{theorem}\label{thm:intro-normal-kahler-uniruled}
Let $(X,B)$ be a compact K\"ahler pair of dimension $n$. Let $\omega$ be a K\"ahler form on $X$.
If $[K_X+B+t\omega]$ is not pseudoeffective for some $t>0$,
then for every $x\in X$ there is a rational curve $R\subset X$
through $x$ such that
\begin{equation}
 [\omega]\cdot R\leq\frac{n+1}{t}.
\end{equation}
In particular, if $K_X+B$ is not pseudoeffective, then $X$ is
uniruled.
\end{theorem}


\vspace{0.25cm}
\subsection{Rational curves on compact K\"ahler space with non-nef canonical class}
Theorem~\ref{thm:pointwise-criterion} also produces rational curves when $K_X+B$ is not nef.



A $\mathbb{R}$-Cartier divisor $H$ is nef on a projective variety $X$ if $H\cdot C\geq0$ for every irreducible curve $C\subset X$.
Over algebraically closed field of arbitrary characteristic, Jovinelly, Lehmann and Riedl \cite{JLR} obtained the Miyaoka--Mori theorem with optimal degree bound:
\begin{theorem}\cite[Theorem~1.1]{JLR}\label{thm:intro-alg-MM}
Let $X$ be a projective variety and $H$ be a nef $\mathbb{R}$-Cartier divisor on $X$. Suppose there exists an irreducible curve $C\subset X$ contained in the smooth locus of $X$ such that $K_X\cdot C<0$, then for every closed point $x\in C$, there exists a rational curve $R$ containing $x$ such that
\[
H\cdot R\leq\frac{(n+1)H\cdot C}{-K_X\cdot C}.
\]
\end{theorem}
Recently, Du--Guo--Xie \cite{DGX} proposed an analytic approach to constructing rational curves on Fano manifold and
Guo--Xie \cite{GXJets} obtained rational curves on projective manifold with non-nef canonical bundle in characteristic zero by a jet technique.

\vspace{0.15cm}
In the analytic setting, 
a $(1,1)$-class $\gamma\in H^{1,1}_{\rm BC}(X)$ on a normal compact K\"ahler space $X$ is called nef if it lies in the closure of the K\"ahler cone. 
Demailly and P\u{a}un \cite{DP04} extended the  Kleiman--Nakai--Moishezon's criterion \cite{Nak63,Moi64,Kle66} to compact K\"ahler manifolds. Das--Hacon--P\u{a}un \cite{DHP24} extended the criterion to compact K\"ahler normal spaces.
$\gamma$ is nef if and only if
\[
 \int_Y \gamma\wedge\omega^{\dim Y-1}\ge 0
\]
for every irreducible subvariety $Y\subset X$ of
positive dimension and every K\"ahler form $\omega$ on $X$.
We call an irreducible subvariety $Y\subset X$ a non-nef subvariety of $\gamma$ if $\int_Y\gamma\wedge\omega^{\dim Y-1}<0$ for some K\"ahler form $\omega$.
This naturally suggests replacing the $K_X$-negative curve in the
Miyaoka--Mori Theorem~\ref{thm:intro-alg-MM} by the non-nef subvariety $Y$ of $K_X+B$.

\begin{theorem}[=Theorem~\ref{thm:normal-kahler-rational-curves}]\label{thm:intro-normal-kahler-rational-curves}
Let $(X,B)$ be a compact K\"ahler pair. Let $Y\subset X$ be a $d$-dimensional subvariety and
$\omega_1,\ldots,\omega_{d-1}$ be nef $(1,1)$-classes on $X$.
Set
\[
 L_Y(\gamma):=
 \int_Y\gamma\wedge\omega_1\wedge\cdots\wedge\omega_{d-1},\qquad\gamma\in H^{1,1}_{\rm BC}(X).
\]
Assume that $L_Y(K_X+B)<0$ and $Y$ is not contained in $X_{\rm sing}\cup\Supp B$.
Then for every nef $(1,1)$-class $\alpha$ on $X$ and every
$x\in Y$, there is a rational curve $R\subset X$
through $x$ such that
\begin{equation}
 \alpha\cdot R\leq
 (n+1)\frac{L_Y(\alpha)}{-L_Y(K_X+B)}.
\end{equation}
\end{theorem}


Together, Theorems~\ref{thm:intro-normal-kahler-uniruled}
and~\ref{thm:intro-normal-kahler-rational-curves} show that
non-pseudoeffectivity of $K_X+B$ yields rational curves through
every point of $X$, while non-nefness yields rational curves
through every point of every non-nef subvariety of $K_X+B$ not contained in $X_{\rm sing}\cup \Supp B$.
The conclusion may fail when $Y\subset X_{\rm sing}\cup\Supp B$.
Indeed, there exists normal projective Gorenstein $X$
such that $K_X$ is not nef and every non-nef subvariety of $K_X$ is contained in $X_{\rm sing}$, yet $X$ contains no rational curves.

\vspace{0.2cm}
For a nef $(1,1)$-class $\alpha$ on a normal compact K\"ahler space $X$, its null locus
$\operatorname{Null}(\alpha)$ is the union of all positive-dimensional
irreducible subvarieties $Y\subset X$ such that
$\int_Y\alpha^{\dim Y}=0$. If $\int_X\alpha^n=0$, then $\operatorname{Null}(\alpha)=X$.
If $\int_X\alpha^n>0$, Collins and
Tosatti \cite{CT15} showed that
$\operatorname{Null}(\alpha)$ coincides with the non-K\"ahler locus
$E_{\mathrm{nK}}(\alpha)$ when $X$ is smooth and Hacon--P\u{a}un \cite{HP24intro} extended this result to normal compact K\"ahler $X$. 

In K\"ahler minimal model program, construction of $K_X$-negative rational curves is crucial.
With Ou's uniruledness criterion \cite{Ou}, Hacon--P\u{a}un
\cite[Theorem~0.5]{HP24intro} and
Hacon--Li--Xie \cite[Corollary~3.6]{HLX26intro} establish the cone theorem with degree bound $2n$. Our Theorem~\ref{thm:intro-normal-kahler-rational-curves} has the following corollary. 

\begin{corollary}[= Corollary~\ref{cor:nef-threshold}]
\label{cor:intro-nef-threshold}
Let $(X,B)$ be a compact K\"ahler pair and $\omega$ be a K\"ahler form.
Assume $K_X+B$ is not nef and set
\[
 \tau:=\inf\big\{t\ge 0: [K_X+B+t\omega]\text{ is nef}\big\},
 \qquad \alpha:=[K_X+B+\tau\omega].
\]
For positive dimensional subvariety $Y\subset \operatorname{Null}(\alpha)$ which is not contained in $X_{\rm sing}\cup\Supp B$ and every
$x\in Y$, there is a rational curve
$R\subset X$ through $x$ satisfying
\[
 \alpha\cdot R=0,\qquad
 0<-(K_X+B)\cdot R=\tau\int_R\omega\le n+1.
\]
\end{corollary}

\subsection{Kobayashi's canonical ampleness conjecture}

Kobayashi \cite{Kob70} conjectured that
the canonical bundle $K_X$ of a compact hyperbolic K\"ahler manifold is ample. 
\cite{OAI26} recently claimed a proof of this conjecture. 
In Section~\ref{sec:hyperbolic-canonical-bigness}, 
we give an immediate proof based on our setting.

\begin{theorem}[= Theorem~\ref{thm:hyperbolic-canonical-bigness}]
\label{thm:intro-hyperbolic-canonical-ampleness}
Let $X$ be a compact K\"ahler manifold. If every holomorphic map
$\mathbb C\to X$ is constant, then $K_X$ is ample.
\end{theorem}

\subsection{Strategy of the proof of Theorem~\ref{thm:pointwise-criterion}}\label{subsec:stategy}


\vspace{0.2cm}
\noindent\textbf{Step 1. Squeeze out a rational curve by Gromov compactness.}
Assume that there are holomorphic discs $f_j:\D\to X$ satisfying $f_j(0)=x$ and $\sup_j A(f_j)<\infty$. 
With uniformly bounded area, compactness theorem \cite{Gro85,ISInvent,ISCompactness} shows that $f_j:\D\to X$ converge to a stable curve $F:C\to X$.
If further $q(f_j)\to+\infty$, 
the logarithmic term $\log(1/|z|)$ in disc functional $q$ "squeezes out" a bubble at $z=0$, and produces a rational
curve $R$ through $x$ with 
\[
 \int_R\omega\le\liminf_j A(f_j).
\]
See Corollary \ref{cor:quantitative-marked-bubbling} for details. Our task is therefore to find holomorphic discs $f_j$ so that $q(f_j)$ diverge while their areas $A(f_j)$ are uniformly bounded.

\noindent\textbf{Step 2. Minimizer of functional $J_{t,\beta}$.}
We construct the desired discs by minimizing the following functional. For $t,\beta>0$, we set
\[
 J_{t,\beta}(f):=A(f)+te^{-q(f)}-\beta\Psi(f),\qquad f\in\mathcal A(X,x).
\]
We do not give precise definition of $\Psi$ here and simply mention that $\Psi$ ensures the minimizing sequence converges to an actual holomorphic disc $f_{t,\beta}:\D\to X$, not a general stable curve with bubbles, see Proposition~\ref{prop:foundation-attainment}.  

\noindent\textbf{Step 3. Family of holomorphic discs and their variations.}
We construct holomorphic map $F:D_{\varepsilon}\times\D\to X$ so that the central fiber $F(0,\cdot)=f_{t,\beta}$ and we say $\{f_\tau=F(\tau,\cdot)\}$ is a family of holomorphic discs. 
Since $f_{t,\beta}$ minimizes the functional $J_{t,\beta}$, we get $\Delta_\tau J_{t,\beta}(f_\tau)|_{\tau=0}\geq0$. With an appropriate choice of the holomorphic family, which is discussed in the next step, this variation inequality gives uniform bound of area of $f_{t,\beta}$ and the divergence of $q(f_{t,\beta})$, see Proposition~\ref{prop:uniform-variation-bound} and Proposition~\ref{lem:vanishing-correction-minima} respectively.

\noindent\textbf{Step 4. Berndtsson--Rosay's boundary-unitary frame and Forstneri\v c's graph theorem.}
The minimizer $f=f_{t,\beta}$ is proved to be $W^{1,p}$, $2<p<4$, in Proposition~\ref{prop:foundation-attainment}. Sobolev embedding $W^{1,p}(\D)\hookrightarrow C^{0,1-\frac{2}{p}}(\D)$ implies $f$ continuously extend to $\overline\D$. 
We construct a $W^{1,p}$
holomorphic frame $e_1,\ldots,e_n$ of $f^*T^{1,0}_X$ with
$h(e_i,e_j)=\delta_{ij}$ on $\partial\D$, where $h$ is the
Hermitian metric induced by $\omega$. The construction follows
Berndtsson--Rosay \cite[Proposition~6.1]{BR}, see
Theorem~\ref{lem:baseline-frame}. 
An easy computation in Lemma~\ref{lem:exact-kernel} gives tangent fields
$\xi_1,\ldots,\xi_n$ of $f^*T^{1,0}_X$ which remain orthonormal on $\partial\D$ but 
vanish at $0$.
Using the fundamental graph theorem of Forstneri\v c
\cite[Theorem~1.2]{Forstneric}, in Lemma~\ref{lem:actual} we realize tangent field $\xi_j$ as the derivative of holomorphic family $f_{\tau}^j$, namely,
$f_{0}^j=f$, $f_{\tau}^j(0)=x$ and $\partial_\tau f_\tau^j|_{\tau=0}=\xi_j$.

The second variation of $J_{t,\beta}$ appears to be very clean when summed over the families $f^1,\cdots,f^n$.
Proposition~\ref{prop:variation-sums} gives the exact formula
for derivatives of $J_{t,\beta}$. As illustrated in Step~3, we obtain
\[
 \sum_{j=1}^n \Delta_\tau J_{t,\beta}(f_{\tau}^j)
 \Big|_{\tau=0}\ge 0.
\]
The upper bound for this sum gives a bound on $W(f_{t,\beta})$
independent of $t$ for sufficiently small $\beta$.
Proposition~\ref{lem:vanishing-correction-minima} then allows us to
choose $t_j\to\infty$ and $\beta_j\downarrow0$ so that the
corresponding minimizers satisfy $q(f_{t_j,\beta_j})\to+\infty$ while $\sup_jA(f_{t_j,\beta_j})<\infty$. 
The optimal constant $2\pi(n+1)$ in Theorem~\ref{thm:pointwise-criterion} arises solely from the second variation of area functional, other variation terms
have no contributions. This finishes the proof.

\vspace{0.25cm}
\medskip
\noindent\textbf{Conventions and organization.}
We write $D_r=\{z\in\mathbb C:|z|<r\}$, $D_r(b)=\{z:|z-b|<r\}$, $S^1=\partial\D$. 
We denote by $h$ the Hermitian metric of $\omega$
and use the conventions
\[
 d^c=i(\bar\partial-\partial),\qquad
 dd^c=2i\partial\bar\partial,\qquad
 dA=\frac i2 dz\wedge d\bar z,\qquad
 \Delta_z=4\partial_z\partial_{\bar z}.
\]
Section~\ref{sec:stable-discs} states the compactness theorem we use and we prove it in Appendix~\ref{app:marked-compactness}.
Sections~\ref{sec:variation}--\ref{sec:exact-fields} establish
the disc functional criterion Theorem~\ref{thm:pointwise-criterion}. 
Section~\ref{sec:nonpsef} proves Theorem~\ref{thm:disc functional} and uniruled criterion Theorem~\ref{thm:intro-normal-kahler-uniruled}.
Section~\ref{sec:mixed-intersections}
proves K\"ahler Miyaoka--Mori Theorem~\ref{thm:intro-normal-kahler-rational-curves}.
Section~\ref{sec:hyperbolic-canonical-bigness} proves Kobayashi canonical ampleness conjecture.

\subsection*{Acknowledgements}
This research is supported by the National Key R\&D Program of China
(Grant Nos.~2021YFA1002600 and 2021YFA1003100). Z.~Wang and X.~Zhou
are partially supported by grants from the National Natural Science
Foundation of China (NSFC) (Nos.~12571085 and 12288201), respectively.
Z.~Wang is also supported by the Fundamental Research Funds for the
Central Universities. 

In developing the variational argument, we used AI to explore different choices of functional. AI suggested the term $te^{-q(f)}$.

\section{Compactness of holomorphic discs with finite area}
\label{sec:stable-discs}



Let $f_j:\D\to X$ be a sequence of holomorphic discs with uniformly bounded area $\int_\D f_j^*\omega$, the space of stable curves is an appropriate moduli space to describe the limit of $f_j$.
The compactness theorem originates from \cite{Gro85} and we adopt the language of stable curves developed in
\cite{ISInvent,ISCompactness}. Throughout this section we fix $x\in X$.

\begin{definition}\label{def:stable-graph-disc}
Let $Y$ be a complex manifold, a stable curve $C$ over $Y$ consists of
\begin{enumerate}
\item A connected nodal curve $C$: the point on $C$ is either smooth or a node: locally biholomorphic to $\{(u,v)\in\mathbb C^2:uv=0\}$.
Moreover, $C$ is obtained in the following manner:
beginning with the principal component $C_0\simeq\D$, we pick points $a_1,\ldots,a_N$ on $\D$ and glue a copy of $\mathbb P^1$ to each $a_i$ to make $a_i$ nodal. Then on each component $C_i\simeq\mathbb P^1$, we choose finitely many points on $C_i\backslash\, a_i$ and glue $\mathbb{P}^1$ at each of these points, $\ldots$. Such constructions of finitely many times produce the nodal curve $C$.
We say $a_1,\cdots,a_N$ are roots of $\mathbb{P}^1$ tree attached to principal component $C_0$.


\item A holomorphic map $F:C\to Y$ so that every spherical component where $F$ is constant has at
least three nodes. $f:=F|_{C_0}$ is called the principal map. 

\end{enumerate}
\end{definition}
Definition~\ref{def:stable-graph-disc}(1) makes $C$ can be parametrized by $\D$ in the sense of \cite[Definition 2.2]{ISCompactness}, see Definition~\ref{def:para curve} and the figure in Appendix~\ref{app:marked-compactness}. 

Definition~\ref{def:stable-graph-disc}(2) is equivalent to $F:C\to X$ being stable in the sense of \cite[Definition 2.4]{ISCompactness}.
The stability of $F:C\to X$ ensures every $\mathbb{P}^1$ tree of $C$ contains a $F$-non-constant
sphere leaf, since the end component only has one node so $F$ is non-constant there. 

\vspace{0.4cm}
We define the area of the disc curve $C$ by
\[
 \qquad A_*(U):=\int_C F^*\omega.
\]



The following is the so-called $\varepsilon$-regularity, one may consult e.g.
\cite[Lemma~1.1]{ISCompactness}, \cite[Proposition~4.1]{Zinger}. We give a simple proof for the reader's convenience.

\begin{lemma}
\label{lem:small-area-regularity}
There exist constants $e_0,C_{\rm reg}>0$, depending only on
$(X,\omega)$, such that every holomorphic map
$v:D_{2r}(b)\to X$ satisfies
\begin{equation}\label{eq:projective-epsilon-regularity}
 \int_{D_{2r}(b)}v^*\omega<e_0
 \quad\Longrightarrow\quad
 \sup_{D_r(b)}|dv|^2
 \le C_{\rm reg}r^{-2}
       \int_{D_{2r}(b)}v^*\omega.
\end{equation}
\end{lemma}

\begin{proof}
By rescaling $u(z):=v(b+rz)$, we may assume that $b=0$ and $r=1$. Since $X$ is compact, there are finitely many holomorphic coordinate
charts and a positive number
$\rho>0$ such that every ball $B_\omega(x,\rho)$ lies
in one of the charts.

Assume $a\in\overline{D_{3/2}}$ is a maximal point of
$(3/2-|z|)|dv(z)|$ in $\overline{D_{3/2}}$. Set $d=3/2-|a|$ and $M=|dv(a)|$.
We may assume that $d,M>0$, and an easy computation gives $|dv|\le 2M$ on $D_{d/2}(a)$.
Choose $s=\min\left\{\frac d2,\frac{\rho}{2M}\right\}$, we can easily see that
the image $v(D_s(a))$ lies in $B_\omega(v(a),\rho)$. Write $V$ for the coordinate representation of $v$
on $D_s(a)$. The submean-value inequality for $|V'|^2$ gives
\begin{equation}\label{eq:epsilon-reg}
     M^2s^2\le C\int_{D_s(a)}|V'|^2\,dA
 \le C'\int_{D_2}v^*\omega.
\end{equation}
Choose $e_0>0$ so that $C'e_0<\rho^2/4$. If $\int_{D_2}v^*\omega<e_0$,
then we deduce from \eqref{eq:epsilon-reg} that $s\ne\rho/(2M)$, so $s=d/2$ and $d^2M^2\le 4C'\int_{D_2}v^*\omega$.
For $z\in D_1$, the definition of $a$ yields
\[
 |dv(z)|^2
 \le \frac{d^2M^2}{(3/2-|z|)^2}
 \le 16C'\int_{D_2}v^*\omega,
\]
which completes the proof.
\end{proof}



Proposition~\ref{prop:interior-marked-compactness} below follows from
\cite[Theorem~1]{ISCompactness},
we include a detailed proof in Appendix~\ref{app:marked-compactness}. 
\begin{proposition}\label{prop:interior-marked-compactness}
Let $M>0$ and $f_j:\D\to X$ be holomorphic disc with $\int_\D f^*_j\omega\leq M$ and $f_j(0)=x$. 
Up to a subsequence, there is a stable
disc curve $(g,F):C\to \mathbb{C}\times X$ with the following properties:
\begin{enumerate}
\item $g|_{C_0}:C_0\to \D$ is biholomorphic and $g(C_i)$, $i\geq1$, are single points on $\D$. They consist of $g(\text{roots of $C$})$;
\item $f_j$ converge smoothly to $F\circ(g|_{C_0})^{-1}$ on every compact subset of $\D\backslash\,g(\text{roots of $C$})$;
\item there is $p\in C$ with $g(p)=0$ and $F(p)=x$;
\item for every smooth closed real $(1,1)$-form $\rho$ on $X$,
\begin{equation}\label{eq:stable-local-measure-convergence}
 \int_\D\chi\,f_j^*\rho\longrightarrow
 \int_{C}g^*\chi\,F^*\rho
                  \qquad \forall\, \chi\in C_c(\D);
\end{equation}
\item $A_*(U)\le\liminf_j A(f_j)$. 
\end{enumerate}
\end{proposition}

The following corollary illustrates that: if we can construct holomorphic discs $f_j$ at $x$ with uniform bounded area $A(f_j)$ and $q(f_j)\to+\infty$, then there is rational curve $R$ at $x$ and the degree of $R$ is bounded by $\liminf_j A(f_j)$.

\begin{corollary}
\label{cor:quantitative-marked-bubbling}
Fix $\delta,c>0$. 
Let $f_j:\D\to X$ be holomorphic discs with $f_j(0)=x$ and
$\sup_j A(f_j)<\infty$. If $q(f_j)\to+\infty$, there is a
rational curve $R$ through $x$ such that
\begin{equation}\label{eq:marked-bubble-area}
[\omega]\cdot R:= \int_R\omega\le\liminf_{j\to\infty}A(f_j).
\end{equation}
\end{corollary}

\begin{proof}
Let $(g,F):C\to \mathbb{C}\times X$ be the stable disc curve supplied by Proposition~\ref{prop:interior-marked-compactness}.
If $p$ does not lie on any $\PP^1$ tree, namely $p\in C_0\backslash\,\{\text{roots of }C\}$, then for some small $r>0$, $f_j$ converge smoothly on $D_r$, and in particular
$\sup_j\sup_{D_r}|df_j|<\infty$. 
Since $X$ is compact, $\theta-\delta\omega\le C_0\omega$ for some $C_0>0$. We see
\[
 q(f_j)\le |Q_\sigma(f_j)|
 \leq \frac{C_0}{2\pi}\int_{\D\backslash D_r}\big(\log\frac{1}{a}\big)f_j^*\omega+\frac{C_0\sup_{D_r}|df_j|^2}{2\pi}\int_{D_a}\log\frac{1}{|z|}\,dA,
\]
the first term is bounded by $\sup_jA(f_j)<+\infty$ and the second term is bounded by $\sup_j\sup_{D_r}|df_j|^2<+\infty$.
This contradicts with $q(f_j)\to+\infty$.

Therefore $p$ lies on some $\mathbb{P}^1$ tree. If $p$ belongs to a component where $F$ is non-constant, we already get a rational curve. If $p$ belongs to a constant component then we can follow the constant components to a nonconstant component $C_i\simeq \PP^1$, which exists because $C$ is stable as we explained before.
Since $F(p)=x$, we get $F$ also maps the connecting node of $C_i$ to $x$. Thus $R=F(C_i)$ is a rational curve through $x$.
Finally, the lower semi-continuity property in Proposition~\ref{prop:interior-marked-compactness} gives us $\int_R\omega\leq \int_{C_i}F^*\omega\leq F_*(C)\leq \liminf_{j\to\infty}A(f_j)$. 
The first inequality uses $\operatorname{deg}(C_i\xrightarrow{F} R)\cdot\int_R\omega=\int_{C_i}F^*\omega$.
The proof is complete.
\end{proof}

\vspace{0.25cm}
\section{Minimizers of functional \texorpdfstring{$J_{t,\beta}$}{1}}\label{sec:variation}
We abbreviate $\mathcal{A}(X,x)$ to $\mathcal A_x$.
Assume \eqref{eq:pointwise-escape}: $\sup_{f\in \mathcal A_x} q(f)=+\infty$, Corollary~\ref{cor:quantitative-marked-bubbling} reduces Theorem~\ref{thm:pointwise-criterion} to find special $f_j\in \mathcal A_x$ so that $\lim q(f_j)=+\infty$ and $\sup A(f_j)<+\infty$. It turns out that such holomorphic discs come from minimizer of the functional $J_{t,\beta}$ introduced below.

For $t>0$ let $J_t(f):=A(f)+te^{-q(f)}$, $f\in \mathcal A_x$ and $j_t:=\inf_{f\in\mathcal A_x}J_t(f)$.
For a minimizing sequence $f_j$ of $J_t$, their area $A(f_j)$ is uniformly bounded. However, the limit could be a stable disc with bubbles.
We have to introduce a functional $\Psi$ so that for $\beta>0$, the minimizer of functional $J_{t,\beta}:=A+te^{-q}-\beta\Psi$ are actual holomorphic discs, see Proposition~\ref{prop:green-compactness}.

Set
\begin{equation}\label{eq:kernel}
 G(z,w):=\log\left|\frac{z-w}{1-z\bar w}\right|^2,\qquad (z,w)\in \D^2.
\end{equation}
On
$\overline\D^2\setminus\D^2$ we set $G:=0$. It is straightforward to see that $G$ is nonpositive, symmetric, and
upper semicontinuous on $\overline\D^2$. 
For fixed $w\in\D$, $G(\cdot,w)$ is the Green function of $\D_z$ with pole at $w$:
\begin{equation}\label{eq:normalization}
 \partial_z\partial_{\bar z}G(z,w)=\pi\delta_w,
 \qquad G(\cdot,w)|_{\partial\D}=0.
\end{equation}

\begin{definition}\label{def:green-quadratic-form}
For a positive Radon measure $\nu$ on $\overline\D$, we define
\[
 \Psi(\nu):=\int_{\overline\D^2}G(z,w)
                         \dd\nu(z)\dd\nu(w)\in[-\infty,0].
\]
For a holomorphic disc $f:\D \to X$ with finite area $A(f)=\int_\D f^*\omega<\infty$, the measure $f^*\omega$ trivially extends to $\overline\D$ and we put $\Psi(f):=\Psi(f^*\omega)$. The measure $\nu$ we will use are only $f^*\omega$ and their limit.
\end{definition}
Let $\mathcal B_x$ be the set of holomorphic discs
$f:\D\to X$ such that
\[
 f(0)=x,\qquad A(f)<\infty,\qquad \Psi(f)>-\infty.
\]
Note that if $f\in\mathcal A_x$ then $f^*\omega$ is smooth in a neighborhood of $\overline\D$, so $A(f),\Psi(f)$ are finite. Therefore $\mathcal A_x\subset\mathcal B_x$.

For $t,\beta>0$ define
\begin{equation}\label{eq:variational-functional}
 J_{t,\beta}(f):=A(f)+t e^{-q(f)}-\beta\Psi(f),
 \qquad f\in\mathcal B_x,
\end{equation}
and
\begin{equation}\label{eq:variational-infimum}
 j_{t,\beta}:=\inf_{f\in\mathcal A_x}J_{t,\beta}(f).
\end{equation}
The goal of this section is to construct $f_{t,\beta}\in \mathcal B_x$ to achieve infimum of $J_{t,\beta}$.

\vspace{0.2cm}
\subsection{Radial approximation}
In this section, we prove that for $f\in\mathcal B_x$, $A(f),W(f),Q(f),\Psi(f)$ can be approximated by its radial scaling $f_r$.
For a holomorphic disc $f:\D\to X$ and $0<r<1$,
write
\[
 f_r(\zeta):=f(r\zeta).
\]
Then $f_r$ is defined in a neighborhood of $\overline\D$.
A changing variables gives
\begin{align}
 A(f_r)&=\int_{D_r}f^*\omega,
                                      \label{eq:radial-area}\\
 W(f_r)&=\int_{D_r}(1-|z|^2/r^2)\,f^*\omega,
                                      \label{eq:radial-weight}\\
 Q_\rho(f_r)&=\frac{1}{2\pi}\int_{D_r}
              \log(r/|z|)\,f^*\rho
                                      \label{eq:radial-green}
\end{align}
for every smooth closed real $(1,1)$-form $\rho$ on $X$.
Set
\[
 \kappa_r:=\sup_{z\in\D}
 \frac{\min\big\{\log(1/|z|),\log(1/r)\big\}}{1-|z|^2}<\infty,
\]
and, for $M\ge0$, put
\[
 D_M(r):=
 M\left\{\frac{C_\sigma}{2\pi}\log\frac1r+c(1-r^2)\right\},
\]
where $\sigma:=\theta-\delta\omega$ and $C_\sigma>0$ so that $-C_\sigma\omega\leq\sigma\leq C_\sigma\omega$ on $X$.
\begin{lemma}\label{lem:radial-ordinary}
Let $f:\D\to X$ be a holomorphic disc with finite area.
\begin{enumerate}
\item As $r\uparrow1$, one has
\[
 A(f_r)\to A(f),\qquad W(f_r)\to W(f),\qquad
 Q_\rho(f_r)\to Q_\rho(f)
\]
for every smooth closed real $(1,1)$-form $\rho$ on $X$.
If $A(f)\le M$, then
\begin{equation}\label{eq:radial-q-error}
 |q(f)-q(f_r)|\le D_M(r).
\end{equation}

\item For every $0<r<1$, one has 
\begin{align}
 A(f_r)&\le\frac{W(f)}{1-r^2},
                                      \label{eq:weighted-restriction-area}\\
 |q(f)-q(f_r)|
 &\le\left(\frac{C_\sigma\kappa_r}{2\pi}+c\right)W(f).
                                      \label{eq:weighted-restriction-q}
\end{align}
\end{enumerate}
\end{lemma}

\begin{proof}
From the definition of $W$ we get
\begin{equation}\label{eq:radial-W-error}
 0\le W(f)-W(f_r)
 \le\min\{A(f)(1-r^2),W(f)\}.
\end{equation}
Similarly, \eqref{eq:radial-green} gives
\[
 Q_\rho(f)-Q_\rho(f_r)
 =\frac{1}{2\pi}\int_\D \min\{\log(1/|z|),\log(1/r)\} f^*\rho.
\]
Recall that $C_\rho$ is the constant such that $-C_\rho\omega\leq\rho\leq C_\rho\omega$ on $X$, so we get 
\begin{equation}\label{eq:radial-Q-error}
 |Q_\rho(f)-Q_\rho(f_r)|
 \le\frac{C_\rho A(f)}{2\pi}\log(1/r).
\end{equation}
These estimates, together with the trivial monotone convergence $\lim_{r\uparrow 1} A(f_r)=A(f)$, prove the convergences in~(1).
Taking $\rho=\sigma$ and using $q:=Q_\sigma-cW$ gives
\eqref{eq:radial-q-error}.

For~(2),
\[
 (1-r^2)A(f_r)\le\int_{D_r}(1-|z|^2) f^*\omega\le W(f),
\]
which gives \eqref{eq:weighted-restriction-area}.
\[
 |Q_\sigma(f)-Q_\sigma(f_r)|
 =\Big|\frac{1}{2\pi}\int_\D\min\{\log(1/|z|),\log(1/r)\}f^*\sigma\Big|
 \leq\frac{C_\sigma \kappa_r}{2\pi}\int_\D(1-|z|^2)f^*\omega
 =\frac{C_\sigma\kappa_r}{\pi}W(f).
\]
Together with \eqref{eq:radial-W-error}, this proves
\eqref{eq:weighted-restriction-q}.

\end{proof}

\begin{lemma}
\label{lem:green-radial}
Let $f:\D\to X$ be holomorphic with $A(f)<\infty$.
Then
\begin{equation}\label{eq:green-radial-monotonicity}
 \Psi(f_r)>-\infty,\qquad
 \Psi(f_r)\downarrow\Psi(f)\quad(r\uparrow1).
\end{equation}
$\Psi(f)$ may be $-\infty$.
\end{lemma}

\begin{proof}
For $0<r\le1$,
we define the scaling Green function on $\overline\D^2$ by
\[
 G_r(z,w):=
 \begin{cases}
  \log\left|\frac{r(z-w)}{r^2-z\bar w}\right|^2,&z,w\in D_r,\ z\ne w,\\
  -\infty,&z=w\in D_r,\\
  0,&(z,w)\notin D_r^2.
 \end{cases}
\]
Then $G_r(z,w)=G(r^{-1}z,r^{-1}w)$ for $(z,w)\in D_r^2$.
Changing variables by $z=r\zeta$ and
$w=r\eta$ gives
\begin{equation}\label{eq:green-radial-kernel}
 \Psi(f_r)=
 \int_{\overline\D^2}G_r(z,w)\,d\mu(z)d\mu(w).
\end{equation}
Since $\mu=f^*\omega$ is smooth on $\D$, we see $\Psi(f_r)>-\infty$.

For fixed $z\ne w$ and $r>\max(|z|,|w|)$, we compute that
\[
 \frac{\partial G_r(z,w)}{\partial r}
 =-\frac2r\,
   \frac{r^4-|z\bar w|^2}{|r^2-z\bar w|^2}<0.
\]
Consequently, if $0<r<s<1$, then $G_r\ge G_s$ on $D_r^2$.
Outside $D_r^2$ clearly one has $G_r=0\ge G_s$. Thus
\begin{equation}\label{eq:extended-Green-monotonicity}
 0\ge G_r\ge G_s\ge G\qquad(0<r<s<1).
\end{equation}

Since $G_r\downarrow G$ pointwise on
$\overline\D^2$.
Monotone convergence theorem now yields
\[
 \Psi(f_r)
 =\int_{\overline\D^2}G_r \,d(\mu\otimes\mu)
 \downarrow
 \int_{\overline\D^2}G\,d(\mu\otimes\mu)
 =\Psi(f).
\]
The proof is complete.
\end{proof}

\begin{corollary}
\label{cor:ordinary-competitors}
For every $f\in\mathcal B_x$, the discs $f_r\in\mathcal A_x$ and satisfy
\[
 J_{t,\beta}(f_r)\longrightarrow J_{t,\beta}(f)
 \qquad(r\uparrow1).
\]
In particular
\begin{equation}\label{eq:ordinary-infimum}
 j_{t,\beta}
 =\inf_{f\in\mathcal B_x}J_{t,\beta}(f).
\end{equation}
If \eqref{eq:pointwise-escape} holds, then for every fixed
$\beta>0$,
\begin{equation}\label{eq:sublinear-infimum}
 \frac{j_{t,\beta}}t\longrightarrow0\qquad(t\to\infty).
\end{equation}
\end{corollary}

\begin{proof}
The map $f_r$ is holomorphic on $D_{1/r}$ and is centered at $x$.
Apply Lemma~\ref{lem:radial-ordinary} and
Lemma~\ref{lem:green-radial} to each term in
\eqref{eq:variational-functional}, we get $J_{t,\beta}(f_r)\to J_{t,\beta}(f)$.
\eqref{eq:ordinary-infimum} follows from this approximation
and $\mathcal A_x\subset\mathcal B_x$.

For every fixed $g\in\mathcal A_x$, $\Psi(g)>-\infty$ and
\[
 0\le\frac{j_{t,\beta}}t
 \le\frac{A(g)-\beta\Psi(g)}t+e^{-q(g)}.
\]
First let $t\to\infty$ and then use
\eqref{eq:pointwise-escape} to choose $g$ with arbitrarily
large $q(g)$, the proof is complete.
\end{proof}

\begin{lemma}
\label{lem:almost-minimizing-collar}
Fix $t,\beta>0$, $M<\infty$, and $\eta\ge0$.
Suppose $f\in\mathcal B_x$ satisfies
\[
 A(f)\le M,\qquad
 J_{t,\beta}(f)\le j_{t,\beta}+\eta.
\]
With $D_M(r)$ defined in
\eqref{eq:radial-q-error}, one has
\begin{equation}\label{eq:foundation-collar}
 \int_{\D\setminus D_r}f^*\omega
 =A(f)-A(f_r)
 \le\eta+te^{-q(f)}\bigl(e^{D_M(r)}-1\bigr)
 \qquad(0<r<1).
\end{equation}
\end{lemma}

\begin{proof}
Lemma~\ref{lem:radial-ordinary} gives
$|q(f)-q(f_r)|\le D_M(r)$.
By Corollary~\ref{cor:ordinary-competitors},
$J_{t,\beta}(f_r)\ge j_{t,\beta}$.
Together with $J_{t,\beta}(f)\le j_{t,\beta}+\eta$, this gives
\[
\begin{aligned}
 A(f)-A(f_r)+\beta\bigl(\Psi(f_r)-\Psi(f)\bigr)
 &\le\eta+te^{-q(f_r)}-te^{-q(f)}\\
 &=\eta+te^{-q(f)}\bigl(e^{q(f)-q(f_r)}-1\bigr)\\
 &\le\eta+te^{-q(f)}\bigl(e^{D_M(r)}-1\bigr).
\end{aligned}
\]
Lemma~\ref{lem:green-radial} allows us to discard the
nonnegative difference $\Psi(f_r)-\Psi(f)$ on the left.
This proves \eqref{eq:foundation-collar}.
\end{proof}

\subsection{Obtain minimizer of $J_{t,\beta}$}

The minimizer holomorphic disc turns out to be $W^{1,p}$, $2<p<4$. We give the precise definition.
\begin{definition}
    A holomorphic disc $f:\D\to X$ is called $W^{1,p}$ if $\int_\D|df|^p_hdA_z<\infty$ and we will say $f\in W^{1,p}(\D,X)$.
\end{definition}
Let $g$ be the induced Riemannian metric of $h$ and $\iota:X\to\mathbb R^N$ be an isometry embedding of $(X,g)$.
Then $\iota\circ f$ corresponds to smooth functions $f^1,\cdots,f^N$ and $|df|^2_h$ coincides with $\sum_{i=1}^N|\frac{\partial f^i}{\partial z}|^2+|\frac{\partial f^i}{\partial\bar z}|^2$.
Therefore $f\in W^{1,p}(\D,X)$ is equivalent to $f^i\in W^{1,p}(\D)$, $i=1,\cdots,N$.
\begin{definition}
    Let $f:\D\to X$ be a $W^{1,p}$ holomorphic disc, $\xi$ be a holomorphic section of $f^*T^{1,0}X$.
    Using $T^{1,0}X\subset TX\otimes\mathbb C\hookrightarrow T\mathbb R^N\otimes\mathbb C\cong \mathbb R^N\times \mathbb{R}^{2N}$.
    $\xi$ corresponds to smooth functions $f^1,\cdots,f^N,\xi^1,\cdots,\xi^{2N}$.
    $\xi$ is called $W^{1,p}$ if every $\xi^i\in W^{1,p}(\D)$ and we will say $\xi\in W^{1,p}(\D,f^*T^{1,0}X)$.
\end{definition}
We remark that if we cover $X$ by finitely many coordinate charts $U_{\alpha}\subset\subset V_{\alpha}$, then $f\in W^{1,p}(\D,X)$ if and only if $f_{\alpha}^i:=f^*z_{\alpha}^i\in W^{1,p}(f^{-1}(U_\alpha))$ for every $\alpha,i$ where $z_{\alpha}^i$ are coordinates on $V_{\alpha}$.
The tangent bundle $T^{1,0}V_{\alpha}\cong V_\alpha\times \mathbb{C}^n$, and $\xi\in W^{1,p}(\D,f^*T^{1,0}X)$ if and only if $\xi^j_{\alpha}:=\xi^*v_{\alpha}^j\in W^{1,p}(f^{-1}(U_\alpha))$ for every $\alpha,j$, where $v_{\alpha}^j$ are coordinates of fiber $\mathbb{C}^n$.

\begin{lemma}
\label{lem:usc}
The following statements hold.
\begin{enumerate}
\item If $f\in W^{1,p}(\D,X)$ is holomorphic and $p>2$, then
$\Psi(f)>-\infty$.
\item Let $\nu$ be a positive Radon measure on $\overline\D$, if $\nu$ has a positive atom at a point of $\D$, then
$\Psi(\nu)=-\infty$.
\item If $\nu_j\rightharpoonup\nu$ weakly as positive Radon
measures on $\overline\D$, then
$\limsup_j\Psi(\nu_j)\le\Psi(\nu)$.
\end{enumerate}
\end{lemma}
\begin{proof}
It is clear that
\[
 0\le-G(z,w)\le \log\frac4{|z-w|^2}.
\]
Put $q:=p/2>1$ and $q':=q/(q-1)$. Writing
$f^*\omega=|df|^2\,dA$, we have $|df|^2\in L^q(\D)$.
H\"older's inequality gives, for every $w\in \D$,
\[
 \int_\D-G(z,w)\,e_f(z)\,dA_z
 \le \|G(\cdot,w)\|_{L^{q'}}\|e_f\|_{L^q}
 \le C_q\|e_f\|_{L^q},
\]
where we use the inequality
\[
\sup_{w\in\D}\int_\D |G(z,w)|^{q'}\,dA_z
 \le 2\pi\int_0^2
       \left(\log\frac4{t^2}\right)^{q'} t\,dt<\infty.
\]
Integrating against $|df|^2(w)\,dA_w$, we obtain
\[
 0\le-\Psi(f)
 \le C_q\|e_f\|_{L^q}\|e_f\|_{L^1}
 \le C_q'\|e_f\|_{L^q}^2<\infty.
\]

If $\nu(\{a\})=m>0$ then $\nu\otimes\nu$ has mass $m^2$ at
$(a,a)$. $G(a,a)=-\infty$ gives $\Psi(\nu)=-\infty$.

For the last statement, the product measures $\nu_j\otimes\nu_j$ converge weakly to
$\nu\otimes\nu$ on $\overline\D^2$. 
Since $G$ is upper semicontinuous on $\overline\D^2$, it is well-known that
\[
 \limsup_j\int_{\overline\D^2}G\,d(\nu_j\otimes\nu_j)
 \le\int_{\overline\D^2}G\,d(\nu\otimes\nu).
\]
This proves $ \limsup_j \Psi(\nu_j)\leq\Psi(\nu)$.
\end{proof}

\begin{proposition}
\label{prop:green-compactness}
Fix $M,L<\infty$.
Let $f_j:\D\to X$ be holomorphic, with $f_j(0)=x$, and suppose
\[
 A(f_j)\le M,\qquad \Psi(f_j)\ge-L.
\]
There are a subsequence of $f_j$ (still denoted by $f_j$), a holomorphic map $f:\D\to X$ with
$f(0)=x$, and a finite positive measure $\nu_\partial$
supported on $S^1$ such that:
\begin{enumerate}
\item $f_j\to f$ in every $C^k$ norm on every compact subset
of $\D$;
\item 
$
 f_j^*\omega\rightharpoonup
       f^*\omega+\nu_\partial$
       on  $\overline\D$;
\item $A(f)\le\liminf_jA(f_j)$ and
$\Psi(f)\ge\limsup_j\Psi(f_j)\ge-L$;
\item for every smooth closed real $(1,1)$-form $\rho$ on $X$,
\[
 \frac{1}{2\pi}\int_\D\log(1/|z|)f_j^*\rho\longrightarrow \frac{1}{2\pi}\int_\D \log(1/|z|)f^*\rho,\qquad
 \int_\D (1-|z|^2) f_j^*\rho\longrightarrow\int_\D (1-|z|^2) f^*\rho.
\]
In particular $W(f_j)\to W(f)$ and $q(f_j)\to q(f)$.
\end{enumerate}
\end{proposition}

\begin{proof}
Write $\nu_j:=f_j^*\omega$, then $\nu_j(\overline\D)\le M$.
After passing to a subsequence, there is a finite positive
measure $\nu$ on $\overline\D$ such that
\begin{equation}\label{eq:Green-compactness-measure-limit}
 \nu_j\rightharpoonup\nu,\qquad
 \nu(\overline\D)=\lim_j\nu_j(\overline\D)\le M.
\end{equation}
By Lemma~\ref{lem:usc}(3),
\[
 0\ge\Psi(\nu)\ge\limsup_j\Psi(\nu_j)\ge-L.
\]
Lemma~\ref{lem:usc}(2) consequently gives
\begin{equation}\label{eq:Green-limit-no-atoms}
 \nu(\{a\})=0,\qquad\forall a\in\D.
\end{equation}
There exists $r=r(a)>0$ such that
\[
 \overline{D_{4r}(a)}\subset\D,\qquad
 \nu\bigl(\overline{D_{4r}(a)}\bigr)<e_0/2,
\]
where $e_0$ is the constant in $\varepsilon$-regularity
Lemma~\ref{lem:small-area-regularity}.

The weak convergence $\nu_j\rightharpoonup\nu$ implies for $j$ large $\int_{D_{2r}(a)}f^*_j\omega<e_0/2$,
Lemma~\ref{lem:small-area-regularity} then yields
\begin{equation}\label{eq:Green-local-derivative-bound}
 \sup_{D_{r}(a)}|df_j|^2
 \le C_{\rm reg}r^{-2}
          \int_{D_{2r}(a)}f_j^*\omega
 \le C_{\rm reg}r^{-2}e_0.
\end{equation}

Since $X$ is compact, after a diagonal subsequence, Montel theorem
\[
 f_j\longrightarrow f
 \quad\text{in every $C^k$ norm on compact subsets of }\D
\]
for a holomorphic map $f:\D\to X$. Since $f_j(0)=x$,
one has $f(0)=x$.

For every $\chi\in C_c(\D)$, statement~(1) and
\eqref{eq:Green-compactness-measure-limit} imply
\[
 \int_{D_r}f^*\omega
 =\lim_j\int_{D_r}f_j^*\omega
 \le\liminf_j A(f_j)\leq M,\qquad
 \int\chi\,d\nu
 =\lim_j\int_\D\chi f_j^*\omega
 =\int_\D\chi f^*\omega.
\]
Letting $r\uparrow1$ yields $A(f)=\int_\D f^*\omega\leq \liminf_j A(f_j)\leq M$. Set $\nu_\partial:=\nu|_{S^1}$ we obtain $\nu=f^*\omega+\nu_\partial$,
which is statement~(2).

By definition \eqref{eq:kernel},
$G=0$ on
$(S^1\times\overline\D)\cup(\overline\D\times S^1)$.
Write $\mu=f^*\omega$ and use $\nu=\mu+\nu_\partial$, we get
\[
 \begin{aligned}
 \Psi(\nu)
 &=\int_{\overline\D^2}G\,
       d\bigl((\mu+\nu_\partial)\otimes(\mu+\nu_\partial)\bigr)\\
 &=\int_{\overline\D^2}G\,d(\mu\otimes\mu)
 =\Psi(f).
 \end{aligned}
\]
Consequently Lemma~\ref{lem:usc}(3) gives
\begin{equation}\label{eq:Green-quadratic-upper-semicontinuity}
 \Psi(f)=\Psi(\nu)
 \ge\limsup_j\Psi(f_j)\ge-L,
\end{equation}
which proves statement~(3).

Finally, fix a smooth closed real $(1,1)$-form $\rho$ on $X$.
By the convergence $f_j\to f$ and $\rho\leq C_\rho\omega$
\[
 \limsup_j\left|
 \int_\D (1-|z|^2)(f_j^*\rho-f^*\rho)\right|
 \le 2C_\rho M(1-r^2),
 \qquad\forall\, r<1.
\]
To deal with $\log(1/|z|)$, we choose $\varepsilon>0$. 
On $\{\varepsilon\leq|z|<1\}$, using convergence $f_j\to f$ we get 
\[
 \limsup_j\left|
 \int_{\{\varepsilon\leq|z|<1\}} \log(1/|z|)(f_j^*\rho-f^*\rho)\right|
 \le 2C_\rho M\log(1/r)\longrightarrow0
\]
as before. On small ball $\{|z|<\varepsilon\}$, we have $|f^*_j\rho|+|f^*\rho|\leq C\,dA$ for some uniform $C>0$. Dominated convergence theorem therefore completes the proof.
\end{proof}

We now establish the main result of this section.
\begin{proposition}
\label{prop:foundation-attainment}
For every $t,\beta>0$, there is $f_{t,\beta}\in\mathcal B_x$
such that $J_{t,\beta}(f_{t,\beta})=j_{t,\beta}$.
Every such minimizer satisfies
\begin{equation}\label{eq:minimum-integrated-ricci}
 I(f_{t,\beta})\ge\delta A(f_{t,\beta}),\qquad
 I_w(f_{t,\beta})\ge\delta W(f_{t,\beta}),
\end{equation}
where $I(f):=\int_\D f^*\theta$ and $I_w(f):=\int_\D (1-|z|^2) f^*\theta$. We alao have
\begin{equation}\label{eq:minimum-q-radial}
 q(f_{t,\beta})\ge q((f_{t,\beta})_r),\qquad
 Q_\sigma(f_{t,\beta})\ge
                 Q_\sigma((f_{t,\beta})_r)
 \quad(0<r\le1).
\end{equation}
The map $f_{t,\beta}$ belongs to $W^{1,p}(\D,X)$ for every $2<p<4$.
\end{proposition}

\begin{proof}
Considering the constant disc $f(z)=x$, we have $J_{t,\beta}(f)=t$ and hence $j_{t,\beta}\le t$.
Choose a sequence $f_j\in\mathcal A_x$ with
$J_{t,\beta}(f_j)\to j_{t,\beta}$.
We may assume without loss of generality that $J_{t,\beta}(f_j)\le t+1$, and get
\[
 A(f_j)\le t+1,\qquad
 0\le-\beta\Psi(f_j)\le t+1.
\]
Thanks to Proposition~\ref{prop:green-compactness}, up to a subsequence, there is a holomorphic map $f:\D\to X$ with
\[
 \begin{gathered}
 f(0)=x,\qquad
 A(f)\le\liminf_j A(f_j)\le t+1,\\
 \Psi(f)\ge\limsup_j\Psi(f_j)\ge-(t+1)/\beta,\qquad
 q(f_j)\longrightarrow q(f).
 \end{gathered}
\]
In particular, $f\in\mathcal B_x$. Now, we can write
\[
 \begin{aligned}
 J_{t,\beta}(f)
 &=A(f)-\beta\Psi(f)+te^{-q(f)}\\
 &\le\liminf_j A(f_j)
       +\liminf_j\bigl(-\beta\Psi(f_j)\bigr)
       +te^{-q(f)}\\
 &\le\liminf_j\bigl(A(f_j)-\beta\Psi(f_j)\bigr)
       +te^{-q(f)}\\
 &=\liminf_j J_{t,\beta}(f_j)
 =j_{t,\beta}.
 \end{aligned}
\]
By Corollary~\ref{cor:ordinary-competitors}, we obtain $J_{t,\beta}(f)=j_{t,\beta}$.

We next prove \eqref{eq:minimum-q-radial}.
By Corollary~\ref{cor:ordinary-competitors}, we can select a minimizer $f\in\mathcal B_x$ of $J_{t,\beta}$. For $0<r\le1$, set $f_r(z):=f(rz)$ and
$\mu:=te^{-q(f)}$. We clearly have
$J_{t,\beta}(f_r)\ge J_{t,\beta}(f)$, and thus
\begin{equation}\label{eq:minimum-radial-comparison}
 t\bigl(e^{-q(f_r)}-e^{-q(f)}\bigr)\geq  A(f)-A(f_r)+\beta\bigl(\Psi(f_r)-\Psi(f)\bigr).
\end{equation}
The right-hand side is nonnegative by
Lemma~\ref{lem:radial-ordinary} and
Lemma~\ref{lem:green-radial}, hence $q(f)\ge q(f_r)$.
Also $W(f_r)\le W(f)$ by \eqref{eq:radial-weight}.
Consequently,
\begin{equation}\label{eq:radial-green-maximum}
 Q_\sigma(f)-Q_\sigma(f_r)
 =q(f)-q(f_r)+c\bigl(W(f)-W(f_r)\bigr)\ge0.
\end{equation}
This proves \eqref{eq:minimum-q-radial}.

We next establish \eqref{eq:minimum-integrated-ricci}.
Set $\nu:=f^*\sigma=f^*\theta-\delta f^*\omega$, it is sufficient to prove that 
\[\nu(\D)\ge0,\qquad \int_\D (1-|z|^2)\,d\nu\ge0.\] 
It is easy to see that $\log(r/|z|)$ is integrable with respect to $|\nu|$. We then use Fubini's theorem together with
\eqref{eq:radial-green} to write
\begin{equation}\label{eq:radial-green-fubini}
 Q_\sigma(f_r)
 =\frac{1}{2\pi}\int_{D_r}\int_{|z|}^r\frac{ds}{s}\,d\nu(z)
 =\frac{1}{2\pi}\int_0^r\frac{\nu(D_s)}{s}\,ds.
\end{equation}
Since $\nu$ has finite total variation,
$\nu(D_s)/s\to\nu(\D)$ as $s\uparrow1$. Consequently,
\[
 \nu(\D)
 =\lim_{r\uparrow1}\frac1{1-r}
       \int_r^1\frac{\nu(D_s)}s\,ds
 =2\pi\lim_{r\uparrow1}\frac{Q_\sigma(f)-Q_\sigma(f_r)}{1-r}
 \ge0.
\]
Similarly, using again Fubini's theorem and noting that $w(z)=\int_{|z|}^1 2s\,ds$, we have
\[
 \begin{aligned}
 \int_\D w\,d\nu
 &=\int_\D\left(\int_{|z|}^1 2s\,ds\right)d\nu(z)\\
 &=2\int_0^1 s\,\nu(D_s)\,ds\\
 &=4\int_0^1 r
       \left(\int_r^1\frac{\nu(D_s)}{s}\,ds\right)dr\\
 &=8\pi\int_0^1 r\bigl(Q_\sigma(f)-Q_\sigma(f_r)\bigr)\,dr
 \ge0.
 \end{aligned}
\]
This proves \eqref{eq:minimum-integrated-ricci}.

Finally we prove $f\in  W^{1,p}(\D,X)$ for $2<p<4$.
If $\frac{1}{2}<r<1$, \eqref{eq:foundation-collar} gives
\[
 \int_{\D\setminus D_r}f^*\omega
 =A(f)-A(f_r)\le te^{-q(f)}\Big(e^{D_M(r)}-1\Big)\leq C(1-r)
\]
where the constant $C$ is independent of $r$.
If $|z|\geq\frac{1}{3}$, then
$ D_{(1-|z|)/2}(z)
 \subset\D\setminus D_{(3|z|-1)/2}$.
Thus for $|z|\geq r_0:=1-\frac{2e_0}{3C}$, it holds that
\[
 \int_{D_{(1-|z|)/2}(z)}f^*\omega
 \leq \int_{\D\backslash D_{(3|z|-1)/2}}f^*\omega
 \leq \frac{3C}{2}(1-|z|)
 \leq e_0
\]
where $e_0$ is the constant in Lemma~\ref{lem:small-area-regularity}, from which we get
\[
 |df(z)|^2
 \le \frac{16C_{\rm reg}}{(1-|z|)^2}
       \int_{D_{(1-|z|)/2}(z)}f^*\omega
 \le\frac{C'}{1-|z|}.
\]
For sufficiently large integer $k$, set
\[
 E_k:=\{z:2^{-k-1}<1-|z|\le2^{-k}\}.
\]
The derivative bound of $f$ and \eqref{eq:foundation-collar} give,
respectively,
\[
 \sup_{E_k}|df|\le \sup_{E_k}\sqrt{\frac{C'}{1-|z|}} =C'^{\frac{1}{2}}2^{\frac{k+1}{2}},
 \qquad
 \int_{E_k}f^*\omega\leq\int_{\D\backslash D_{1-2^{-k-1}}}f^*\omega  \le C2^{-k-1}.
\]
For $2<p<4$, it follows that
\[
 \begin{aligned}
 \int_{E_k}|df|^p\,dA
 &\le \bigl(\sup_{E_k}|df|\bigr)^{p-2}
           \int_{E_k}|df|^2dA\\
 &\le C''2^{\frac{k(p-2)}{2}}2^{-k}
 =C2^{-k(2-\frac{p}{2})}.
 \end{aligned}
\]
This is summable when $p<4$. This proves $f\in W^{1,p}(\D,X)$.
\end{proof}

\vspace{0.25cm}
\section{Family of holomorphic discs of class \texorpdfstring{$W^{1,p}$}{2} and their variations}\label{sec:green-levi}

Let $F(\tau,z):D_{\varepsilon}\times \D\to X$ be a holomorphic map, then $\{f_\tau:=F(\tau,\cdot)\}_{\tau\in D_{\varepsilon}}$ is called a family of holomorphic discs.

The differential $dF$ is a holomorphic section of $T^{*1,0}(D_{\varepsilon}\times\D)\otimes F^*T^{1,0}X$ and 
$
\partial_\tau F:=dF(\frac{\partial}{\partial\tau}),\partial_zF:=dF(\frac{\partial}{\partial z})
$ are holomorphic sections of $F^*T^{1,0}X$.
The following commutativity is very useful.

\begin{lemma}\label{lem:closed}
Let $h$ be the Hermitian metric induced by K\"ahler form $\omega$, then it holds that
\begin{equation}
\begin{aligned}
&\partial_\tau h(\partial_zF,\partial_zF)=\partial_zh(\partial_\tau F,\partial_zF),\qquad
\partial_{\bar\tau}h(\partial_\tau F,\partial_zF)=\partial_{\bar z}h(\partial_\tau F,\partial_\tau F),\\
&\partial_\tau\partial_{\bar\tau}h(\partial_zF,\partial_zF)=
\partial_z\partial_{\bar z}h(\partial_\tau F,\partial_\tau F).
\end{aligned}
\end{equation}
More generally, let $h^\rho$ be the Hermitian form induced by a closed real $(1,1)$-form $\rho$ on $X$, then the same equation holds:
\begin{equation}\label{eq:actual-local-variations}
\partial_\tau h^\rho(\partial_zF,\partial_zF)
 =\partial_z h^\rho(\partial_\tau F,\partial_zF),
 \qquad
 \partial_\tau\partial_{\bar\tau}h^\rho(\partial_zF,\partial_zF)
 =\partial_z\partial_{\bar z}h^\rho(\partial_\tau F,\partial_\tau F).
\end{equation}
\end{lemma}

\begin{proof}
The identities basically follows from $d\omega=0$: in local coordinates write $\omega=\frac{i}{2}(g_{i\bar j}dz^i\wedge d\bar{z}^j)$, then $\partial_kg_{i\bar j}=\partial_i g_{k\bar j}$, $\partial_{\bar k}g_{i\bar j}=\partial_{\bar j}g_{i\bar k}$ and $\partial_k\partial_{\bar l}g_{i\bar j}=\partial_i\partial_{\bar j}g_{k\bar l}$.

A direct calculation shows that
\[
 \begin{aligned}
 F^*\omega=\frac i2\bigl(&
h(\partial_zF,\partial_zF)\,dz\wedge d\bar z
 +h(\partial_\tau F,\partial_zF)\,d\tau\wedge d\bar z\\
 &+\overline{h(\partial_\tau F,\partial_zF)}\,dz\wedge d\bar\tau
 +h(\partial_\tau F,\partial_\tau F)\,d\tau\wedge d\bar\tau\bigr).
 \end{aligned}
\]
Since $d\omega=0$, we have $\partial(F^*\omega)=0$, which gives 
\[
 \begin{aligned}
 0=\frac i2\,d\tau\wedge dz\wedge\bigl[
 &(\partial_\tau h(\partial_zF,\partial_zF)-\partial_z h(\partial_\tau F,\partial_zF))\,d\bar z\\
 &+(\partial_\tau\overline{h(\partial_\tau F,\partial_zF)}-\partial_z h(\partial_\tau F,\partial_\tau F))
       \,d\bar\tau\bigr].
 \end{aligned}
\]
We then infer that $\partial_\tau\partial_{\bar\tau}h(\partial_zF,\partial_zF)=\partial_\tau\partial_{\bar z}h(\partial_zF,\partial_\tau F)=\partial_z\partial_{\bar z}h(\partial_\tau F,\partial_\tau F)$. The proof is complete.
\end{proof}

\begin{definition}
    Let $p>2$, $F:D_\varepsilon\times\D\to X$ be a holomorphic map, $\iota:X\to \mathbb{R}^N$ be an embedding.
    $F$ corresponds to smooth functions $F^i,\cdots,F^N$ on $D_\varepsilon\times\D$. 
    We say $F\in C^2(D_\varepsilon,W^{1,p}(\D,X))$ if
    $F^i\in C^2(D_\varepsilon,W^{1,p}(\D))$, $i=1,\cdots,N$.
\end{definition}
    Since $p>2$, the Sobolev embedding $W^{1,p}(\D)\hookrightarrow C^{0,1-\frac{2}{p}}(\D)$ and the multiplication map $W^{1,p}(\D)\times W^{1,p}(\D)\to W^{1,p}(\D)$ is continuous. 
    For another embedding $\iota':X\to \mathbb{R}^{N^\prime}$, there is a smooth mapping $\Phi:\mathbb{R}^N\to \mathbb{R}^{N'}$ so that $\Phi\circ\iota=\iota'$.
    Using chain rule we can check that $F^i\in C^2(D_\varepsilon;W^{1,p}(\D))$, $i=1,\cdots,N$ implies $(F')^j\in C^2(D_\varepsilon;W^{1,p}(\D))$, $j=1,\cdots,N'$. Therefore $F\in C^2(D_\varepsilon;W^{1,p}(\D,X))$ is independent of the choice of the embedding.
    
\begin{remark}\label{rem:C2 family}
    $F^i\in C^2(D_\varepsilon,W^{1,p}(\D))$ and $F^i$ is smooth on $D_\varepsilon\times\D$ imply that, writing $\tau=x+iy$, we have
    $\partial_xF^i,\partial_yF^i\in C(D_{\varepsilon},W^{1,p}(\D))$ and
    \[
    (|x|+|y|)^{-1}\|F^i(\tau+x+iy,\cdot)-F^i(\tau,\cdot)-(x\partial_xF^i(\tau,\cdot)+y\partial_yF^i(\tau,\cdot))\|_{W^{1,p}(\D)}\to 0,\qquad(x,y)\to0
    \]
    Similar differentiability holds for $\partial_x F^i$ (and $\partial_yF^i$):  $\partial_x\partial_xF^i,\partial_y\partial_xF^i\in C(D_\varepsilon,W^{1,p}(\D))$ so that
    \[
    (|x|+|y|)^{-1}\|\partial_xF^i(\tau+x+iy,\cdot)-\partial_xF^i(\tau,\cdot)-(x\partial_x\partial_xF^i(\tau,\cdot)+y\partial_y\partial_xF^i(\tau,\cdot))\|_{W^{1,p}(\D)}\to 0,\qquad(x,y)\to0.
    \]
    If $\iota:X\to \mathbb R^N$ is an isometry embedding of $(X,h)$, then the norm of the differential of $F$ coincides with that of $\iota\circ F$, so we have
    \begin{equation}\label{eq:embedding}    |d_zF|^2_h:=h(\partial_zF,\partial_zF)=\sum_{i=1}^N\big(\big|\frac{\partial F^i}{\partial z}\big|^2+\big|\frac{\partial F^i}{\partial\bar{z}}\big|^2\big).
    \end{equation}
    Since $F\in C^2(D_{\varepsilon},W^{1,p}(\D,X))$, $\xi:=\partial_\tau F(0,\cdot)$ is a holomorphic section of $f_0^*T^{1,0}X$ of class $W^{1,p}$.
    By Sobolev embedding $W^{1,p}(\D)\hookrightarrow C^{0,1-\frac{2}{p}}(\overline{\D})$, we get $\xi$ continuously extends to $\overline{\D}$.
\end{remark}

\vspace{0.25cm}
\subsection{A lower bound of the second variation $\partial_\tau\partial_{\bar\tau}\Psi(f_\tau)|_{\tau=0}$}

\begin{proposition}
\label{prop:levi}
Assume
\[
 F:D_\varepsilon\times\D\to X
\]
is holomorphic. Assume either:
\begin{enumerate}
\item $F$ extends holomorphically to
$D_\varepsilon\times D_R$ for some $R>1$;
\item $F(\tau,\cdot)\in W^{1,p}(\D,X)$ for some $p>2$ and 
\[
\tau\longmapsto
 F(\tau,\cdot)\in W^{1,p}(\D,X)
\]
is $C^2$ on $D_{\varepsilon}$.
\end{enumerate}
Let $f_\tau:=F(\tau,\cdot)$ and
$ \xi(z):= \partial_\tau F(0,z)
       \in T^{1,0}_{f_0(z)}X$
 be the variation vector field of $f_0$,
 then 
\[\tau\mapsto\Psi(f_\tau):=\int_{\D^2}G(z_1,z_2)h(\partial_{z_1}F,\partial_{z_1}F)(\tau,z_1)h(\partial_{z_2}F,\partial_{z_2}F)(\tau,z_2)dA_{z_1}dA_{z_2}\] 
is a $C^2$ function and 
\begin{equation}\label{eq:levi}
 \left. \partial_\tau\partial_{\bar\tau}\Psi(f_\tau)\right|_{\tau=0}
 \ge -2\pi A(f_0)\max_{z\in \partial\D}h(\xi,\xi)(z).
\end{equation}

\end{proposition}

We start with a simple statment.
\begin{lemma}
\label{lem:green-bilinear}
For every $q>1$, the integral
\begin{equation}\label{eq:bilinear}
 B(a,b):=\int_{\D^2}G(z_1,z_2)a(z_1)b(z_2)\,dA_{z_1}\,dA_{z_2}
\end{equation}
is absolutely convergent for $a,b\in L^q(\D)$ and defines
a symmetric bilinear form satisfying
\[
 |B(a,b)|\le C_q\|a\|_{L^q}\|b\|_{L^q}.
\]
\end{lemma}
\begin{proof}
Put $q':=q/(q-1)$. The bound
$0\le-G(z_1,z_2)\le\log(4/|z_1-z_2|^2)$ gives
\[
 \sup_{z_2\in\D}\|G(\cdot,z_2)\|_{L^{q'}}\le C_q.
\]
Hence H\"older's inequality yields
\[
 \begin{aligned}
 \int_{\D^2}|G(z_1,z_2)a(z_1)b(z_2)|
 &\le C_q\|a\|_{L^q}\|b\|_{L^1}\\
 &\le C_q'\|a\|_{L^q}\|b\|_{L^q}.
 \end{aligned}
\]
Bilinearity is immediate, and symmetry follows from
$G(z_1,z_2)=G(z_2,z_1)$.
\end{proof}

\begin{proof}[Proof of Proposition~\ref{prop:levi}]

\textbf{Case 1. Family of holomorphic discs up to the boundary.}

Using $G$ is symmetric, we compute that
\begin{equation}\label{eq:two-terms}
\begin{aligned}
\partial_\tau\partial_{\bar\tau}\Psi(f_\tau)
=&\partial_{\tau}\partial_{\bar\tau}\int_{\D^2}G(z_1,z_2)h(\partial_{z_1}F,\partial_{z_1}F)(\tau,z_1)h(\partial_{z_2}F,\partial_{z_2}F)(\tau,z_2)\\
=&2\int_{\D^2}G(z_1,z_2)\partial_{\tau}\partial_{\bar\tau}h(\partial_{z_1}F,\partial_{z_1}F)(\tau,z_1)h(\partial_{z_2}F,\partial_{z_2}F)(\tau,z_2)\\
+&2\int_{\D^2}G(z_1,z_2)\partial_\tau h(\partial_{z_1}F,\partial_{z_1}F)(\tau,z_1)\partial_{\bar\tau}h(\partial_{z_2}F,\partial_{z_2}F)(\tau,z_2).
\end{aligned}
\end{equation}
By commutativity Lemma~\ref{lem:closed} and letting $\tau=0$ we get
\begin{align*}
\partial_\tau\partial_{\bar\tau}\Psi(f_\tau)|_{\tau=0}
    =&2\int_{\D^2}G(z_1,z_2)\partial_{z_1}\partial_{\bar {z}_1 }h(\xi,\xi)(z_1)h(\partial_{z_2}f_0,\partial_{z_2}f_0)(z_2)\\
    +&2\int_{\D^2}G(z_1,z_2)\partial_{z_1} h(\xi,\partial_{z_1}f_0)(z_1)\partial_{\bar {z}_2 }h(\partial_{z_2}f_0,\xi)(z_2).
\end{align*}
Write $S:=h(\xi,\xi)$, which is a smooth function in a neighborhood of $\overline{\D}$, and let $T\in C^{\infty}(\overline{\D})$ be such that $T|_{S^1}=0$ and
$\partial_{z_1}\partial_{\bar z_1} T=\partial_{z_1}\partial_{\bar z_1}S$. For fixed $z_2\in\D$, using $\partial_{z_1}\partial_{\bar z_1}G(z_1,z_2)=\pi\delta_{z_2}$ we get
\[
 \begin{aligned}
 \int_\D G(z_1,z_2)\,\partial_{z_1}\partial_{\bar z_1}S(z_1)\,dA_{z_1}
 &=\int_\D G(z_1,z_2) \partial_{z_1}\partial_{\bar z_1}T(z_1)\,dA_{z_1}\\
 &=\int_\D T(z_1)\partial_{z_1}\partial_{\bar z_1}G(z_1,z_2)\,dA_{z_1}
   =\pi T(z_2),
 \end{aligned}
\]
where the integration by parts uses the fact that $T,G$ are both smooth up to boundary and $T|_{\partial\D}=G|_{\partial \D}=0$.
The first term therefore becomes
\begin{equation}\label{eq:diagonal}
 \int_{\D^2}G(z_1,z_2)\partial_{z_1}\partial_{\bar{z}_1 }h(\xi,\xi)(z_1)h(\partial_{z_2}f_0,\partial_{z_2}f_0)(z_2)
 =\pi\int_\D T(z_2)h(\partial_{z_2}f_0,\partial_{z_2}f_0)(z_2)\,dA_{z_2}.
\end{equation}
For the second term, write $g(z):=\partial_{z}h(\xi,\partial_{z}f_0)(z)$, which is a smooth function in a neighborhood of $\overline{\D}$, and
\[
 u(z_2):=\frac1\pi\int_\D G(z_1,z_2)g(z_1)\dA_{z_1}.
\]
is smooth up to boundary, satisfying $u|_{\partial\D}=0$ and
$\partial_{z_2}\partial_{\bar {z}_2}u=g$ on $\D$.
The second term becomes
\begin{equation}\label{eq:energy-identity}
\begin{aligned}
 \int_{\D^2}G(z_1,z_2)g(z_1)\overline{g(z_2)}dA_{z_1}dA_{z_2}
   &=\pi\int_\D u(z_2)\overline{g(z_2)}\dA_{z_2}\\
&=\pi\int_{\D}u(z_2)\overline{\partial_{z_2}\partial_{\bar z_2}u(z_2)}dA_{z_2}\\
     &=-\pi\int_\D|\partial_{\bar {z}_2}u|^2\dA_{z_2},
       \end{aligned}
       \end{equation}
       On the other hand,
       \begin{equation}   \label{eq:cs}
       \begin{aligned}
 \int_\D|\partial_{\bar z}u|^2\dA
   &=-\int_\D u\bar g\dA\\
   &=-\int_{\D}u\overline{\partial_z h(\xi,\partial_zf_0)}dA
     =\int_\D(\partial_{\bar z}u)\overline{h(\xi,\partial_zf_0)}\dA
     \le\|\partial_{\bar z}u\|_{L^2}\|h(\xi,\partial_zf_0)\|_{L^2},
\end{aligned}
\end{equation}
where we use the fact that the integral $\int_\D(\partial_{\bar z}u)\overline{h(\xi,\partial_zf_0)}\dA$ is real-valued. Consequently
\[
\int_{\D^2}G(z_1,z_2)g(z_1)\overline{g(z_2)}dA_{z_1}dA_{z_2}\ge-\pi\int_\D|h(\xi,\partial_zf_0)|^2(z)\dA_z.
\]
Since $|h(\xi,\partial_zf_0)|^2\le h(\xi,\xi)h(\partial_zf_0,\partial_zf_0)$ point-wise, we deduce that
\begin{align*}
 \partial_\tau\partial_{\bar\tau}\Psi(f_\tau)\Big|_{\tau=0}
 \ge2\pi\int_\D\Big(h(\xi,\xi)h(\partial_zf_0,\partial_zf_0)-|h(\xi,\partial_zf_0)|^2\Big)\dA
       &-2\pi\int_\D (h(\xi,\xi)-T) h(\partial_zf_0,\partial_zf_0)\dA\\
       \ge &-2\pi\int_\D(h(\xi,\xi)-T) h(\partial_zf_0,\partial_zf_0)\dA.
\end{align*}
Since $h(\xi,\xi)-T$ is harmonic on $\D$ and $(h(\xi,\xi)-T)|_{\partial\D}=h(\xi,\xi)|_{\partial\D}$, we get to know
\[
 \partial_\tau\partial_{\bar\tau}\Psi(f_\tau)\Big|_{\tau=0}
 \ge-2\pi A(f_0)\sup_{\partial\D}h(\xi,\xi).
\]
The proof is complete.

\medskip
\textbf{Case 2. Family of $W^{1,p}$ holomorphic discs.}

We use the isometry embedding $\iota:X\to \mathbb R^N$, $\iota\circ F$ corresponds to $F^1,\cdots,F^N\in C^2(D_\varepsilon,W^{1,p}(\D))$.
By equality \eqref{eq:embedding} we get
\begin{equation}\label{eq:embed Psi}
\Psi(f_\tau)=\int_{\D^2}G(z_1,z_2)\sum_{i=1}^N\Big(|\frac{\partial F^i}{\partial z_1}|^2+|\frac{\partial F^i}{\partial \bar z_1}|^2\Big)(\tau,z_1) \sum_{i=1}^N\Big(|\frac{\partial F^i}{\partial z_2}|^2+|\frac{\partial F^i}{\partial \bar z_2}|^2\Big)(\tau,z_2)
\end{equation}
$F^i\in C^2(D_{\varepsilon},W^{1,p}(\D))$ implies $\tau\mapsto \sum_{i=1}^N\Big(|\frac{\partial F^i}{\partial z}|^2+|\frac{\partial F^i}{\partial \bar z}|^2\Big)(\tau,z)$ belongs to $C^2(D_\varepsilon,L^{\frac{p}{2}}(\D))$. This in particular shows that $\Psi(f_\tau)$ is a $C^2$ function and equality~\eqref{eq:two-terms} still holds in this case:
\begin{align*}
    \partial_\tau\partial_{\bar\tau}\Psi(f_\tau)
    =&2\int_{\D^2}G(z_1,z_2) \partial_\tau\partial_{\bar\tau} h(\partial_{z_1}F,\partial_{z_1}F)(\tau,z_1)h(\partial_{z_2}F,\partial_{z_2}F)(\tau,z_2)\\
    +&2\int_{\D^2}G(z_1,z_2) \partial_\tau h(\partial_{z_1}F,\partial_{z_1}F)(\tau,z_1) \partial_{\bar\tau}h(\partial_{z_2}F,\partial_{z_2}F)(\tau,z_2).
\end{align*}
For $0<r<1$, set
\[
 F^{(r)}(\tau,z):=F(\tau,rz),
\]
which is holomorphic on
$D_\varepsilon\times D_{1/r}$, so Case~1. may apply.
Since $\partial_zF^{(r)}(\tau,z)=r\partial_zF(\tau,rz)$, 
we get $h(\partial_zF^{(r)},\partial_zF^{(r)})(\tau,z)=r^2h(\partial_zF,\partial_zF)(\tau,rz)$.
It follows that when $\tau$ is fixed, in $L^q(\D_z)$ we have
\begin{align*}
    &\lim_{r\uparrow1}h(\partial_zF^{(r)},\partial_zF^{(r)})(\tau,\cdot)=h(\partial_zF,\partial_zF)(\tau,\cdot),\\
     &\lim_{r\uparrow1}\partial_\tau h(\partial_zF^{(r)},\partial_zF^{(r)})(\tau,\cdot)=\partial_\tau h(\partial_zF,\partial_zF)(\tau,\cdot),  \\
    &\lim_{r\uparrow1}\partial_\tau \partial_{\bar\tau}h(\partial_zF^{(r)},\partial_zF^{(r)})(\tau,\cdot)=\partial_\tau\partial_{\bar\tau} h(\partial_zF,\partial_zF)(\tau,\cdot).
\end{align*}
Lemma~\ref{lem:green-bilinear} then implies
\[
\lim_{r\to 0} \partial_\tau\partial_{\bar\tau}\Psi(f_\tau^{(r)})=\partial_\tau\partial_{\bar\tau}\Psi(f_\tau).
\]
At $\tau=0$, Case~1. gives us 
\[
\partial_\tau\partial_{\bar\tau}\Psi(f_\tau^{(r)})\Big|_{\tau=0}\geq -2\pi A(f_0^{(r)})\sup_{z\in\partial\D}h(\partial_\tau F^{(r)},\partial_\tau F^{(r)})(0,z)
=-2\pi A(f_0^{(r)})\sup_{z\in \partial\D}h(\xi,\xi)(rz)
\]
In Remark~\ref{rem:C2 family} it is proved that $\xi\in W^{1,p}(\D,f_0^*T_X)$ extends continuously to $\overline\D$.
So we get
\[
 \sup_{z\in\partial\D}h(\xi,\xi)(rz)
 \longrightarrow\sup_{z\in\partial\D}h(\xi,\xi)(z).
\]
The convergence $ A(f^{(r)}_0)=\int_{D_r}f_0^*\omega\longrightarrow A(f_0),
 \qquad$ is clear. 
 Letting $r\uparrow1$, we complete the proof.
\end{proof}

\begin{remark}
We discuss an equality case for \eqref{eq:levi}.  
Consider the family $f_\tau:\D\to\mathbb{P^1},\,z\mapsto[\tau z:1]$. One obtains
$\Psi(f_\tau)=-\frac{\pi^2}{2}|\tau|^4$,
$\partial_\tau\partial_{\bar\tau}\Psi(f_\tau)=-2\pi^2|\tau|^2$, $A(f_\tau)=\pi|\tau|^2$ and
$\sup_{\partial\D}|\partial_\tau f_\tau|^2=1$.
\end{remark}

\vspace{0.25cm}
\subsection{Variations for $A(f_\tau),W(f_\tau),Q_\sigma(f_\tau),J_{t,\beta}(f_\tau)$}

For a given holomorphic disc $f:\D\to X$ of class $W^{1,p}$, we consider family of holomorphic discs $F:D_\varepsilon\times\D\to X $ so that $F(0,\cdot)=f$, $\partial_\tau F|_{\tau=0}=\xi$ and compute the variations of functionals $A,W,Q_\sigma,J_{t,\beta}$ at $\tau=0$. 
\begin{lemma}
\label{lem:variation-identities}
For every smooth closed real $(1,1)$-form $\rho$ and holomorphic $F:D_{\varepsilon}\times\D\to X$ so that $F\in C^2(D_\varepsilon; W^{1,p}(\D,X))$,
the functions $A(f_\tau)$, $W(f_\tau)$, $Q_\rho(f_\tau)$, and $J_{t,\beta}(f_\tau)$ are $C^2$ on $D_\varepsilon$.
Moreover, the second variations of these functions at $\tau=0$ only depend on $f=F(0,\cdot)$ and $\xi=\partial_\tau F|_{\tau=0}$.
\end{lemma}

\begin{proof}
We prove $Q_\rho(f_\tau)=\frac{1}{\pi}\int_\D\log(1/|z|)f_\tau^*\rho$ is $C^2$ and the proof of that of $A(f_\tau),W(f_\tau),J_{t,\beta}(f_\tau)$ is similar.
We can write $\rho=\omega_1-\omega_2$ for some K\"ahler forms $\omega_1,\omega_2$, by linearity we may assume $\rho$ is K\"ahler.
Let $\iota_\rho:X\to \mathbb{R}^M$ be an isometry embedding of $(X,\rho)$ and $F^i$ be the component of $\iota_\rho\circ F$, then
\[
|d_zF|^2=h^\rho(\partial_zF,\partial_zF)=\sum_{i=1}^M\big|\frac{\partial F^i}{\partial z}\big|^2+\big|\frac{\partial F^i}{\partial\bar z}\big|^2
\]
and therefore
\[Q_\rho(f_\tau)=\frac{1}{\pi}\int_\D\log(1/|z|)|d_zF|^2=\frac{1}{\pi}\int_\D\log(1/|z|)\sum_{i=1}^M\Big(\big|\frac{\partial F^i}{\partial z}\big|^2+\big|\frac{\partial F^i}{\partial\bar z}\big|^2\Big)(\tau,z)
\]
$F^i\in C^2(D_\varepsilon;W^{1,p}(\D))$ implies $\tau\longrightarrow |\frac{\partial F^i}{\partial z}|^2+|\frac{\partial F^i}{\partial\bar z}|^2\in C^2(D_{\varepsilon};L^{\frac{p}{2}}(\D))$ and $\log|z|\in L^q(\D)$ for any $q>1$, we can infer $Q_\rho(f_\tau)$ is $C^2$ on $D_\varepsilon$ and by equation~\eqref{eq:actual-local-variations},
\begin{align*}
\partial_\tau\partial_{\bar\tau}Q_\rho(f_\tau)\Big|_{\tau=0}
=\frac{1}{\pi}\int_\D \log(1/|z|)\Big(\partial_\tau\partial_{\bar\tau}h^\rho(\partial_zF,\partial_zF)\Big)(0,z)
&=\frac{1}{\pi}\int_\D \log(1/|z|)\Big(\partial_z\partial_{\bar z}h^\rho(\partial_\tau F,\partial_\tau F)\Big)(0,z)\\
&=\frac{1}{\pi}\int_\D \log(1/|z|)\partial_z\partial_{\bar z}h^\rho(\xi,\xi)(z).
\end{align*}
This completes our proof.
\end{proof}

Recall that
\[
 \sigma:=\theta-\delta\omega,\qquad
 q(f):=Q_\sigma(f)-cW(f)
 =\frac{1}{2\pi}\int_\D\log(1/|z|)\,f^*\sigma
  -c\int_\D (1-|z|^2)\,f^*\omega.
\]
We denote the first variation of $q(f_\tau)$ by
\[
 \ell_f(\xi)
 :=\left.\partial_\tau q(f_\tau)\right|_{\tau=0}.
\]

\begin{lemma}
\label{lem:first-variation-q}
Assume $\xi(0)=0$. The first variation $\ell_f$ satisfies
\begin{equation}\label{eq:bulk-ell}
 \ell_f(\xi)
 =\frac1{4\pi}\int_\D\frac{h^\sigma(\xi,\partial_zf)}z\,\dA
  -c\int_\D\bar z\,h(\xi,\partial_zf)\,\dA.
\end{equation}
In particular, $\ell_f(\xi)$ depends only on $f=F(0,\cdot)$ and $\xi=\partial_\tau F|_{\tau=0}$.
\end{lemma}

\begin{proof}
Differentiation under the integral and 
\eqref{eq:actual-local-variations} give
\[
 \ell_f(\xi)
 =\frac1\pi\int_\D\lambda\,\partial_z h^\sigma(\xi,f_z)\,\dA
  -c\int_\D w\,\partial_z h(\xi,f_z)\,\dA.
\]
Since $\partial_z\lambda=-\frac1{2z}$ and $\partial_zw=-\bar z$, integration by parts yields
\[
 \begin{aligned}
 \int_\D\lambda\,\partial_z h^\sigma(\xi,f_z)\,\dA
 &=-\int_\D(\partial_z\lambda)h^\sigma(\xi,f_z)\,\dA
   =\frac12\int_\D\frac{h^\sigma(\xi,f_z)}z\,\dA,\\
 \int_\D w\,\partial_z h(\xi,f_z)\,\dA
 &=-\int_\D(\partial_zw)h(\xi,f_z)\,\dA
   =\int_\D\bar z\,h(\xi,f_z)\,\dA.
 \end{aligned}
\]

For example, to justify integration by parts for $\lambda$,
choose smooth cutoff functions $0\le\eta_\varepsilon,\chi_r\le1$
such that $\eta_\varepsilon=1$ on $D_{1-2\varepsilon}$
and vanishes for $|z|\ge1-\varepsilon$, while $\chi_r=0$
on $D_r$ and equals $1$ outside $D_{2r}$, with
\[
 |D\eta_\varepsilon|\le C/\varepsilon,
 \qquad |D\chi_r|\le C/r.
\]
Since $\eta_\varepsilon\chi_r$ is compactly supported
away from zero, we can write
\[
 \begin{aligned}
 \int_\D\eta_\varepsilon\chi_r\lambda\,
                    \partial_z h^\sigma(\xi,f_z)\,\dA
 ={}&\frac12\int_\D\eta_\varepsilon\chi_r
                    \frac{h^\sigma(\xi,f_z)}z\,\dA\\
 &-\int_\D\chi_r\lambda h^\sigma(\xi,f_z)\,
                    \partial_z\eta_\varepsilon\,\dA\\
 &-\int_\D\eta_\varepsilon\lambda h^\sigma(\xi,f_z)\,
                    \partial_z\chi_r\,\dA.
 \end{aligned}
\]
On $A_\varepsilon:=\{1-2\varepsilon<|z|<1\}$,
$\lambda\le C\varepsilon$ and
$|h^\sigma(\xi,f_z)|\le C|df|\,|\xi|$. Thus
\[
 \left|\int_\D\chi_r\lambda h^\sigma(\xi,f_z)\,
                    \partial_z\eta_\varepsilon\,\dA\right|
 \le C\int_{A_\varepsilon}|df|\,|\xi|\,\dA
 \longrightarrow0,
\]
since $df\in L^p$ and $\xi$ is bounded.

Near $0$, we have
$\xi(z)=O(|z|)$ and $f_z$ is bounded, hence
$h^\sigma(\xi,f_z)=O(|z|)$. Consequently,
\[
 \begin{aligned}
 \left|\int_\D\eta_\varepsilon\lambda h^\sigma(\xi,f_z)\,
                    \partial_z\chi_r\,\dA\right|
 &\le\frac Cr\int_{D_{2r}\setminus D_r}\lambda|z|\,\dA\\
 &\le Cr^2\log(1/r)\longrightarrow0.
 \end{aligned}
\]
Letting $\varepsilon,r\downarrow0$, we deduce that
\[
 \int_\D\lambda\,\partial_z h^\sigma(\xi,f_z)\,\dA
 =\frac12\int_\D\frac{h^\sigma(\xi,f_z)}z\,\dA.
\]
Similar arguments, using $\partial_zw=-\bar z$,
gives $\int_\D w\,\partial_z h(\xi,f_z)\,\dA
=\int_\D\bar z\,h(\xi,f_z)\,\dA$, and the proof is concluded.
\end{proof}

\section{Boundary-unitary tangent frame}
\label{sec:second-variation}

\subsection{Coordinates of graph and realization of variation fields}

Let $f:\D\to X$ be a holomorphic disc of class $W^{1,p}$, $p>2$ and $\xi$ be a holomorphic section of $f^*T^{1,0}X$ of class $W^{1,p}$. We are interested in constructing a holomorphic family $F:D_{\varepsilon}\times \D\to X$ of class $W^{1,p}$ so that $\partial_\tau F(\tau,\cdot)=\xi$.

We need the fundamental graph theorem of Forstneri\v c in \cite[Theorem~1.2]{Forstneric}.
\begin{theorem}
    Assume that $h:Y \to S$ is a holomorphic submersion of a complex manifold $Y$ onto a Stein manifold S, $D \Subset S$ is a strongly pseudoconvex domain with $C^2$ boundary in $S$, and $f: \overline D \to Y $ is a continuous section of h which is holomorphic in $D$. 
    
    There exists a holomorphic vector bundle $\pi: E \to U$ over an open set $U \subset S$ containing $\overline D$, and for every open set $\Omega_0 \subset Y$ containing $f(\overline D)$, there exist a Stein open set $\Omega$ in $Y$ with $f(\overline  D) \subset \Omega \subset\Omega_0$ and a biholomorphic map $\Theta: \Omega\to \Theta(\Omega)\subset E$ which maps the fiber $\Omega_z=h^{-1}(z)\cap\Omega$ over any point $z \in h(\Omega) \subset U$ biholomorphically onto an open convex set $\Theta(\Omega_z)$ in the fiber $E_z = \pi^{-1}(z)$.
\end{theorem}

We use Forstneri\v c's theorem to construct (global) coordinate of the graph of a given holomorphic disc.
\begin{lemma}
\label{lem:graph-coordinates}
Let $f:\overline\D\to X$ be continuous and holomorphic on $\D$.
There are $\varepsilon>0$, a Stein neighborhood
$\Omega_f\subset D_{1+\varepsilon}\times X$ of
\[
 G_f:=\{(z,f(z)):z\in\overline\D\},
\]
and a fibre-preserving biholomorphism
\[
 \Theta_f:\Omega_f\longrightarrow
       \Omega_f'\subset D_{1+\epsilon}\times\mathbb C^n.
\]
\end{lemma}

\begin{proof}
Apply Forstneri\v c's theorem to the projection
$\mathbb C\times X\to\mathbb C$ and the section
$z\mapsto(z,f(z))$ over $\overline\D$.
It gives a holomorphic vector bundle $E\to U$, where $U$ is
an open neighborhood of $\overline\D$,
and we can choose $\varepsilon>0$ such that
$\overline{D_{1+\varepsilon}}\subset U$.
Over $D_{1+\varepsilon}$ there is a Stein neighbourhood $\Omega_f$ of $G_f$ and a fiber-preserving biholomorphism $\Theta_f$ so that $\Omega_f\xrightarrow{\Theta_f}\Omega_f':=\Theta_f(\Omega_f)\subset E$.

Since $D_{1+\varepsilon}$ is contractible, the restricted bundle
$E|_{D_{1+\varepsilon}}$ is topologically trivial.
By Oka--Grauert principle
\cite{Grauert} this bundle is holomorphically trivial.
It is clear that $\operatorname{rank}E=\dim X=n$. The proof is complete.
\end{proof}

\begin{lemma}
\label{lem:actual}
Let $f:\D\to X$ be holomorphic and of class $W^{1,p}$,
$p>2$, and let
$\xi\in W^{1,p}(\D,f^*T^{1,0}X)$ be a holomorphic section
with $\xi(0)=0$. There are $\varepsilon>0$ and a holomorphic map
\[
 F:D_\varepsilon\times\D\to X
\]
such that
\[
 F(0,\cdot)=f,\qquad
 \left.\partial_\tau F(\tau,z)\right|_{\tau=0}=\xi(z),
 \qquad F(\tau,0)=f(0),
\]
and
\[
 \tau\longmapsto F(\tau,\cdot)
 \in C^2(D_\varepsilon;W^{1,p}(\D,X)).
\]

\end{lemma}

\begin{remark}\label{rem:basic of W1p}
    We will frequently use the following facts:
    \begin{itemize}
        \item $W^{1,p}(\D)\hookrightarrow C^{0,1-\frac{2}{p}}(\overline{\D})$ and the multiplication $W^{1,p}(\D)\times W^{1,p}(\D)\to W^{1,p}(\D)$ is continuous,
        \item Let $f^1,\cdots,f^N\in W^{1,p}(\D)$, then $f=(f^1,\cdots,f^N):\D\to\mathbb{R}^N$ continuously extends to $\overline\D$.
        If $G$ is a smooth mapping defined in a neighborhood of $f(\overline\D)$, then $G\circ f\in W^{1,p}(\D,\mathbb{R}^N)$.
    \end{itemize}
\end{remark}

\begin{proof}
Lemma~\ref{lem:graph-coordinates} gives a biholomorphism
\[
 \Theta:\Omega\longrightarrow\Omega',
\]
where $\Omega\subset\mathbb C\times X$ is a neighborhood
of the graph $G_f:=\{(z,f(z)):z\in\overline\D\}$, $\Omega'\subset\mathbb C\times\mathbb C^n$, Moreover, $\Theta$ maps fiber $\Omega_z$ biholomorphically to $\Omega^\prime_z$, $\forall z\in D_{1+\varepsilon_0}$.
So we can write
\[
 \Theta(z,x)=(z,\mathcal G(z,x)),\qquad
 \Theta^{-1}(z,v)=(z,\mathcal K(z,v)).
\]
Define
\[
 a(z):=\mathcal G(z,f(z)),\qquad
 b(z):=d_x\mathcal G(z,x)\big|_{x=f(z)}\xi(z).
\]
Note that
$d_x\mathcal G(z,x):T^{1,0}X\to T^{1,0}\mathbb{C}^n$ is the differential of $x$-direction and $\xi(z)\in T^{1,0}_{f(z)}X$, so
$a,b$ are holomorphic $\mathbb C^n$-valued functions on $\D$.
Since $\mathcal G$ is defined in a neighborhood of $\{(z,f(z)):z\in\overline{\D}\}$, 
it follows from $f\in W^{1,p}(\D,X)$, $\xi\in W^{1,p}(\D,f^*T^{1,0}X)$ and
Remark \ref{rem:basic of W1p} that
\[
 a,b\in W^{1,p}(\D,\mathbb C^n)
       \cap C^0(\overline\D,\mathbb C^n).
\]
Since $\{(z,a(z)):z\in\overline\D\}=\Theta(G_f)$ is a compact subset of $\Omega'$
and $b$ is continuous on $\overline\D$, we can choose $\varepsilon>0$ so that
\[
 \big(z,a(z)+\tau b(z)\big)\in\Omega'
 \qquad(z\in\overline\D,\ \tau\in \overline D_\varepsilon).
\]
Set
\[
 F(\tau,z):=\mathcal K(z,a(z)+\tau b(z)).
\]
Note that $d_v\mathcal K(z,v)|_{v=a(z)}$ is the inverse of $d_x\mathcal G(z,x)|_{x=f(z)}$, we see
\[
 F(0,z)=f(z),\qquad
 \partial_\tau F(0,z)
 =d_v\mathcal K(z,v)\big|_{v=a(z)}\,b(z)=\xi(z).
\]
Moreover, $\xi(0)=0$ implies $b(0)=0$ and consequently $F(\tau,0)=f(0)$.

Finally, $\mathcal{K}$ is defined in a neighborhood of $\{(z,a(z)+\tau b(z)):(\tau,z)\in \overline{D_{\varepsilon}}\times\overline\D\}$,

Remark~\ref{rem:basic of W1p} then implies
\[
 \tau\longmapsto F(\tau,\cdot)
 \in C^2(D_\varepsilon;W^{1,p}(\D,X)).
\]
Indeed, let $\iota:X\to \mathbb{R}^N$ be an embedding, $(\iota\circ F)(z,a(z)+\tau b(z))$ corresponds to $N$ functions $f^i(z,a(z)+\tau b(z))$, where $f^i$ are smooth functions defined in a neighborhood of $\{(z,a(z)+\tau b(z))\in\mathbb C\times\mathbb C^n:(\tau,z)\in \overline{D_{\varepsilon}}\times\overline\D\}$.
A direct computation yields that $f^i\in C^{\infty}(D_{\varepsilon};W^{1,p}(\D))$. Thus $\tau\mapsto F(\tau,\cdot)\in W^{1,p}(\D,X)$ is actually of class $C^{\infty}$. The proof is complete.
\end{proof}

\vspace{0.25cm}
\subsection{Construct boundary-unitary tangent frame of class $W^{1,p}$}

Let $f:\D\to X$ be a holomorphic disc of class $W^{1,p}$,
where $p>2$. A holomorphic frame $e=(e_1,\ldots,e_n)$ of
$f^*T^{1,0}_X$ consists of holomorphic sections which form a basis
at each point. 
Clearly under a bundle isomorphism
$\D\times\mathbb C^n\to f^*T^{1,0}_X$, the standard basis $(\textbf{e}_1,\cdots,\textbf{e}_n)$ of $\mathbb{C}^n$ corresponds to a frame $e=(e_1,\cdots,e_n)$

At each $z\in\D$, the pullback metric $f^*h$ on the fiber over $z$ is $h_{f(z)}$. We obtain the Hermitian form
$(u,v)\mapsto h_{f(z)}(e(z)u,e(z)v)$ on $\mathbb C^n$.
The associated Gram matrix $G_e(z)$ is defined as
\[
 h_{f(z)}(e(z)u,e(z)v)=v^*G_e(z)u,
 \qquad (G_e)_{ij}(z)=h_{f(z)}(e_j(z),e_i(z)),
\]
where $^*$ denote the conjugate transpose.
We say that $e\in W^{1,p}$ if each
$e_j\in W^{1,p}(\D,f^*T^{1,0}X)$, and that $G_e\in W^{1,p}$
if each of its entries belong to $W^{1,p}(\D)$.

\begin{theorem}
\label{lem:baseline-frame}
Let $f:\D\to X$ be holomorphic and of class $W^{1,p}$, $p>2$. There is a holomorphic frame $e$ of $f^*T^{1,0}_X$
such that
\[
 e\in W^{1,p},
\]
the Gram matrix and its inverse $G_e,\,G_e^{-1}\in W^{1,p}$ and $G_e|_{S^1}=I_n$.
Such a frame will be called \emph{boundary-unitary}.
\end{theorem}

\begin{proof}
If $f:\D\to X$ extends smoothly to $\overline\D$, then the desired result can be deduced from a result of Wiener--Masani \cite{WM57}, see \cite[Section~2]{BR}.
The situation that $f$ is merely $W^{1,p}$ is quite subtle, we give a detailed proof which mainly follows the argument in \cite[Section~6]{BR}. 
While preparing this version of the paper, we learned that \cite[Proposition~4.2]{OAI26} obtained the same result with a different proof.

\medskip
\noindent\textbf{Step 1. Obtain a $W^{1,p}$ holomorphic frame by Forstneri\v c's theorem.}

Use Forstneri\v c's graph theorem for $G_f=\{(z,f(z)):z\in\overline\D\}$ as in Lemma~\ref{lem:actual}, we obtain isomorphisms
\[
 \Theta(z,x)=(z,\mathcal G(z,x)),
 \qquad
 \Theta^{-1}(z,v)=(z,\mathcal K(z,v)).
\]
between some neighborhood $\Omega\subset \mathbb C\times X$ of $G_f$ and $\Omega'\subset D_{1+\varepsilon}\times\mathbb{C}^n$.
Set
\begin{align*}
 a(z):&=\mathcal G(z,f(z)),\\
 T(z):&=d_v\mathcal K(z,v)\big|_{v=a(z)}:
       T_{a(z)}\mathbb C^n\longrightarrow T_{f(z)}X.
\end{align*}
$T(z)$ is isomorphism for every $z$, hence $T(\textbf{e}_j)$ form a frame of $f^*T^{1,0}_X$. Moreover, since $\mathcal G,H$ are defined in neighborhoods of $\{(z,f(z)):z\in\overline\D\}$, $\{(z,a(z)):z\in\overline\D\}$ respectively,
Remark~\ref{rem:basic of W1p} implies
\[
 a\in W^{1,p}\cap C^0(\overline\D),
 \qquad
 T(\textbf{e}_j)\in W^{1,p}\cap C^0(\overline\D).
\]
We use $G:=G_T$ to denotes the Gram matrix of $T(\textbf{e}_j)$, then $G$ is continuous up to boundary.
Let $(t_1,\cdots,t_n)$ be the corresponding sections of the frame $T$, we have
\[
 \partial_zG_{ij}
 :=\partial_z(h(t_j,t_i))
 =(\partial_zh)(t_j,t_i)+h(\partial_z t_j,t_i),
 \qquad 
 \partial_{\bar z}G_{ij}
 :=\partial_{\bar z}(h(t_j,t_i))
 =(\partial_{\bar z}h)(t_j,t_i)+h( t_j,\partial_{ z}t_i),
\]
where we abbreviate the notation $h_{f(z)}$ to $h$.
Each term has one $L^p$ factor and two bounded factors.
Hence
\[
 G\in W^{1,p}(\D,\operatorname{Mat}_{n\times n}(\mathbb C)).
\]
Since $a\in C^0(\overline\D)$ and $H$ is defined in a neighborhood of $\{(z,a(z)):z\in\overline\D\}$, there are constants
$0<m\le M<\infty$ such that
\[
 mI_n\le G(z)\le MI_n
 \qquad(z\in\overline\D).
\]
We write $G_b:=G|_{S^1}$, which is continuous on $S^1$.

\vspace{0.25cm}
The following strategy, which uses Hardy space to construct holomorphic matrix $E,E^{-1}\in W^{1,p}(\D)$ so that $E^*GE=I_n$ on $S^1$, follow the line in \cite[Proposition~6.1]{BR}.

\medskip
\noindent\textbf{Step 2. Hardy space with an inner product induced by $G_b$.}

We consider the Hardy space:
\[
 \mathcal H:=H^2(\D,\mathbb C^n)
 :=\left\{u(z)=
\begin{pmatrix}
\sum_{k=0}^{\infty} a_{1k}z^k\\
\sum_{k=0}^{\infty} a_{2k}z^k\\
\vdots\\
\sum_{k=0}^{\infty} a_{nk}z^k
\end{pmatrix},
z\in\mathbb D,\qquad
\|u\|_{\mathcal H}^2:= \sum_{i=1}^n\sum_{k=0}^{\infty}|a_{ik}|^2<\infty\right\}.
\]
with $a_k:=
\begin{pmatrix}
a_{1k}\\ a_{2k}\\ \vdots\\ a_{nk}
\end{pmatrix},
$
we can write $u=\sum_{k\geq0} a_kz^k$ and
$u_b(e^{i\theta}):=\sum_{k\ge0}a_ke^{ik\theta}\in L^2(S^1,\mathbb{C}^n)$
is called the \emph{boundary value} of $u$. A useful property is that
\begin{equation}\label{eq: L2 appro}
u(re^{i\theta})\to u_b(e^{i\theta}) \qquad \text{in $L^2(S^1,\mathbb C^n)$}.
\end{equation}
Since every component of $u$ is holomorphic on $\D$, this approximation property shows Poisson formula holds for $u\in\mathcal H$ component-wise:
\begin{equation}\label{eq:Poisson}
u(z)=\frac{1}{2\pi}\int_0^{2\pi}\frac{1-|z|^2}{|e^{i\theta}-z|^2} u_b(e^{i\theta})\, d\theta.
\end{equation}
Moreover, if $w$ is a bounded holomorphic function on $\D$, then $w\in H^2(\D,\mathbb C)$ and $w_b$ is also bounded. For every $u\in\mathcal H$, $wu\in\mathcal H$ and $(wu)_b=w_bu_b$. 

Equip $\mathcal H$ with the inner product
\[
 \langle u,v\rangle_G
 :=\frac1{2\pi}\int_0^{2\pi}
 u_b(e^{i\theta})^*G_b(e^{i\theta})
 v_b(e^{i\theta})\,d\theta.
\]
It is easy to see that
\[
 \langle u,v\rangle_G= \langle zu,zv\rangle_G,\qquad u,v\in\mathcal H.
\]
Since $mI_n\le G_b\le MI_n$,
\[
 m\|u\|_{\mathcal H}^2\le\|u\|_G^2\le M\|u\|_{\mathcal H}^2.
\]
Thus $(\mathcal H,G)$ is a Hilbert space.

The subspace $z\mathcal H=\{u\in\mathcal H:u(0)=0\}$ is clearly closed and has codimension $n$.
Consider its orthocomplement in $(\mathcal H,G)$
\[
 \mathcal W:=(z\mathcal H)^{\perp_G},
 \qquad
 \mathcal H=\mathcal W\oplus_G z\mathcal H.
\]
Then $\dim_{\mathbb C}\mathcal W=n$; we may choose a $G$-orthonormal basis
$u_1,\ldots,u_n$ of $\mathcal W$. We prove that $\{z^ku_j:k\geq0,j=1,\cdots,n\}$ is a complete $G$-orthonormal basis of $\mathcal H$.

Let $u\in\mathcal H$, there is a unique orthogonal decomposition
\begin{align*}
&u=(c_{0,1}u_1+\cdots+c_{0,n}u_n)+zr_1,\qquad c_i\in\mathbb C, r_1\in\mathcal H,\\
&\|u\|^2_G=\|c_{0,1}u_1+\cdots+c_{0,n}u_n\|^2_G+\|zr_1\|^2_G.
\end{align*}
An induction construction yields
\begin{align*}
    u=(c_{0,1}u_1+\cdots+c_{0,n}u_n)+\cdots+z^{k-1}(c_{k-1,1}u_1+\cdots+c_{k-1,n}u_n)+z^kr_k,\\
    \|u\|^2_G=\|c_{0,1}u_1+\cdots+c_{0,n}u_n\|^2_G+\cdots+\|z^{k-1}(c_{k-1,1}u_1+\cdots+c_{k-1,n}u_n)\|_G^2+\|z^kr_k\|^2_G
\end{align*}
and the terms on right-hand side are orthogonal to each other.
This orthogonality makes $z^kr_k$ a Cauchy sequence in $\mathcal H$. So we let $r:=\lim z^kr_k$ and obtain that
\begin{align*}
    u=\sum_{k=0}^\infty z^k(c_{k,1}u_1+\cdots+c_{k,n}u_n)+r
\end{align*}
However $z^{k_0}\mathcal H=\{u=\sum_{k\geq k_0}a_kz^k\in\mathcal H\}$ is closed in $\mathcal H$, so we see $r\in z^{k_0}\mathcal H$ for every $k_0$. This proves $r=0$.
Consequently
\[
 \{z^ku_j:k\ge0,\ 1\le j\le n\}
\]
is a complete $G$-orthonormal basis of $\mathcal H$.

\medskip
\noindent\textbf{Step 3. Construct boundary-unitary frame.}

Set
\[
 E(z):=(u_1(z),\ldots,u_n(z)),
 \qquad Q:=E_b^*G_bE_b,
\]
where $E_b:=(u_{1,b},\cdots,u_{n,b})$.
Each entry of $Q$ belongs to $L^1(S^1)$, since the entries
of $E_b$ belong to $L^2(S^1)$ and $G_b$ is bounded on $S^1$.
By definition
$Q_{ij}(e^{i\theta})=u_{i,b}(e^{i\theta})^*G_b(e^{i\theta})u_{j,b}(e^{i\theta})$.
Orthogonality of $\{z^ku_i:k\geq0,1\leq j\leq\}$ implies that, for $k,\ell\geq0$ and $1\leq i,j\leq n$,
\[
 \begin{aligned}
 \frac1{2\pi}\int_0^{2\pi}
 e^{i(\ell-k)\theta}Q_{ij}(e^{i\theta})\,d\theta
 &=\frac1{2\pi}\int_0^{2\pi}e^{i(\ell-k)\theta}u_{i,b}(e^{i\theta})^*G_b(e^{i\theta})u_{j,b}(e^{i\theta})\,d\theta\\
 &=\frac1{2\pi}\int_0^{2\pi}(z^ku_{i,b})(e^{i\theta})^*G_b(e^{i\theta})(z^lu_{j,b})(e^{i\theta})\,d\theta\\
 &=\langle z^ku_i,z^\ell u_j\rangle_G\\
 &=\delta_{k\ell}\delta_{ij}.
 \end{aligned}
\]
Since $\ell-k$ ranges over all integers, these identities
say that the Fourier coefficients of $Q_{ij}$ and $\delta_{ij}$ are identical.
Uniqueness of Fourier coefficients in $L^1(S^1)$ yields
\[
 E_b^*G_bE_b=I_n
 \quad\text{in $L^1(S^1)$}.
\]
Combining this with $G_b\ge mI_n$, we obtain every diagonal entry
\[
 1=u_{j,b}^*G_bu_{j,b}
 \ge m|u_{j,b}|^2
 \qquad\text{a.e. on $S^1$}.
\]
The Poisson representation \eqref{eq:Poisson} then gives
\[
 \begin{aligned}
 |u_{j}(z)|
 &\le\frac1{2\pi}\int_0^{2\pi}
 \frac{1-|z|^2}{|e^{i\theta}-z|^2}
 |u_{j,b}(e^{i\theta})|\,d\theta\\
 &\le m^{-1/2},
 \qquad z\in\D.
 \end{aligned}
\]
Hence every entry of $E$ is bounded and we obtain
\begin{align*}
    E_b^*G_bE_b=I_n\qquad\text{in $L^\infty(S^1)$}.
\end{align*}

\medskip
\noindent\textbf{Step 4. Prove $E^{-1}\in \mathcal H$ and $E^{-1}$ is bounded.}

Since entry of holomorphic matrix $E$ are bounded, multiplication by $E$ maps $\mathcal H$ into itself and $E_b^*G_bE_b=I_n$ gives
\[
 \|Eu\|_G^2
 =\frac1{2\pi}\int_0^{2\pi}
 u_b(e^{i\theta})^*E_b(e^{i\theta})^*
 G_b(e^{i\theta})E_b(e^{i\theta})
 u_b(e^{i\theta})\,d\theta
 =\frac{1}{2\pi}\int_0^{2\pi}u_b(e^{i\theta})^*u_b(e^{i\theta})\,d\theta
 =\|u\|_{H^2}^2.
\]
The range of $E:\mathcal H\to \mathcal H$ is therefore closed.
Since $E(z^k\mathbf e_j)=z^ku_j$, where $\mathbf e_j$ is the standard basis, the range contains
the complete basis constructed in Step~2 and hence equals $\mathcal H$.

For each $\mathbf e_j$, there is
$a_j\in\mathcal H$ so that $Ea_j=\mathbf e_j$.
Setting $A:=(a_1,\ldots,a_n)$ we obtain
\[
 EA=I_n\quad\text{on }\D.
\]
Thus $E$ is holomorphically invertible and $A=E^{-1}$.

Since entries of $E,A$ are in Hardy space, entries of $E$ are bounded, taking boundary values gives $E_bA_b=I_n$ in $L^2(S^1)$.
The boundary identity $E_b^*G_bE_b=I_n$ in $L^{\infty}(S^1)$ from Step~3 multiplied by $A_b^*,A_b$ yields
\begin{align*}
 G_b=A_b^*A_b,
 \qquad \text{in  $L^1(S^1)$}
 \end{align*}
 and therefore by $G_b\leq MI_n$ on $S^1$, the diagonal entry satisfies
 \begin{align*}
 \sum_{i=1}^n|A_{i j,b}|^2
 =(G_{b})_{jj}\le M
\end{align*}
Hence $|A_{i j,b}|\le M^{1/2}$ a.e.\ on $S^1$.
Applying the Poisson representation \eqref{eq:Poisson} as in Step~3 gives
\[
 |A_{i j}(z)|\le M^{1/2}
 \qquad(z\in\D).
\]
Therefore $E^{-1}=A$ is bounded.

Up to now, using $G_b\in L^{\infty}(S^1,\text{Mat}_{n\times n}(\mathbb{C}))$ with $mI_n\leq G\leq MI_n$ we have constructed $E\in H^2(\D,\text{Mat}_{n\times n}(\mathbb{C}))$ so that $|E_{i,j}|\leq m^{-\frac{1}{2}}$, $|E^{-1}_{ij}|\leq M^{\frac{1}{2}}$ and $E_b^*G_bE_b=I_n$ on $S^1$.

\medskip
\noindent\textbf{Step 5. Prove $E,E^{-1}\in W^{1,p}$ using $G\in W^{1,p}$.}

On $A_+:=\{1<|z|<2\}$, consider the reflection of $G$ and $E$.
\[
 G_{\rm ref}(z):=G(1/\bar z),\qquad
 E_+(z):=G_{\rm ref}(z)^{-1}[E(1/\bar z)^*]^{-1}.
\]
$E(1/\bar z)^*$ and its inverse are clearly holomorphic.
Since $G_{\rm ref}$ is smooth on $A_+$, we have $\bar\partial G_{\rm ref}^{-1}=-G_{\rm ref}^{-1}\bar\partial G_{\rm ref} G_{\rm ref}^{-1}$. Moreover $G_{\rm ref}\in W^{1,p}(A^+)$ gives
\[
 \bar\partial E_+
 =-G_{\rm ref}^{-1}(\bar\partial G_{\rm ref}) G_{\rm ref}^{-1}[E(1/\bar{z})^*]^{-1}
 =-G_{\rm ref}^{-1}(\bar\partial G_{\rm ref})E_+
 =:Q\in L^p(A_+).
\]

The fact that $E_b^*G_bE_b=I_n$ gives
$G_b^{-1}(E_b^{-1})^*=E_b$.
$G$ is continuous up to the boundary, Hardy boundary
convergence \eqref{eq: L2 appro} of $E,E^{-1}$ therefore implies
\begin{equation}\label{eq:frame-matching-boundary-limits}
 \begin{aligned}
 E(re^{i\theta})&\longrightarrow E_b(e^{i\theta}),\\
 E_+(r^{-1}e^{i\theta})
 &=G(re^{i\theta})^{-1}[E^{-1}(re^{i\theta})]^*
 \longrightarrow G_b(e^{i\theta})^{-1}(E_b^{-1}(e^{i\theta}))^*=E_b(e^{i\theta})
 \end{aligned}
\end{equation}
in $L^2(S^1)$ as $r\uparrow1$.

On $A:=\{1/2<|z|<2\}$, we extend $E$ by reflection:
\[
 \widetilde E(z):=
 \begin{cases}
 E(z),&|z|<1,\\
 E_+(z),&|z|>1.
 \end{cases}
\]
For smooth function $\varphi\in C_c^\infty(A)$ and $3/4<r<1$, put
\[
 A_r^-:=\{1/2<|z|<r\},\qquad
 A_r^+:=\{r^{-1}<|z|<2\}.
\]
Stokes' formula gives
\[
 \begin{aligned}
 &\int_{A_r^-}E\,\partial_{\bar z}\varphi\,dA
 +\int_{A_r^+}E_+\,\partial_{\bar z}\varphi\,dA\\
 &\qquad=-\int_{A_r^+}Q\varphi\,dA
 +\frac1{2i}\left(
   \int_{S_r}E\varphi\,dz
   -\int_{S_{r^{-1}}}E_+\varphi\,dz\right).
 \end{aligned}
\]
By \eqref{eq:frame-matching-boundary-limits}, both boundary
term converge to $\int_{S^1}E_b\varphi\,dz$.
Letting $r\uparrow1$, we obtain that
\[
 \int_A\widetilde E\,\partial_{\bar z}\varphi\,dA
 =-\int_{A_+}Q\varphi\,dA,
 \qquad
 \text{which means }\quad \bar\partial\widetilde E=\mathbf1_{A_+}Q
 \quad\text{in distributions on }A.
\]
The standard interior
$\bar\partial$ estimate on $\mathbb C$ gives, for \(A'=\{3/4<|z|<3/2\}\subset\subset A\),
\[
 \|\widetilde E\|_{W^{1,p}(A')}
 \le C_p\left(
 \|\widetilde E\|_{L^p(A)}+\|Q\|_{L^p(A_+)}
 \right)<\infty.
\]
We therefore control the $W^{1,p}$-norm of $E$ in a neighborhood of $\partial\D$. Since $E$ is holomorphic on $\D$, it follows that $E\in W^{1,p}(\D)$.
$E^{-1}$ is also holomorphic on $\D$, we obtain
\[
 \partial_j(E^{-1})
 =-E^{-1}(\partial_jE)E^{-1}\in L^p(\D),
 \qquad j=x,y,
\]
so $E^{-1}\in W^{1,p}(\D)$ as well.

Note that if $G$ has higher regularity, the estimates of $\bar\partial$-operator will give higher regularity of $E,E^{-1}$, while we do not persist this issue here.

\vspace{0.2cm}
We finish the proof. Define holomorphic frame $e$ of $f^*T^{1,0}_X$ by
\[
 e:=TE:\D\times\mathbb C^n\to f^*T^{1,0}_X.
\]
$W^{1,p}$ is closed under multiplication, so $e$ is of class  $W^{1,p}$.
Moreover by construction, $G_e=E^*GE$ and $G_e^{-1}=E^{-1}G^{-1}(E^{-1})^*$ are of class $W^{1,p}$, moreover
\[
 G_e|_{S^1}
 =E_b^*G_bE_b
 =I_n.
\]
This completes the proof.
\end{proof}

\section{Proof of Theorem~\ref{thm:pointwise-criterion}}
\label{sec:exact-fields}

Let $f_{t,\beta}$ be a minimizer of $J_{t,\beta}$. 
The key ingredient to establish Theorem~\ref{thm:pointwise-criterion} is to give a uniform bound of the weighted area $W(f_{t,\beta})$ independent of $t,\beta$.

Let $e_1,\cdots,e_n\in W^{1,p}(\D,f^*T_X)$ be a holomorphic frame so that the Gram matrix and its inverse $G_e,G_e^{-1}\in W^{1,p}$ and $G_e|_{S^1}=I_n$.
For $b\in\D$, consider
\[
 B_b(z):=\frac{z-b}{1-\bar bz},\qquad z\in\D.
\]
It is holomorphic near $\overline\D$, has modulus one on $S^1$. We shall use the following two kinds of diagonal matrices:
\[
 \Phi(z):=zI_n,
 \qquad\text{or}\qquad
 \Phi(z):=\operatorname{diag}(zB_b(z),z,\ldots,z).
\]
The holomorphic sections $\xi:=e\Phi$ have the form
\[
 (\xi_1,\cdots,\xi_n)=(ze_1,\cdots,ze_n),
 \qquad\text{or}\qquad
 (\xi_1,\cdots,\xi_n)=(zB_b(z)e_1,ze_2,\cdots,ze_n)
\]
Then $\xi_j(0)=0$ and their Gram matrix $\Phi^*G_e\Phi=:H$ still has boundary value $I_n$.
Let $f$ be a holomorphic disc of class $W^{1,p}$ and $f^j$, $j=1,\cdots,n$ be the family of holomorphic disc associated to $\xi_j$ produced by Lemma~\ref{lem:actual}.
We derive an explicit expression of the second variation of $\sum_{j=1}^nJ_{t,\beta}(f^j_\tau)$, which is an extension of \cite[Proposition~3.4]{DGX} and \cite[Proposition~5.2]{OAI26}.

\begin{proposition}\label{prop:variation-sums}
Set
\begin{equation}\label{eq:blaschke-degree-data}
\begin{aligned}
 H:=\Phi^*G_e\Phi,\qquad &S:=\tr H=\sum_jh(\xi_j,\xi_j),
 \qquad
 I:=\int_\D f^*\theta,\quad I_w:=\int_\D w f^*\theta,\\
 &d_\Phi:=\deg(\det\Phi),\qquad
 D_\Phi:=\sum_{b\in\D,\,\det\Phi(b)=0}(1-|b|^2),
 \end{aligned}
\end{equation}
 where the zeros are counted with multiplicity.
Then
\begin{align*}
 \tau_A:&=\sum_j\Delta_{\tau}A(f^j_\tau)|_{\tau=0}=4\pi d_\Phi-2I,\\
 \tau_W:&=\sum_j\Delta_{\tau}W(f^j_\tau)|_{\tau=0}=4\pi n-4\int_\D S\,\dA,\\
 \tau_{Q_\sigma}:&=\sum_j\Delta_{\tau}Q_{\sigma}(f^j_\tau)|_{\tau=0}
   =\frac{1}{2\pi}\int_{S^1}\tr_\omega\sigma(f)\,ds,\\
 \int_\D S\,\dA
   &\geq\pi(n-D_\Phi)+\frac12I_w.
\end{align*}
Here $\Delta_\tau:=4\partial_\tau\partial_{\bar\tau}$ and $\tr_\omega\sigma$ denotes the trace of $h^\sigma$ with respect to
$h$. With $\mu:=t e^{-q(f)}$ one has
\begin{equation}\label{eq:second-variation-sum}
 \sum_j\Delta_{\tau}J_{t,\beta}(f^j_\tau)|_{\tau=0}
   =\tau_A+\mu(c\tau_W-\tau_{Q_\sigma})
          +4\mu\sum_j|\ell_f(\xi_j)|^2
          -\beta\sum_j\Delta_{\tau}\Psi(f^j_\tau)|_{\tau=0}.
\end{equation}
\end{proposition}

\begin{proof}

Put $G:=G_e$, $g:=\log\det G$, $h_j^\sigma:=h^\sigma(\xi_j,\xi_j)$, we first prove $G,g,S\in W^{2,\frac p2}(\D)$ and $\Delta h^\sigma_j\in L^{\frac p2}(\D)$.

We check that $\Delta G,\Delta h^\sigma_j\in L^{p/2}(\D)$.
Let $G_{ij}$ be the entry of Gram matrix $G$, so we have
\[
 G_{ij}(z)=h(f(z))(e_j,e_i).
\]
In a local coordinate chart of $X$, $h$ is expressed as a Hermitian matrix and we set $h_f:=h\circ f$.
Since $e_i,e_j$ are holomorphic, 
\begin{equation}\label{eq:lapG}
 \begin{aligned}
 \partial_z\partial_{\bar z}G_{ij}
 ={}&(\partial_z\partial_{\bar z}h_f)(e_j,e_i)+(\partial_{\bar z}h_f)(\partial_z e_j,e_i)\\
 &+(\partial_zh_f)(e_j,\partial_ze_i)+h_f(\partial_ze_j,\partial_ze_i).
 \end{aligned}
\end{equation}
Since $h$ is smooth on $X$ and $f$ is holomorphic, a similar computation as \eqref{eq:lapG} shows there is a constant $C>0$ so that on $\D$,
\begin{equation}\label{eq:laphf}
 |\partial_zh_f|\le C|df|,
 \qquad
 |\partial_z\partial_{\bar z}h_f|\le C|df|^2.
\end{equation}
Since $e_i,e_j$ is continuous up to the boundary, \eqref{eq:lapG},\eqref{eq:laphf} then give us
\[
 |\Delta G|
 \le C'\left(
 |De|^2+|De|\,|df|+|df|^2
 \right).
\]
As $|De|,|df|\in L^p(\D)$, we obtain
$\Delta G\in L^{\frac p2}(\D)$.

Since $\xi_j$ is holomorphic and of class $W^{1,p}$, the same computations show that $\Delta h^\sigma_j\in L^{\frac p2}(\D)$.

Since $e$ is a holomorphic frame of $f^*TX$ and $G=I_n$ on $S^1$, 
from the definition of $\theta:= \operatorname{Ric}(\omega)=-i\partial\bar\partial\log\det(\omega^n)$ we get
\[
 \Delta g\,\dA
 =2i\partial\bar\partial 
 \log\det G
 =-2f^*\theta,
 \qquad g|_{S^1}=0.
\]
Therefore $
 |\Delta g|\le C|df|^2$ and $\Delta g\in L^{p/2}(\D)$.

Since $G-I_n,g\in W^{1,p}_0(\D)$, and $p/2>1$,
standard regularity theory of elliptic operators
\cite[Theorem~9.15]{GT} yields that
\[
 G,g\in W^{2,p/2}(\D).
\]
Since
$\Phi\in C^\infty(\overline\D)$, it follows that
\[
S=\tr(\Phi^*G\Phi)\in W^{2,p/2}(\D).
\]

\smallskip
\noindent\textbf{Second variations of $A$.}
Differentiation of $g=\log\det G$ gives
\[
 \partial_jg=\tr(G^{-1}\partial_jG),
 \qquad j=x,y.
\]
These derivatives are in $W^{1,p/2}(\D)$ hence the outward normal derivatives $\partial_\nu g,\partial_\nu S \in L^{p/2}(S^1)$ are well-defined in the sense of Sobolev traces.

On $S^1$, we have $G=\Phi\Phi^*=I_n$. Consequently, using $S:=\tr(\Phi^*G\Phi)=\tr(G\Phi\Phi^*)$ we get
\[
 \begin{aligned}
 \partial_\nu g
 &=\tr(G^{-1}\partial_\nu G)
   =\tr(\partial_\nu G),\\
 \partial_\nu S
 &=\tr\left((\partial_\nu G)\Phi\Phi^*
             +G\partial_\nu(\Phi\Phi^*)\right)\\
 &=\partial_\nu g+\partial_\nu\tr(\Phi\Phi^*).
 \end{aligned}
\]
This together with integrating by parts gives
\[
 \begin{aligned}
 \int_\D\Delta S\,\dA
 &=\int_{S^1}\partial_\nu S\,ds\\
 &=\int_{S^1}\partial_\nu\,g\,ds+\int_{S^1}\partial_\nu\tr(\Phi\Phi^*)\,ds\\
 &=\int_\D\Delta g\,\dA
   +\int_\D\Delta\tr(\Phi\Phi^*)\,\dA.
 \end{aligned}
\]
Holomorphic functions $\Phi_{jj}$ are smooth up to the boundary and $|\Phi_{jj}|=1$ on $S^1$, so we have
\[
 \begin{aligned}
 \int_\D\Delta|\Phi_{jj}|^2\,\dA
 =4\pi\deg \Phi_{jj}.
 \end{aligned}
\]
Since $\xi_j=(e\Phi)_j$, we have
$S=\sum_{j=1}^n|\xi_j|_h^2$. 
\eqref{eq:actual-local-variations} 
and $\Delta g\,\dA=-2f^*\theta$, $I=\int_\D f^*\theta$ gives us
\[
 \begin{aligned}
 \tau_A
 :=\sum_{j=1}^n4\partial_\tau\partial_{\bar\tau}|_{\tau=0}
 \Big(\int_\D h(\partial_zf^j,\partial_zf^j)\,dA_z\Big)
 &=\sum_{j=1}^n\int_\D4\partial_z\partial_{\bar z}h(\xi_j,\xi_j)\,dA_z\\
 &=\int_\D\Delta S\,\dA\\
 &=\int_\D\Delta g\,\dA
   +\int_\D\Delta\tr(\Phi\Phi^*)\,\dA\\
 &=-2I+\sum_{j=1}^n\int_\D\Delta|\Phi_{jj}|^2\,\dA\\
 &=-2I+4\pi\sum_{j=1}^n\deg\Phi_{jj}
   =4\pi d_\Phi-2I.
 \end{aligned}
\]
Here $\det\Phi=\prod_{j=1}^n\Phi_{jj}$, so
$d_\Phi=\sum_{j=1}^n\deg\Phi_{jj}$.

\smallskip
\noindent\textbf{Second variations of $W$ and lower bound for $\int_\D S\,\dA$.}

We need the integration by parts formula
\begin{equation}\label{eq:IBP1}
    \int_\D w\Delta v\,\dA
 =-4\int_\D v\,\dA+2\int_{S^1}v\,ds,
 \qquad w(z)=1-|z|^2,\quad v\in W^{2,\frac{p}{2}}(\D).
\end{equation}
Indeed, we know that $C^{\infty}(\overline\D)$ is dense in $W^{2,\frac{p}{2}}(\D)$ and the trace map $W^{1,\frac{p}{2}}(\D)\to L^{\frac{p}{2}}(S^1)$ is continuous. Therefore, we can use smooth approximation to reduce \eqref{eq:IBP1} to the case that $v$ is smooth up to the boundary, where classical Green formula applies.

Since $S$ and $g$ are in $W^{2,\frac{p}{2}}(\D)$ and their boundary traces are $n$ and $0$, respectively.
Commutativity \eqref{eq:actual-local-variations} and integration by parts \eqref{eq:IBP1} give
\begin{align*}
\tau_W:&=\sum_{j=1}^n4\partial_\tau\partial_{\bar\tau}|_{\tau=0}\Big(\int_\D wh(\partial_zf^j,\partial_zf^j)\,dA\Big)
=\int_\D w\Delta S
 =4\pi n-4\int_\D S\,\dA,\\
 I_w:&=\int_\D wf^*\theta=-\frac{1}{2}\int_\D w\Delta g\,dA=2\int_\D g\,\dA.
\end{align*}
We have $S=\tr H\geq n+\log\det H$ by diagonalizing $H$.
$\log\det H=\log\det(\Phi^*G\Phi)=g+\log\det(\Phi^*\Phi)=g+\log|\det\Phi|^2$. So we obtain
\begin{align*}
    \int_\D S\,dA&\geq \pi n+\int_\D g\,dA+\int_\D\log|\det\Phi|^2\,dA\\
    &=\pi n+\frac{1}{2}I_w+\int_\D\log|\det\Phi|^2\,dA.
\end{align*}
If $\Phi=zI_n$ then $\int_\D\log|\det\Phi|^2\,dA=-n\pi$.
If $\phi=\operatorname{diag}(zB_b(z),z,\cdots,z)$ then $\int_\D\log|\det\Phi|^2\,dA=-\pi(n+1-|b|^2)$ using
\[
 \begin{aligned}
 \int_\D\log|B_b(z)|^2\,\dA
 &=4\pi\int_0^1 r\log\max\{r,|b|\}\,dr\\
 &=-\pi(1-|b|^2).
 \end{aligned}
\]
In either case, with $D_\Phi=\sum_{b\in\D,\,\det\Phi(b)=0}(1-|b|^2)$ we get
\[
\int_\D S\,dA\geq\pi(n-D_\Phi)+\frac{1}{2}I_w.
\]

\smallskip
\noindent\textbf{Second variations of $Q_\sigma$.}

Let $v\in W^{1,p}\cap C^{\infty}(\D)$ so that $\Delta v\in L^{\frac p2}(\D)$.
We claim that following integration by parts holds:
\begin{equation}\label{eq:lambda intbp}
     \frac{1}{2\pi}\int_\D\log(1/|z|)\Delta v\,\dA
 =\frac{1}{2\pi}\int_{S^1}v\,ds-v(0).
\end{equation}

Let $\lambda_{\varepsilon}(z):=-\max_{\varepsilon}\{\log|z|,\log\varepsilon\}$ where $\max_\varepsilon$ is the standard regularized maximum function. Then $\lambda_\varepsilon$ is smooth on $\D$ and $\lambda_{\varepsilon}=0,\,\partial_\nu \lambda_{\varepsilon}=-1$ on $S^1$.
Let $\chi_r$ be the cut-off function equal to $1$ on $D_{1-2r}$ and vanishes if $|z|\geq1-r$, $|D\chi_r|\leq \frac{C}{r}$.
\[
\int_\D\chi_r\lambda_\varepsilon \Delta v\,dA
=-\int_\D \chi_rD\lambda_\varepsilon\cdot Dv\,dA
-\int_\D\lambda_\varepsilon D\chi_r\cdot Dv\,dA. 
\]
Using $Dv\in L^p(\D)$, $\lambda_\varepsilon |D\chi_r|\leq C$ if $|z|\geq 1-2r$, we see $\lim_{r\to0}|\int_\D\lambda_\varepsilon D\chi_r\cdot Dv\,dA|\leq \lim_{r\to0}C\int_{\{1-2r\leq |z|<1\}}|Dv|\,dA=0$.
Therefore by letting $r\downarrow0$ we get
\[
\int_\D \lambda_\varepsilon \Delta v\,dA
=-\int_\D  D\lambda_\varepsilon\cdot Dv\,dA.
\]
Since $v\in W^{1,p}(\D)$, using density of $C^{\infty}(\overline\D)$ in $W^{1,p}(\D)$ we get
\[
-\int_\D  D\lambda_\varepsilon\cdot Dv\,dA
=\int_\D v\Delta\lambda_\varepsilon\,dA
-\int_{S^1}v\partial_\nu\lambda_{\varepsilon}\,ds
=\int_\D v\Delta\lambda_\varepsilon\,dA+\int_{S^1}v\,ds.
\]
Finally, using $v$ is smooth and $\Delta \lambda_\varepsilon$ is supported in $\frac12\D$ and weakly converge to $\Delta\lambda=-2\pi\delta_0$, we get the desired formula \eqref{eq:lambda intbp}.

Note that $h^\sigma(\xi_j,\xi_j)$ satisfies these hypotheses. Together with commutativity \eqref{eq:actual-local-variations}, this integration by parts formula gives
\begin{align*}
\tau_{Q_\sigma}:=\frac{1}{2\pi}\sum_{j=1}^n4\partial_\tau\partial_{\bar\tau}|_{\tau=0}\Big(\int_\D \log(1/|z|)h^\sigma(\partial_zf^j,\partial_zf^j)\,dA\Big)
&=\frac{1}{2\pi}\int_\D\log(1/|z|)\Delta_z\Big(\sum_{j=1}^nh^\sigma(\xi_j,\xi_j)\Big)\,dA\\
&=\frac{1}{2\pi}\int_{S^1}{h^\sigma}(\xi_j,\xi_j)\,ds
  -{h^\sigma}(\xi_j,\xi_j)(0).
\end{align*}
Since $\xi_j(0)=0$ and $\xi_j$ is boundary-unitary, we get
\[
\tau_{Q_\sigma}=\frac{1}{2\pi}\int_{S^1}\tr_\omega\sigma(f)\,ds.
\]

\smallskip
\noindent\textbf{Second variation of $J_{t,\beta}$.}
Chain rule gives, for $\mu:=te^{-q(f)}$, $\ell_f(\xi_j):=\partial_\tau q(f^j_\tau)|_{\tau=0}$,
\[
 \left.\Delta_\tau(t e^{-q(f^j_\tau)})\right|_{\tau=0}
 =-\mu\Delta_\tau q(f^j_\tau)\Big|_{\tau=0}
  +4\mu|\ell_f(\xi_j)|^2.
\]
Since $q=Q_\sigma-cW$, \eqref{eq:second-variation-sum} follows from the previous computations.
\end{proof}

We prove that the variation fields $\xi=(\xi_1,\cdots,\xi_n)$ can be chosen so that $\ell_f(\xi_j)=0$.

\begin{lemma}
\label{lem:exact-kernel}
Let $f\in W^{1,p}(\D,X)$ be holomorphic for some $p>2$, and let $e=(e_1,\ldots,e_n)$ be a boundary-unitary frame from
Lemma~\ref{lem:baseline-frame}.
There are a unitary matrix $U\in U(n)$ and a holomorphic matrix $\Phi$ having the form
\[
 \Phi(z)=zI_n
 \quad\text{or}\quad
 \Phi(z)=\operatorname{diag}(zB_b(z),z,\ldots,z),
 \qquad b\in\D,
\]
such that, writing $\xi:=(\xi_1,\cdots,\xi_n):=eU\Phi$, one has
\[
 \xi_j\in W^{1,p}(\D,f^*T^{1,0}X),\qquad
 \xi_j(0)=0,\qquad \ell_f(\xi_j)=0
 \quad(1\le j\le n),
\]
and
\[
 h(\xi_j,\xi_k)|_{S^1}=\delta_{jk}
 \quad(1\le j,k\le n).
\]
\end{lemma}

\begin{proof}
Let $\xi\in W^{1,p}(\D,f^*T_X)$ be a holomorphic section vanishing at $0$, let $F:D_\varepsilon\times\D\to X$ be the holomorphic map produced by Lemma~\ref{lem:actual}, it satisfies $F(0,\cdot)=f$, $\partial_\tau F|_{\tau=0}=\xi$.
Lemma~\ref{lem:first-variation-q} then shows 
\begin{equation}\label{eq:first var of q}
\ell_f(\xi):=\partial_\tau q(f_\tau)|_{\tau=0}=\frac{1}{4\pi}\int_\D\frac{h^\sigma(\xi,\partial_zf)}{z}-c\int_Dzh(\xi,\partial_zf).
\end{equation}
This in particular shows that $\ell_f(\xi)$ is linear in $\xi$.

Put
\[
 a_j:=\ell_f(ze_j),\qquad
 a:=(a_1,\ldots,a_n),\qquad
 |a|:=\left(\sum_{j=1}^n|a_j|^2\right)^{1/2},
\]
If $a=0$ then we are done, taking $U=I_n$, $e'=e$, and $\Phi=zI_n$.
Suppose $a\ne0$.
Complete the unit vector
$|a|^{-1}(\bar a_1,\ldots,\bar a_n)^T$
to the column of a unitary matrix $U$.
Then
\[
 aU=(|a|,0,\ldots,0).
\]
Set $e':=eU$, we have
\begin{equation}\label{eq:kernel-unitary-reduction}
 \ell_f(ze'_j)
 =\sum_{k=1}^n U_{kj}\ell_f(ze_k)
 =\sum_{k=1}^n U_{kj}a_k
 =
 \begin{cases}
  |a|,&j=1,\\
  0,&2\le j\le n.
 \end{cases}
\end{equation}

We only need to choose $B_b(z)$ for some $b\in\D$ so that $\ell_f(zB_b(z)e_1')=0$.
Since $B_b$ is holomorphic on a neighbourhood of $\overline\D$, $zB_b(z)e'_1\in W^{1,p}(\D,f^*TX)$ and vanishes at $0$. \eqref{eq:first var of q} then gives
\begin{equation}\label{eq:kernel-scalar-functional}
 \begin{aligned}
 \ell_f(zB_b(z)e'_1)
 =\frac1{4\pi}\int_\D B_b(z)\,h^\sigma(e'_1,\partial_zf)\,dA
       -c\int_\D |z|^2B_b(z)\,h(e'_1,\partial_zf)\,dA.
 \end{aligned}
\end{equation}
We define $\Gamma:\overline\D\to\mathbb C$ by
\[
 \Gamma(b):=
 \begin{cases}
  \ell_f(zB_be_1'),&|b|<1,\\
  -b|a|,&|b|=1.
 \end{cases}
\]
We check $\Gamma$ is continuous on $\overline\D$.
Using $f$ and $e_1'$ are of class $W^{1,p}$ on $\D$, we get $h^\sigma(e_1',\partial_zf), h(e_1',\partial_zf)\in L^p(\D)$. 
It is also clear that $|B_b(z)|\leq1$, $\forall z,b\in\D$.
If sequence $b_\nu\to b\in\D$, then
$B_{b_\nu}(z)\to B_b(z)$ for every $z\in\D$. Dominated convergence theorem proves $\lim \Gamma(b_\nu)=\Gamma(b)$.
If sequence $b_\nu\to b\in S^1$, then
\[
 B_{b_\nu}(z)\longrightarrow
       \frac{z-b}{1-\bar b z}=-b,
 \qquad z\in\D.
\]
Then dominated convergence gives
\begin{align*}
    \lim \Gamma(b_\nu)&=\lim\Big(\frac1{4\pi}\int_\D B_{b_\nu}(z)\,h^\sigma(e'_1,\partial_zf)\,dA
       -c\int_\D |z|^2B_{b_\nu}(z)\,h(e'_1,\partial_zf)\,dA\Big)\\
       &=(-b)\Big(\frac{1}{4\pi}\int_\D \,h^\sigma(e'_1,\partial_zf)\,dA
       -c\int_\D |z|^2B_b(z)\,h(e'_1,\partial_zf)\,dA\Big)=-b|a|.
\end{align*}
The continuous function $\frac{\Gamma}{-|a|}$ is identity on $S^1$, fixed-point theorem proves $\Gamma$ has zero on $\D$: there is $b\in\D$ such that
\[
 \Gamma(b)=\frac1{2\pi}\int_\D B_{b}(z)\,h^\sigma(e'_1,\partial_zf)\,dA
       -c\int_\D |z|^2B_{b}(z)\,h(e'_1,\partial_zf)\,dA=0.
\]
Take
\[
 \Phi=\operatorname{diag}(zB_{b},z,\ldots,z),
\]
$\xi:=(eU\Phi)$ is what we want.
\end{proof}

We combine Proposition~\ref{prop:variation-sums} with the estimate of the second variation of $\Psi(f_\tau)$ in Proposition~\ref{prop:levi} to get a uniform bound of weighted area of minimizer of $J_{t,\beta}$.

\begin{proposition}
\label{prop:uniform-variation-bound}
Suppose $f\in W^{1,p}(\D,X)$ is holomorphic, $p>2$, and
\[
 I(f):=\int_\D f^*\theta\ge\delta A(f),\qquad I_w(f):=\int_\D(1-|z|^2)f^*\theta \ge\delta W(f).
\]
Fix $0<\beta\le\delta/(8\pi n)$ and put
$\mu:=te^{-q(f)}$, set
\[
 S_0:=\max_X(-\operatorname{tr}_\omega\sigma)_+,\qquad
 C_0:=4\pi(n+1),\qquad
 C_1:=4\pi c(n+1)+S_0.
\]
Let $\xi=(\xi_j)$ be the sections obtained in Lemma~\ref{lem:exact-kernel}, then
\begin{align}
 \Delta_\tau\sum_{j=1}^n{J_{t,\beta}}(f^j_\tau)|_{\tau=0}
 &\le C_0-(2\delta-8\pi n\beta)A(f)
       +\mu\bigl(C_1-2c\delta W(f)\bigr)
                                      \label{eq:sharp-variation-bound}\\
 &\le C_0-\delta A(f)+\mu\bigl(C_1-2c\delta W(f)\bigr)
                                      \label{eq:variation-bound}\\
 &\le C_0+C_1\mu-\delta(1+2c\mu)W(f).
                                      \label{eq:weighted-area-variation-bound}
\end{align}
Consequently every global minimum $f_{t,\beta}$ from
Proposition~\ref{prop:foundation-attainment} satisfies
\begin{equation}\label{eq:uniform-W-bound}
 W(f_{t,\beta})\le
 W_{\max}:=\frac1\delta\max\{C_0,C_1/(2c)\}.
\end{equation}
\end{proposition}

\begin{proof}
Since $\ell_f(\xi_j)=0$, $\forall j=1,\cdots,n$, Proposition~\ref{prop:variation-sums} gives
\begin{align*}
\Delta_\tau\sum_{j=1}^n{J_{t,\beta}}(f^j_\tau)|_{\tau=0}
=\tau_A+\mu(c\tau_W-\tau_{Q_\sigma})-\beta\Delta_\tau\Big(\sum_{j=1}^n\Psi(f^j_\tau)\Big)|_{\tau=0}.
\end{align*}
We have
\[
\tau_A=4\pi d_\Phi-2I(f)\leq 4\pi(n+1)-2I(f).
\]
The lower bound for $\int_\D S\,dA$ gives
\[
 \tau_W=4\pi n-4\int_\D S\,dA
       \le4\pi D_\Phi-2I_w(f)\leq 4\pi(n+1)-2I_w(f).
\]
The formula $\tau_{Q_\sigma}=\frac{1}{2\pi}\int_{S^1}\tr_{\omega}\sigma$ gives $-\tau_{Q_\sigma}\le S_0$.
Proposition~\ref{prop:levi} yields the estimate 
\[
\Delta_\tau\Psi(f^j_\tau)|_{\tau=0}\geq-8\pi A(f)\max_{S^1}h(\xi_j,\xi_j)=-8\pi A(f),
\]
the equality uses the fact that $\xi_j$ is boundary-unitary.
Together with the two assumed inequalities $I(f)\geq\delta A(f), I_w(f)\geq\delta W(f)$ therefore imply
\[
\Delta_\tau\sum_{j=1}^n{J_{t,\beta}}(f^j_\tau)|_{\tau=0}
 \le C_0-(2\delta-8\pi\beta n)A(f)
             +\mu(C_1-2c\delta W(f)).
\]
This is \eqref{eq:sharp-variation-bound}.
The bound on $\beta$ is equivalent to $2\delta-8\pi\beta n\ge\delta$,
which proves \eqref{eq:variation-bound}.
Since $A\ge W$, it also proves
\eqref{eq:weighted-area-variation-bound}.

Let $f_{t,\beta}$ be a minimizer obtained in
Proposition~\ref{prop:foundation-attainment}, then $f_{t,\beta}\in W^{1,p}(\D,X)$
for every $2<p<4$ and automatically satisfies
\[
 I(f_{t,\beta})\ge\delta A(f_{t,\beta}),
 \qquad
 I_w(f_{t,\beta})\ge\delta W(f_{t,\beta}).
\]
Denote the sections obtained in Lemma~\ref{lem:exact-kernel} by $\xi_1,\ldots,\xi_n$.
For each $j$, the holomorphic
family $f^{(j)}_\tau$ satisfies
\[
 f^{(j)}_0=f_{t,\beta},\qquad
 \left.\partial_\tau f^{(j)}_\tau\right|_{\tau=0}=\xi_j,
 \qquad
 f^{(j)}_\tau(0)=x.
\]
Each $f^{(j)}_\tau$ belongs to $W^{1,p}(\D,X)$ and hence,
by Lemma~\ref{lem:usc}(1) is a competitor of $\inf_{f\in \mathcal B_x} J_{t,\beta}(f)$. Therefore,
\[
 J_{t,\beta}(f^{(j)}_\tau)
 \geq
 J_{t,\beta}(f^{(j)}_0),
\]
and hence
$\Delta_\tau\varphi_j(0)\ge0$. \eqref{eq:weighted-area-variation-bound} gives us
\[
 0\leq  \Delta_\tau\sum_{j=1}^n{J_{t,\beta}}(f^j_\tau)|_{\tau=0}\leq C_0+C_1\mu-\delta(1+2c\mu)W(f_{t,\beta})
\]
which says
\[
 W(f_{t,\beta})\le\frac{C_0+C_1\mu}{\delta(1+2c\mu)}
   \le\frac1\delta\max\{C_0,C_1/(2c)\},
\]
and completes the proof.
\end{proof}
We will choose appropriate $t_j,\beta_j$ so that $\mu_j=t_je^{-q(f)},\beta_j\to0$, therefore we see that  only the second variation of area functional contributes to the optimal constant $2\pi(n+1)$ in Theorem~\ref{thm:pointwise-criterion}.

We now assume \eqref{eq:pointwise-escape}. Radial scaling first
turns the uniform bound of weighted area $W(f)$ into uniform bound of area $A(f)$.
with $q\to+\infty$.
We then select minimizers for which $te^{-q}$ and $-\beta\Psi$ tend to zero,
so that \eqref{eq:sharp-variation-bound} keeps its sharp
area coefficient $4\pi(n+1)$.

\begin{proposition}
\label{lem:vanishing-correction-minima}
Assume \eqref{eq:pointwise-escape}. There are $t_j\to\infty$,
$\beta_j\downarrow0$, and global minima
$f_j=f_{t_j,\beta_j}\in\mathcal B_x$ such that, with
$\mu_j:=t_j e^{-q(f_j)}$,
\begin{equation}\label{eq:vanishing-correction-limits}
 \mu_j\to0,\qquad -\beta_j\Psi(f_j)\to0,\qquad
 q(f_j)\to+\infty,\qquad A(f_j)\to L,
\end{equation}
where
\begin{equation}\label{eq:sharp-limiting-area}
 0\le L\le\frac{2\pi(n+1)}{\delta}.
\end{equation}
\end{proposition}
\begin{proof}
We first fix $\beta_*:=\delta/(16\pi n)$ and take minimizers
$f_{t,\beta_*}$ of  $J_{t,\beta_*}$ obtained in Proposition~\ref{prop:foundation-attainment}.
By Proposition~\ref{prop:uniform-variation-bound}, we have
$W(f_{t,\beta_*})\le W_{\max}$. For any fixed $f\in\mathcal A_x$,
Corollary~\ref{cor:ordinary-competitors} gives
$J_{t,\beta_*}(f_{t,\beta_*})\le J_{t,\beta_*}(f)$.
Since $te^{-q(f_{t,\beta_*})}\le J_{t,\beta_*}(f_{t,\beta_*})$,
we obtain
\[
 \limsup_{t\to\infty}e^{-q(f_{t,\beta_*})}
 \le \lim_{t\to\infty}\frac{J_{t,\beta_*}(f)}{t}
 =e^{-q(f)}.
\]
It follows that
$q(f_{t,\beta_*})\to+\infty$ since we can let $q(f)\to\infty$ thanks to \eqref{eq:pointwise-escape}. 
Consider radial scaling $(f_{t,\beta_*})_r(z):=f_{t,\beta_*}(rz)$ and
apply Lemma~\ref{lem:radial-ordinary} with $r=1/2$, we get discs $(f_{t,\beta_*})_{1/2}\in\mathcal A_x$ and the estimates
\[
 A((f_{t,\beta_*})_{1/2})
 \le \frac43 W(f_{t,\beta_*})
 \le \frac43 W_{\max},
 \qquad
 q((f_{t,\beta_*})_{1/2})\longrightarrow+\infty
 \quad(t\to\infty).
\]

Set
\[
 L:=\lim_{s\to\infty}
 \inf\{A(f): f\in\mathcal A_x,\ q(f)>s\}.
\]
The above discs $(f_{t,\beta_*})_{1/2}$ shoes $L\le\frac43W_{\max}$. 
By definition of $L$ we can choose $g_j\in\mathcal A_x$ with
$q(g_j)>2j$ and $A(g_j)\le L+j^{-1}$. We then take $t_j=e^j$ and
choose $0<\beta_j\le\beta_*$, $\lim_j\beta_j=0$ so that
$-\beta_j\Psi(g_j)\le j^{-1}$, which is possible since $\Psi(g_j)$ is finite for each $j$ (see e.g. Lemma~\ref{lem:usc}(1)). Now, take a minimizer $f_j:=f_{t_j,\beta_j}$
from Proposition~\ref{prop:foundation-attainment} and write
$\mu_j=t_je^{-q(f_j)}$. By Corollary~\ref{cor:ordinary-competitors}, we deduce that
\[
 A(f_j)+\mu_j-\beta_j\Psi(f_j)
 =J_{t_j,\beta_j}(f_j)
 \le J_{t_j,\beta_j}(g_j)
 \le L+2j^{-1}+e^{-j}.
\]
In particular, $\mu_j=t_je^{-q(f_j)}$ is bounded, thus $q(f_j)\to+\infty$.
By Lemma~\ref{lem:radial-ordinary}, we can choose $r_j<1$
with $q((f_j)_{r_j})\ge q(f_j)-1$. Since
$(f_j)_{r_j}\in\mathcal A_x$ and $A((f_j)_{r_j})\le A(f_j)$,
the definition of $L$ gives $\liminf_j A(f_j)\ge L$, which forces $A(f_j)\to L$ and
$\mu_j\to0,\,-\beta_j\Psi(f_j)\to0$.

Finally, variation estimates \eqref{eq:sharp-variation-bound} in Proposition~\ref{prop:uniform-variation-bound} gives
\[
 0
 \le 4\pi(n+1)-(2\delta-8\pi n\beta_j)A(f_j)+C_1\mu_j.
\]
Passing to the limit yields $2\delta L\le4\pi(n+1)$.
\end{proof}

\begin{proof}[Proof of Theorem~\ref{thm:pointwise-criterion}]
Proposition~\ref{lem:vanishing-correction-minima} produces holomorphic discs $f_j$ with uniformly bounded area centered at $x$, with
\[
 q(f_j)\to+\infty,\qquad
 A(f_j)\to L\le\frac{2\pi(n+1)}{\delta}.
\]
Corollary~\ref{cor:quantitative-marked-bubbling} gives a rational curve $R$ through $x$, with $\omega$-area $[\omega]\cdot R:=\int_R\omega\leq L\leq\frac{2\pi(n+1)}{\delta}$. The proof is complete.
\end{proof}

\vspace{0.2cm}
\section{Disc functional, pseudoeffectivity, and uniruledness of compact K\"ahler spaces}
\label{sec:nonpsef}
In this section we prove Theorem~\ref{thm:disc functional} by combining Poletsky's graph theorem \cite{PoletskyGraphs} and an argument of Magn\'usson's \cite{Mag11,Magnusson} using Poisson-envelope theorem
\cite{Ros03}. 
As an application, we give an alternative proof of Ou's uniruledness theorem \cite{Ou}.

\begin{theorem}[= Theorem~\ref{thm:disc functional}]
\label{lem:smooth-disc-envelope}
Let $X$ be a connected complex manifold and $\rho$ be a smooth closed real $(1,1)$-form.
Then for every $x\in X$,
\[
 x\longrightarrow  E(x):= \inf_{f\in\mathcal A(X,x)}\frac{1}{2\pi}\int_\D\log(1/|z|)f^*\rho
\]
is either $-\infty$ or $\rho$-plurisubharmonic.
\end{theorem}
\begin{proof}
We first prove $E$ is upper semicontinuous.
The idea is that, for $x_0\in X$ and $f:D_R\to X$ be a holomorphic disc with $f(0)=x_0$, we construct a family of discs parametrized by a neighborhood $U$ of $x_0$ so that
\[
 F:D_s\times U\longrightarrow X,
 \qquad F(\cdot,x_0)=f,\qquad F(0,x)=x
\]
for some $1<s<R$.

Forstneri\v c's graph theorem plays a key role.
There are some $1<s<R$ and neighborhood $\Omega\subset \mathbb C\times X$ of graph $G=\{(z,f(z)):|z|\leq s\}$ and an fiber-preserving isomorphism
\[
\Theta:\Omega\longrightarrow \Omega'\subset D_R\times \mathbb{C}^n
\]
We can write $\Theta(z,x)=(z,G(z,x))$, $\Theta^{-1}(z,v)=(z,H(z,v))$, then the set $\{(z,G(z,f(z))):|z|\le s\}$ is compact in $\Omega'$,
so there is $\eta>0$ such that
\[
 \big(z,G(z,f(z))+b\big)\in\Omega_f'
 \qquad(|z|\le s,\ |b|\leq \eta).
\]
Choose a neighborhood $U$ of $f(0)$ such that
$(0,x)\in\Omega$ and $|G(0,x)-G(0,f(0))|<\eta$ if $x\in U$.
Then define
\[
 F(z,x)=H\bigl(z,G(z,f(z))+G(0,x)-G(0,f(0))\bigr),
\]
we obtain $F(z,f(0))=H(z,G(z,f(z)))=f(z)$ and $F(0,x)=H(0,G(0,x))=x$.

Using $F(\cdot,x)\to F(\cdot,x_0)$ in $C^1(\overline\D)$ as $x\to x_0$, we then obtain for $x_0\in U$,
\[
 \left|\int_\D\log(1/|z|)F^*(\cdot,x)\rho-\int_\D\log(1/|z|)F(\cdot,x_0)^*\rho\right|
 \longrightarrow0,\qquad x\to x_0.
\]
Since $f=F(\cdot,x_0)$ is arbitrary, we get $E$ is upper semicontinuous.

\vspace{0.15cm}
We prove the sub-mean value inequality for $E$. The key idea is to use Poletsky--Rosay theorem of the Poisson envelope \cite{Pol91,Pol93,Ros03}. 
Let $U\subset X$ be a small ball centered at $x_0$ and $dd^c\psi=\rho$ for some smooth potential $\psi$. 
Let $h:D_r\hookrightarrow U$ be a holomorphic disc with $r>1$ centered at $x_0$. We need to prove
\begin{equation}\label{eq:disc-envelope-submean}
 E(x_0)+\psi(x_0)
 \le\frac{1}{2\pi}\int_{S^1}(E +\psi )(h(e^{i\theta}))\,d\theta.
\end{equation}
We may assume $E(x_0)>-\infty$, fix $\varepsilon>0$ and a continuous function
$v:S^1\to\mathbb R$ with $v\ge E+\psi$ on $S^1$.
At each $t_0\in S^1$ choose
a disc $f_0:D_R\to X$ with $f_0(0)=h(t_0)$ with
$Q_\rho(f_0)+\psi(h(t_0))<v(t_0)+\varepsilon/2$,
and we construct a family $F:D_s\times U_{t_0}\to X$, $s>1$ as before. For an open set $B_0\subset\subset D_r$ centered at $t_0$, we define $F_0:D_s\times B_0\to X,\,\,(z,t)\mapsto F(z,h(t))$.
Continuity then gives, up to shrinking $B_0$,
\[
Q_\rho(F_0(\cdot,t))+\psi(h(t))\leq v(t)+\varepsilon,\quad  \text{$t\in S^1\cap \overline{B_0}$}.
\]
We can choose finitely many such
holomorphic maps $F_j:D_s\times B_j\to X$, $s>1$ satisfying
$B_j\subset\subset D_r$ are pairwise disjoint, $J_j\subset S^1\cap B_j$ are compact arcs and
\begin{equation}\label{eq:disc-envelope-arc-choice}
 \begin{gathered}
 F_j(0,t)=h(t),\qquad
 Q_\rho(F_j(\cdot,t))+\psi(h(t))\le v(t)+\varepsilon
       \quad \text{if}\,\,\, t\in J_j,\\
 \frac{1}{2\pi}\int_{S^1\backslash\cup_j J_j}\bigl(|v|+|\psi\circ h|\bigr)(e^{i\theta})\,d\theta<\varepsilon.
 \end{gathered}
\end{equation}
However, we cannot take \emph{a single potential} of $\rho$ on these discs, since they may not lie in the small ball $U$.
To overcome this difficulty, on $\mathbb C^2\times X$ we consider the graph of these discs
\[
 \begin{aligned}
 K_0&:=\{(0,t,h(t)):t\in\overline\D\},\\
 K_j&:=\{(z,t,F_j(z,t)):
                  t\in J_j,\ z\in\overline\D\},\qquad
 K:=\bigcup_{j=0}^mK_j.
 \end{aligned}
\]
Let $\pi:\mathbb C^2\times X\to X$ be the projection.
Using Lemma~\ref{lem:disc-graph-potential} proved below,
there is a connected Stein neighborhood $V\supset K$ so that $\pi^*\rho$ has potential on $V$, namely
\begin{equation}\label{eq:disc-envelope-global-potential}
 dd^c\Psi=\pi^*\rho.
\end{equation}

Define on $V$
\[
 \widehat E(y):=
 \inf_{g\in\mathcal A(V,y)}Q_{\pi^*\rho}(g):=\inf_{g\in\mathcal A(V,y)}\frac{1}{2\pi}\int_\D\log(1/|z|)g^*(\pi^*\rho).
\]
It is clear that
\begin{equation}\label{eq:disc-envelope-lift-bound}
 E\circ\pi\le\widehat E\le0.
\end{equation}
Green's formula gives, for every $g\in\mathcal A(V,y)$,
\[
 Q_{\pi^*\rho}(g)
 =\frac{1}{2\pi}\int_{S^1}\Psi\circ g(e^{i\theta})\,d\theta-\Psi(y).
\]
Hence $\widehat E+\Psi=\inf_{g\in \mathcal{A}(V,y)}\frac{1}{2\pi}\int_{S^1}\Psi\circ g(e^{i\theta}) d\theta$, it is exactly the Poisson envelope of function
$\Psi$. Since $V$ is connected, Poletsky--Rosay theorem
\cite{Ros03} says
it is plurisubharmonic or identically $-\infty$.
For $\widehat h(t):=(0,t,h(t))$, inequality
\eqref{eq:disc-envelope-lift-bound} and $E(h(0))> -\infty$ ensure $(\hat{E}+\Psi)(\hat{h}(0))>-\infty$, hence $\widehat{E}+\Psi$ is plurisubharmonic.
Moreover, 
\[
 dd^c\bigl(\psi\circ h-\Psi\circ\widehat h\bigr)
 =h^*\rho-\widehat h^*\pi^*\rho=0.
\]
Since $\widehat{h}(\overline\D)\subset V$, We get $\widehat E\circ\hat h+\psi\circ h$ is subharmonic in a neighborhood of $\overline\D$.

For $t\in J_j$, the holomorphic disc
$\widehat{h}_j(z):=(z,t,F_j(z,t))$ lies in $V$ and is centered at
$(0,t,F_j(0,t))=(0,t,h(t))=\widehat h(t)$. Therefore from definition of $\widehat E$ we get
\[
 \widehat E(\widehat h(t))\le \frac{1}{2\pi}\int_\D\log(1/|z|)\widehat{h}_j^*(\pi^*\rho)
 =\frac{1}{2\pi}\int_\D\log(1/|z|)(F_j(\cdot,t))^*\rho
 =
 Q_\rho(F_j(\cdot,t)),
 \quad t\in J_j.
\]
The sub-mean value inequality for
$\widehat E\circ\widehat h+\psi\circ h$ and \eqref{eq:disc-envelope-lift-bound}, \eqref{eq:disc-envelope-arc-choice} now yield
\[
 \begin{aligned}
 E(h(0))+\psi(h(0))
 &\le\widehat E(\widehat h(0))+\psi(h(0))\\
 &\le\frac{1}{2\pi}\int_{S^1}
        (\widehat E\circ\widehat h+\psi\circ h)(e^{i\theta})\,d\theta\\
 &\le\sum_j\frac{1}{2\pi}\int_{J_j}(v+\varepsilon)(e^{i\theta})\,d\theta
       +\frac{1}{2\pi}\int_{S^1\backslash\cup_{j} J_j}\psi\circ h\,d\theta\\
 &\le\frac{1}{2\pi}\int_{S^1}v\,d\theta+2\varepsilon.
 \end{aligned}
\]
Since $\varepsilon>0$ and continuous function $v\geq (E+\psi)|_{S^1}$ are arbitrary, \eqref{eq:disc-envelope-submean} is proved.
We therefore obtain that $E$ is $\rho$-plurisubharmonic or identically $-\infty$.
\end{proof}


We complete the proof of \eqref{eq:disc-envelope-global-potential} by establishing the statement below.
\begin{lemma}
\label{lem:disc-graph-potential}
Let $r>1$ and $h:D_r\to X$ be a holomorphic disc.
Let $J_1,\ldots,J_m\subset S^1$ be pairwise disjoint compact arcs, and let $B_j\subset\subset D_r$ be pairwise disjoint
open neighborhoods of $J_j$.
Suppose $s>1$, $F_j:D_s\times B_j\to X$ are holomorphic and
$F_j(0,t)=h(t)$ for $t\in B_j$.
In $\mathbb C^2\times X$, set
\[
 \begin{aligned}
 K_0&:=\{(0,t,h(t)):t\in\overline\D\},\\
 K_j&:=\{(z,t,F_j(z,t)): z\in\overline\D,\,\,t\in J_j\},
 \qquad K:=\bigcup_{j=0}^mK_j.
 \end{aligned}
\]
For every smooth closed real $(1,1)$-form $\alpha$ on $\mathbb{C}^2\times X$,
there are a connected Stein neighborhood $V\supset K$ and
$\Psi\in C^\infty(V,\mathbb R)$ such that on $V$,
\[
 dd^c\Psi=\alpha.
\]
\end{lemma}

\begin{proof}
We explain the idea of the proof. $J:=\cup J_j\neq S^1$ ensures $L:=(\{0\}\times\overline\D)\cup (\overline\D\times J)$ is polynomial convex, hence has Stein neighborhood basis.
The continuous map $G:L\to X$ obtained by gluing $h,F_1,\cdots,F_m$ is locally holomorphically approximable in the sense of \cite[Theorem~3.1]{PoletskyGraphs}, hence $K$, the graph of $G$, has Stein neighborhood basis in $\mathbb{C}^2\times X$.
It is easy to see $K$ has a contractible neighborhood, so $\rho=d\eta$ on this neighborhood. Finally pick a smaller Stein neighborhood of $K$ we use standard decomposition $\eta=\eta^{1,0}+\eta^{0,1}$ and solve $\eta^{0,1}=\bar{\partial}a$, $\eta^{1,0}=\partial \bar{a}$. Then $dd^c(\operatorname{Im}(a))=\partial\eta^{0,1}+\bar\partial\eta^{1,0}=d\eta=\rho$.

We give the details.
Put $J:=\bigcup_jJ_j$ and
\[
 L:=(\{0\}\times\overline\D)\cup (\overline\D\times J)
 \subset\mathbb C^2.
\]
Recall that the polynomial hull of $L$ is
\[
 \widehat L
 :=\left\{(z_0,t_0)\in\mathbb C^2:
 |Q(z_0,t_0)|\le\max_{(z,t)\in L}|Q(z,t)|
 \text{ for every }Q\in\mathbb C[z,t]\right\}.
\]
We prove that $L$ is polynomial convex, namely $L=\widehat L$.
Fix $(z_0,t_0)\notin L$. If $|z_0|>1$ or
$|t_0|>1$, the coordinate $z$ or $t$, respectively,
provides the required separation. 
We may assume that
$|z_0|\leq1$ and $|t_0|\leq1$. Since $(z_0,t_0)\notin L$, this means $z_0\neq0$ and $t_0\notin J$. If $|t_0|<1$, choose any
$a\in S^1\setminus J$; if $|t_0|=1$, choose $a=t_0$. For
$0<\lambda<1$, we consider
\[
q_\lambda(t)=\frac{1+\lambda\bar a t}{1-\lambda\bar a t}.
\]
Then
\[
\operatorname{Re}q_\lambda(t)
=\frac{1-\lambda^2|t|^2}{|1-\lambda\bar a t|^2}.
\]
Since $a\notin J$ and $J$ is compact, it follows that
\[
\max_{t\in J}\operatorname{Re}q_\lambda(t)\longrightarrow0
\qquad\text{as }\lambda\uparrow1.
\]
On the other hand, if $|t_0|<1$, then
\[
\operatorname{Re}q_\lambda(t_0)
\longrightarrow
\frac{1-|t_0|^2}{|1-\bar a t_0|^2}>0,\qquad\text{as }\lambda\uparrow1
\]
if $t_0=a\in S^1$, then
\[
q_\lambda(t_0)=\frac{1+\lambda}{1-\lambda}
\longrightarrow+\infty\qquad\text{as }\lambda\uparrow1.
\]
We can thus fix $\lambda<1$ sufficiently close to $1$ so that
\[
c:=\operatorname{Re}q_\lambda(t_0)
>
b:=\max_{t\in J}\operatorname{Re}q_\lambda(t).
\]
Choose an integer $N\geq1$ such that $|z_0|e^{N(c-b)}>1$, and put
\[
H(z,t)=z\,e^{Nq_\lambda(t)}.
\]
The function $H$ vanishes on $\{0\}\times\overline{\mathbb D}$,
and hence
\[
|H(z_0,t_0)|=|z_0|e^{Nc}
>
e^{Nb}
=\max_L|H|.
\]
The function $t\mapsto e^{Nq_\lambda(t)}$ is holomorphic on
$D_{1/\lambda}$. Its Taylor
polynomials $P_k$ therefore converge uniformly on
$\overline{\mathbb D}$. Consequently, for sufficiently large $k$ we still have
\[
|z_0P_k(t_0)|>\max_L|zP_k|.
\]
This proves $\widehat L=L$. 
We remark that $J\neq S^1$ is necessary, otherwise the polynomial hull is obviously $\overline\D^2$. 


By e.g. \cite[Lemma~2.7.4]{Hormander}, 
$L$ has a neighborhood basis consisting of polynomial
polyhedron. In particular $L$ has a basis
of open Stein neighborhoods.

Define $G:L\to X$ by
\[
 G(0,t):=h(t),\qquad G(z,t):=F_j(z,t)\quad(t\in J_j).
\]
The identities $F_j(0,t)=h(t)$ make $G$ well defined and
continuous. $K\subset\mathbb C^2\times X$ is its graph.
Our assumption $F_j(0,t)=h(t)$ for $t\in B_j$ yields that near a point $(z_0,t_0)\in \overline{\D}\times J_j$, the map $G$ extends
holomorphically to $(z,t)\mapsto F_j(z,t)$.
Remaining points of $L$ satisfies $(0,t_0)\in \{0\}\times(\overline\D\,\backslash\, J)$, $G$ extends to
$(z,t)\mapsto h(t)$.
Thus $G$ is locally holomorphic, in particular is locally holomorphically approximable.
Poletsky's graph theorem \cite[Theorem~3.1]{PoletskyGraphs}
therefore gives a basis of open Stein neighbourhoods of
$K$ in $\mathbb{C}^2\times X$.

Choose $\chi_j\in C_c^\infty(B_j,[0,1])$ with
$\chi_j=1$ near $J_j$.
We define
\[
 \widetilde G(z,t):=
 \begin{cases}
 F_j(\chi_j(t)z,t),&t\in B_j,\\
 h(t),&t\notin\bigcup_jB_j
 \end{cases}
 \qquad (z,t)\in  D_s\times D_r,
\]
which a smooth map agreeing with $G$ on $L$. The graph
\[
 \Gamma:=
 \{(z,t,\widetilde G(z,t)):(z,t)\in  D_s\times D_r\}
\]
is a smooth embedded submanifold of $ D_s\times D_r\times X$ containing $K$.
The graph $\Gamma$ is diffeomorphic to
$D_r\times D_s$, which is contractible.
$\Gamma$ has a tubular neighborhood $\Omega\subset D_s\times D_r\times X$ so that $\Gamma$ is a strong deformation retract of \(\Omega\) (see e.g. \cite[Theorem~5.25]{LeeRiemannian}).
Hence $\Omega$ is contractible.

Since $d\alpha=0$, Poincare's lemma yields that there is a smooth real $1$-form $\eta$
on $\Omega$ such that
\[
 d\eta=\alpha.
\]

Finally we choose a Stein neighborhood $V$ of $K$ so that
$V\subset\Omega$.
Since $K$ is connected, if $V$ is not connected, we can simply take the connected component containing $K$.
Decompose $\eta=\eta^{1,0}+\eta^{0,1}$ on $V$.
The $(0,2)$-part of $d\eta=\alpha$ vanishes, so we get
\[
 \bar\partial\eta^{0,1}=0.
\]
Since $V$ is Stein, there exists
$a\in C^\infty(V,\mathbb C)$ with
$\bar\partial a=\eta^{0,1}$.
Since $\eta$ is real, $\partial\bar a=\eta^{1,0}$.
Thus, for $\Psi:=\operatorname{Im}a$,
\[
 \begin{aligned}
 dd^c\Psi
 &=\partial\bar\partial(a-\bar a)\\
 &=\partial\eta^{0,1}+\bar\partial\eta^{1,0}
 =d\eta=\alpha.
 \end{aligned}
\]
\end{proof}


\begin{theorem}[=Theorem~\ref{thm:intro-normal-kahler-uniruled}]\label{thm:normal-kahler-uniruled}
Let $(X,B)$ be a compact K\"ahler pair, where $X$ is normal of dimension $n$, $B\geq0$, and $K_X+B$ is
$\mathbb R$-Cartier. Let $\omega$ be a K\"ahler form on $X$.
If $[K_X+B]+t[\omega]$ is not pseudoeffective for some $t>0$,
then for every $x\in X$ there is a rational curve $R\subset X$
through $x$ such that
\begin{equation}\label{eq:normal-uniruled-degree-bound}
 \int_R\omega\leq\frac{n+1}{t}.
\end{equation}
In particular, if $K_X+B$ is not pseudoeffective, then $X$ is
uniruled.
\end{theorem}

\begin{proof}
By \cite{BM97}, there is a log resolution $\pi:M\to X$ of $(X,B)$,
which is an isomorphism over
$X_{\mathrm{reg}}\setminus\Supp B$, obtained by a sequence of blow-ups along smooth center. Let $\widetilde B$ be the strict transform of $B$. Since $K_X+B$ is $\mathbb{R}$-Cartier, we can write
\begin{equation}\label{eq:canonical-resolution}
 K_M+\widetilde B=\pi^*(K_X+B)+\sum_i a_iE_i,
\end{equation}
where the $E_i$ are $\pi$-exceptional prime divisors.

We first show that $[K_M+t\pi^*\omega]$ is not
pseudoeffective. Suppose, on the contrary, that this class
contains a closed positive current $T$.
Choose a smooth representative $\sigma$ of
$[K_X+B+t\omega]$. The current
\[
 S:=T+[\widetilde B]-\sum_i a_i[E_i]
\]
represents $\pi^*[\sigma]$, and hence
\[
 S=\pi^*\sigma+\ddc v
\]
for an integrable function $v$ on $M$, which is actually the difference of quasi-plurisubharmonic
functions.
There is an analytic subset $A_\pi\subset X$ of codimension at least $2$ such that $\pi$ is an isomorphism over
$X\setminus A_\pi$ and $E_i\subset \pi^{-1}(A_\pi)$ for every $i$.
On $X\setminus A_\pi$, the function
$v_0:=v\circ\pi^{-1}$ is $\sigma$-plurisubharmonic,
since $\pi^*(\sigma+dd^cv_0)=T+[\widetilde B]\geq0$ on $\pi^{-1}(X\setminus A_\pi)$. The extension
theorem \cite{GR56} for plurisubharmonic functions on normal spaces extends $v_0$ to a $\sigma$-psh function on $X$. Thus $[K_X+B+t\omega]$ contains a
closed positive current, a contraction.

Fix a K\"ahler form $\kappa$ on $M$. Since the pseudoeffective
cone of $M$ is closed, we can choose $\varepsilon>0$ sufficiently
small and $s>t$ sufficiently close to $t$ such that
\[
 [K_M]+s\bigl(\pi^*[\omega]+\varepsilon[\kappa]\bigr)
 \quad\text{is not pseudoeffective}.
\]
Put
\[
 \Omega_\varepsilon:=\pi^*\omega+\varepsilon\kappa,
 \qquad
 \theta_\varepsilon:=\operatorname{Ric}(\Omega_\varepsilon).
\]
Then $\Omega_\varepsilon$ is a K\"ahler form on $M$. Set
\[
 \rho_s:=-\theta_\varepsilon+2\pi s\Omega_\varepsilon,
 \qquad
 E_s(p):=\inf_{f\in\mathcal A(M,p)}Q_{\rho_s}(f),
 \quad p\in M.
\]
Observe that $\rho_s$ represents the class $
 2\pi[K_M+s\Omega_\varepsilon]$, which cannot be pseudoeffective by the discussions above. By Theorem~\ref{thm:disc functional}, we obtain that for every $p\in M$,
\[
 \sup_{f\in\mathcal A(M,p)}
 Q_{\theta_\varepsilon-2\pi s\Omega_\varepsilon}(f)
 =-E_s(p)=+\infty.
\]

Put $\delta=2\pi t$ and $c=(s-t)/2>0$.
The inequality $2\log(1/|z|)\geq1-|z|^2$ on $\mathbb D$
gives $4\pi Q_{\Omega_\varepsilon}(f)
\geq W_{\Omega_\varepsilon}(f)$, and therefore
\begin{align*}
 Q_{\theta_\varepsilon-\delta\Omega_\varepsilon}(f)
       -cW_{\Omega_\varepsilon}(f)
 &=Q_{\theta_\varepsilon-2\pi s\Omega_\varepsilon}(f)
   +c\bigl(4\pi Q_{\Omega_\varepsilon}(f)
                  -W_{\Omega_\varepsilon}(f)\bigr)\\
 &\geq Q_{\theta_\varepsilon-2\pi s\Omega_\varepsilon}(f).
\end{align*}
Theorem~\ref{thm:pointwise-criterion} now gives, for every
$p\in M$, a rational curve $C\subset M$ through $p$ such that
\[
 \int_C\Omega_\varepsilon
 \leq\frac{2\pi(n+1)}{\delta}
 =\frac{n+1}{t}.
\]

Fix $x\in X$. choose $y_j\in X_{\rm reg}\setminus \Supp B$
with $y_j\to x$. For each $j$, let $C_j\subset M$ be a
rational curve through $\pi^{-1}(y_j)$ with the degree bound $\int_{C_j}\Omega_\varepsilon
 \leq\frac{n+1}{t}$. Since $\pi$ is an isomorphism near
$\pi^{-1}(y_j)$, each $\pi(C_j)$ is a rational curve containing $y_j$.
Since $C_j$ has uniformly bounded $\Omega_\varepsilon$-area, Lemma~\ref{lem:rational-specialization} yields, after passing to
a subsequence, the currents $[C_j]$ converge weakly to an
effective rational $1$-cycle
\[
 Z=\sum_{\nu=1}^N m_\nu C_\nu,
 \qquad m_\nu\in\mathbb Z_{>0},
\]
and there is a component $C_\nu$ such that $R:=\pi(C_\nu)$
is a rational curve through $x$.
It follows that
\begin{align*}
 \int_R\omega
 &\leq\int_{C_\nu}\pi^*\omega
 \leq\int_{C_\nu}\Omega_\varepsilon
 \leq\int_Z\Omega_\varepsilon\\
 &=\lim_{j\to\infty}\int_{C_j}\Omega_\varepsilon
 \leq\frac{n+1}{t}.
\end{align*}
The first inequality uses $\deg(C_\nu\to R)\int_R\omega=\int_{C_\nu}\pi^*\omega$. 
This proves the degree bound \eqref{eq:normal-uniruled-degree-bound} and the proof is complete.
\end{proof}



\vspace{0.25cm}
\section{Rational curves on compact K\"ahler spaces with non-nef canonical class}
\label{sec:mixed-intersections}
We first establish a technical statement which allows us to take the limit of a sequence of rational curves with uniformly bounded area.
\begin{lemma}\label{lem:rational-specialization}
Let $(M,\omega)$ be a compact K\"ahler manifold, and let
$C_j\subset M$ be rational curves such that
\[
 \sup_j\int_{C_j}\omega<\infty.
\]
After passing to a subsequence, there is a holomorphic $1$-cycle
\[
 Z=\sum_{\nu=1}^N m_\nu C_\nu,
 \qquad m_\nu\in\mathbb Z_{>0},
\]
such that 
\begin{enumerate}
    \item the integration currents $[C_j]$ converge weakly to $[Z]$;
    \item $C_j$ converge to $|Z|$ in the Hausdorff topology: $x\in|Z|$ precisely when $x=\lim_{j\to\infty}x_{j}$ for some $x_{j}\in C_{j}$, and for any $\varepsilon>0$, $C_j$ belongs to the $\varepsilon$-neighborhood of $|Z|$ for $j$ large.
\end{enumerate}
Let $X$ be a complex space and $\pi:M\to X$ a holomorphic map.
Suppose that each $\pi(C_j)$ is a curve and that there are points $x_j\in\pi(C_j)$
with $x_j\to x$ in $X$. Then $\pi(|Z|)$ contains a rational curve
through $x$.
\end{lemma}

\begin{proof}
Since the volume of $C_j$ is uniformly bounded, Bishop's compactness theorem \cite{Bishop} (see also \cite[Theorems~3.9,\,3.10]{HS74}) gives an effective $1$-cycle $Z$ satisfying the desired convergence properties. We only need to prove that coefficient $m_\nu$ of $C_\nu$ is a positive integer.
We pick some point $p\in (C_\nu)_{\rm reg}\setminus\bigcup_{\mu\neq\nu}C_\mu$. In a local coordinate neighborhood $U\simeq\D^n$ we have $C_\nu\cap U=\D_{z_1}\times\{0\}$. Let $q:\mathbb{D}^n\to\mathbb{D}_{z_1}$ be the projection.
The Hausdorff convergence $C_j\to|Z|$ implies that $C_j\cap U\subset\D\times\frac{1}{2}\D^{n-1}$ when $j$ is large enough. Therefore $q|_{C_j\cap U}:C_j\cap U\to\D$ is finite, surjective. We denote its degree by $d_j$.
We fix a $(1,1)$-form $\eta(z_1)$ compactly supported in $\D$ and $\chi(z_2,\cdots,z_n)$ be a cut-off function on $\D^{n-1}$ which is identically one on $\frac{1}{2}\D^{n-1}$.
On one hand we have $d_j\int_{\D}\eta=\int_{C_j}\chi q^*\eta$. On the other hand the weak convergence $\lim[C_j]=[\sum m_{\nu}C_\nu]$ gives us
\[
\lim_{j\to \infty}\int_{C_j}\chi q^*\eta=m_\nu\int_{C_\nu}\chi q^*\eta=m_\nu\int_{\D}\eta.
\]
This proves $m_\nu=\lim d_j$ is a positive integer.

\vspace{0.1cm}
We prove that every component $C_\nu$ of $Z$ is rational.
Suppose it is not the case, then the normalization of some $C_\nu$ has
positive genus. We claim that there are smooth closed
real one-forms $a,b$ on $C_\nu$, which are compactly supported on $T:=(C_\nu)_{\mathrm{reg}}\setminus\bigcup_{\mu\ne\nu}C_\mu$ and satisfy
\begin{equation}\label{eq:handle-intersection}
 \int_T a\wedge b=1.
\end{equation}
Indeed, identify $T$ with its inverse image in the normalization
$\widehat C_\nu$, so that $\widehat C_\nu\setminus T$ is a finite set.
Since $\widehat C_\nu$ is a Riemann surface with positive genus, it admits a nonzero
holomorphic $1$-form $\eta$. For each $p\in\widehat C_\nu\setminus T$, choose small discs $D_p^\prime\Subset D_p\Subset\widehat C_\nu$ of $p$ and cut-off functions $\chi_p$ supported in $D_p$ with $\chi|_{D_p^\prime}\equiv1$. 
There is holomorphic function $f_p$ on $D_p$ so that $df_p=\eta$, then the closed form
\[
 \eta':=\eta-\sum_{p\in \widehat C_\nu\setminus T} d(\chi_p f_p)
\]
is then compactly supported in $T$. By Stokes' theorem, we have
\[
 \int_T \operatorname{Re}\eta'\wedge\operatorname{Im}\eta'
 =\frac{i}{2}\int_{\widehat C_\nu}\eta\wedge\overline{\eta}>0,
\]
and the claim follows by taking $a=\operatorname{Re}\eta',b=\operatorname{Im}\eta'$.

Choose an open set $T_0\Subset T$ containing $\Supp a\cup\Supp b$.
Let $N$ be a tubular neighborhood of $T_0$ in $M$ with
$N\cap |Z|=T_0$, and $p:N\to T_0$ be the retraction.
Choose $\chi\in C_c^\infty(N)$ equal to $1$ in a neighborhood $N_0\Subset N$ of
$\Supp a\cup\Supp b$. Define
\[
 \widetilde a=\chi p^*a,\qquad \widetilde b=\chi p^*b
\]
on $N$ and extend both forms by zero to $M$.
They are closed in $N_0$.
They vanish on a neighborhood of $|Z|\setminus N_0$.
Therefore $\widetilde a,\widetilde b$ are
closed in an open neighborhood of $|Z|$.
This neighborhood contains $C_j$
for sufficiently large $j$ by Hausdorff convergence.

Let $\nu_j:\PP^1\to C_j$ be the normalization. Since
$H^1(\PP^1,\mathbb R)=0$, $\nu_j^*\widetilde a$ and $\nu_j^*\widetilde b$ are exact.
It follows from Stokes' theorem that
\[
 \int_{C_j}\widetilde a\wedge\widetilde b=0.
\]
Thanks to the weak convergence $[C_j]\to [Z]$ and
\eqref{eq:handle-intersection}, we get
\[
 0=\lim_{j\to\infty}\int_{C_j}\widetilde a\wedge\widetilde b=\int_Z\widetilde a\wedge\widetilde b
   =m_\nu\int_T a\wedge b=m_\nu>0,
\]
a contradiction. Thus every component of $Z$ is rational.

\vspace{0.1cm}
Finally, we descend rational curves on $M$ to $X$. Since $\pi$ is uniformly continuous on the compact manifold $M$, we see $\pi(C_j)$ converge in the Hausdorff topology to $\pi(|Z|)$. $C_j$ is irreducible and $C_j\ni x_j\to x$ implies that $\pi(|Z|)$ is connected and contains
$x$. Let $U\subset X$ be an open Stein neighborhood of $x$, then every analytic set $\pi(C_j)$ meets
$X\setminus U$. The Hausdorff convergence shows that $\pi(|Z|)$ also
meets $X\setminus U$. If $\pi$ is constant on every component of $Z$ whose image
contains $x$, then $x$ would be an isolated point
of $\pi(|Z|)$, which contradicts its connectedness.
So there is a component $C_\nu$ so that $R=\pi(C_\nu)$ is a curve. $R$ is rational since $C_\nu$ is. The proof is complete.
\end{proof}

We can now prove Theorem~\ref{thm:intro-normal-kahler-rational-curves}.

\begin{theorem}[=Theorem~\ref{thm:intro-normal-kahler-rational-curves}]\label{thm:normal-kahler-rational-curves}
Let $(X,B)$ be a compact K\"ahler pair. Let $Y\subset X$ be an
irreducible subvariety of dimension $d\geq1$, and let
$\omega_1,\ldots,\omega_{d-1}$ be nef $(1,1)$-classes on $X$.
Set
\[
 L_Y(\gamma):=
 \int_Y\gamma\wedge\omega_1\wedge\cdots\wedge\omega_{d-1},\qquad\gamma\in H^{1,1}_{\rm BC}(X).
\]
Assume that $L_Y(K_X+B)<0$ and
\begin{equation}\label{eq:good-generic-point}
 Y\cap\bigl(X_{\mathrm{reg}}\setminus\Supp B\bigr)
 \ne\varnothing.
\end{equation}
Then for every nef $(1,1)$-class $\alpha$ on $X$ and every
$x\in Y$, there is a rational curve $R\subset X$
through $x$ such that
\begin{equation}\label{eq:main-degree-bound}
 \alpha\cdot R\leq
 (n+1)\frac{L_Y(\alpha)}{-L_Y(K_X+B)}.
\end{equation}
\end{theorem}

\begin{proof}
Fix a K\"ahler form $\omega_0$ on $X$. Set
\[
 A=L_Y(\alpha),\qquad A_0=L_Y([\omega_0]),\qquad k=-L_Y(K_X+B).
\]
It is clear that $A,A_0\geq0$ and $k>0$. We can actually deduce that $A_0>0$. Since $K_X+B$ is $\mathbb R$-Cartier, the class $[K_X+B+C\omega_0]$ is K\"ahler for sufficiently
large $C>0$, and hence $0\leq L_Y(K_X+B+C\omega_0)= -k+CA_0$.

By \cite{BM97}, there is a projective log
resolution $\pi:M\to X$, which is an isomorphism over $U:=X_{\mathrm{reg}}\setminus\Supp B$.
Then $M$ is a compact K\"ahler manifold. Let $\widetilde B$ be
the strict transform of $B$. We write
\[
 K_M+\widetilde B
 =\pi^*(K_X+B)+\sum_i a_iE_i,
\]
where the $E_i$ are $\pi$-exceptional prime divisors.
It is well-known that there are an effective
$\pi$-exceptional $\mathbb R$-divisor $F$ and a K\"ahler form
$\kappa$ on $M$ such that
\begin{equation}\label{eq:resolution-kahler-class}
 [\kappa]=\pi^*[\omega_0]-[F].
\end{equation}
 For $\varepsilon>0$, choose a K\"ahler form $\Omega_\varepsilon$ in
the K\"ahler class
\begin{equation}\label{eq:upstairs-metric}
 H_\varepsilon:=\pi^*\alpha+\varepsilon[\kappa]
 =\pi^*(\alpha+\varepsilon[\omega_0])-\varepsilon[F].
\end{equation}

We first construct rational curve through $y\in Y\cap(X_{\mathrm{reg}}\setminus\Supp B)$. Let
$\widetilde y:=\pi^{-1}(y)$. Suppose that
\begin{equation}\label{eq:choice-t}
 0<t<\frac{k}{A+\varepsilon A_0}.
\end{equation}
The key observation is that: a closed positive current in
$[K_M+tH_\varepsilon]$ cannot have a local potential which is finite
at $\widetilde y$.

Suppose to the contrary, let $T$ be such a current. Let $\sigma$
be a smooth representative of
\[
 [K_X+B+t(\alpha+\varepsilon\omega_0)].
\]
By \eqref{eq:canonical-resolution} and
\eqref{eq:upstairs-metric}, the current
\[
 S:=T+[\widetilde B]-\sum_i a_i[E_i]+t\varepsilon[F]
\]
represents $\pi^*[\sigma]$, so there is a locally integrable function $v$, which is a difference of quasi-plurisubharmonic functions, such that
\begin{equation}\label{eq:signed-current-potential}
 S=\pi^*\sigma+\ddc v.
\end{equation}
Since $E_i,F,\widetilde B$ lie in $\pi^{-1}(X_{\rm sing}\cup \Supp B)$, they do not pass through $\widetilde y=\pi^{-1}(y)$. Moreover, $T$ has a local potential which is finite at $\tilde y$ by assumption, and we therefore deduce that $v$ is finite at $\widetilde y$.

As in the proof of Theorem~\ref{thm:normal-kahler-uniruled}, by the extension theorem of plurisubharmonic functions on normal spaces, we obtain a $\sigma$-psh function $v_0$ on $X$ so that $v=v_0\circ\pi$ and $v_0(y)>-\infty$.

Let $\nu:\widehat Y\to Y$ be a resolution
with $\widehat Y$ compact K\"ahler. Since $v_0$ is not identically $-\infty$ on $Y$, $[\nu^*\sigma]$ is pseudoeffective
on $\widehat Y$. This implies
$\int_{\widehat Y}[\nu^*\sigma]\wedge\nu^*\omega_1\wedge\cdots\wedge\nu^*\omega_{d-1}\geq0$.
On the other hand, using \eqref{eq:choice-t} we get
\[
 \int_{\widehat Y}[\nu^*\sigma]\wedge
       \nu^*\omega_1\wedge\cdots\wedge\nu^*\omega_{d-1}
 =\int_{Y}[\sigma]\wedge\omega_1\wedge\cdots\wedge\omega_{d-1}
 =L_Y(K_X+B+t(\alpha+\varepsilon \omega_0))
 =-k+t(A+\varepsilon A_0)<0.
\]
A contradiction. 

\vspace{0.2cm}
Put
\[
 E(p):=\inf_{f\in\mathcal A(M,p)}Q_{-\operatorname{Ric}(\Omega_\varepsilon)+2\pi t\Omega_\varepsilon}(f),\qquad p\in M.
\]
By Theorem~\ref{thm:disc functional}, $E_t$ is either identically
$-\infty$ or $\rho_t$-psh. However $[-\operatorname{Ric}(\Omega_\varepsilon)+2\pi t\Omega_\varepsilon]=2\pi[K_M+tH_\varepsilon]$, even in the latter case, we still get to know $E_t(\tilde{y})=-\infty$ by the previous analysis. Consequently, if we set $ \theta_\varepsilon:=\operatorname{Ric}(\Omega_\varepsilon)$, we can then obtain
\begin{equation}\label{eq:canonical-envelope-divergence}
 \sup_{f\in\mathcal A(M,\widetilde y)}
 Q_{\theta_\varepsilon-2\pi t\Omega_\varepsilon}(f)
 =- E(\widetilde y)=+\infty.
\end{equation}
Set $ s_\varepsilon=\frac{2\pi k}{A+\varepsilon A_0}$. For $0<\delta<s_\varepsilon$, choose $t$ satisfying
\eqref{eq:choice-t} and $\delta<2\pi t$, and put
$c=(2\pi t-\delta)/(4\pi)>0$. Then we have
\begin{align*}
 Q_{\theta_\varepsilon-\delta\Omega_\varepsilon}(f)
       -cW_{\Omega_\varepsilon}(f)
 &=Q_{\theta_\varepsilon-2\pi t\Omega_\varepsilon}(f)
       +c\bigl(4\pi Q_{\Omega_\varepsilon}(f)
                    -W_{\Omega_\varepsilon}(f)\bigr)\\
 &\geq Q_{\theta_\varepsilon-2\pi t\Omega_\varepsilon}(f).
\end{align*}
\eqref{eq:canonical-envelope-divergence} and
Theorem~\ref{thm:pointwise-criterion} give a rational curve
$C\subset M$ through $\widetilde y$ with
\begin{equation}\label{eq:upstairs-disc-bound}
 \int_C\Omega_\varepsilon
 \leq\frac{2\pi(n+1)}{\delta}.
\end{equation}
Since $\pi$ is isomorphic over a neighborhood of $y$, 
$\pi(C)$ is a rational curve through $y$.

Now we deal with general $x\in Y$.
By \eqref{eq:good-generic-point}, there are
$y_j\in Y\cap(X_{\mathrm{reg}}\setminus\Supp B)$ with
$y_j\to x$. Choose
$\delta_j\uparrow s_\varepsilon$, and take rational curves
$C_j$ with degree bound \eqref{eq:upstairs-disc-bound} through point
$\pi^{-1}(y_j)$. 
In particular $\int_{C_j}\Omega_\varepsilon$ are uniformly
bounded. Lemma~\ref{lem:rational-specialization} then gives an integral rational $1$-cycle $\sum m_{\nu,\varepsilon} C_{\nu,\varepsilon}$ so that $R_\varepsilon=\pi(C_{\nu,\varepsilon})$ is a rational curve through $x$ for some $\nu$.
Moreover, using weak convergence of currents $\lim[C_j]=[\sum m_{\nu,\varepsilon} C_{\nu,\varepsilon}]$ we obtain that
\begin{equation}\label{eq:epsilon-degree-bound}
 \int_{R_\varepsilon}\alpha\leq\int_{C_{\nu,\varepsilon}}(\pi^*\alpha+\varepsilon\kappa)
 \leq\lim_{j\to\infty}\int_{C_j}(\pi^*\alpha+\varepsilon\kappa)
 \leq \frac{2\pi(n+1)}{s_{\varepsilon}}
 = (n+1)\frac{A+\varepsilon A_0}{k}.
\end{equation}

Suppose that no rational curve $R\subset X$ through $x$
satisfies $\alpha\cdot R\leq \frac{(n+1)A}{k}$. 
We simplify our previous notations by $R_\varepsilon=\pi(C_{\varepsilon})$, then
\[
 \pi^*\alpha\cdot C_\varepsilon
 \geq\alpha\cdot R_\varepsilon>\frac{(n+1)A}{k}.
\]
Together with \eqref{eq:epsilon-degree-bound}, this implies
\begin{equation}\label{eq:uniform-area-last-limit}
 [\kappa]\cdot C_\varepsilon<\frac{(n+1)A_0}{k},
 \qquad
 \pi^*\alpha\cdot C_\varepsilon
       \leq \frac{(n+1)(A+\varepsilon A_0)}{k}.
\end{equation}
Rational curves $C_\varepsilon$ have uniformly bounded $\kappa$-area,
apply Lemma~\ref{lem:rational-specialization} again we get a sequence $C_\varepsilon$ converge to an integral $1$-cycle $\sum m_\nu C_\nu$ and $R=\pi(C_\nu)$ is a rational curve through $x$ for some $\nu$. Moreover,
\[
 \alpha\cdot R\leq
 \pi^*\alpha\cdot C_\nu
 \leq\lim_{\varepsilon\to0}\pi^*\alpha\cdot C_\varepsilon
 =\frac{(n+1)A}{k},
\]
a contradiction. This proves the degree bound
\eqref{eq:main-degree-bound} and completes the proof.
\end{proof}

As a direct consequence of Theorem~\ref{thm:normal-kahler-rational-curves}, we can prove the following result.
\begin{corollary}[= Corollary~\ref{cor:intro-nef-threshold}]\label{cor:nef-threshold}
Let $(X,B)$ be a compact K\"ahler pair, where $X$ is normal and of dimension $n$, $B\geq0$, and $K_X+B$ is
$\mathbb R$-Cartier. Suppose that $K_X+B$ is not nef. Fix a
K\"ahler form $\omega$ on $X$ and set
\[
 \tau:=\inf\{t\geq0:[K_X+B]+t[\omega]\text{ is nef}\},
 \qquad
 \alpha:=[K_X+B]+\tau[\omega].
\]
Let $Z\subset X$ be an irreducible analytic subvariety of
dimension $d\geq1$ such that
\[
 \int_Z\alpha^d=0,
 \qquad
 Z\cap\bigl(X_{\mathrm{reg}}\setminus\operatorname{Supp}B\bigr)
 \ne\varnothing.
\]
Then every point $x\in Z$ lies on a rational curve $R\subset X$
satisfying
\begin{equation}\label{eq:normal-null-locus-degree}
 \alpha\cdot R=0,
 \qquad
 0<-(K_X+B)\cdot R
 =\tau\int_R\omega\leq n+1.
\end{equation}
\end{corollary}

\begin{proof}
Since $K_X+B$ is not nef, we have $\tau>0$, and $\alpha$ is
nef by closeness of the nef cone. Fix $x\in Z$. Let
$q\in\{0,\ldots,d-1\}$ be the largest integer such that
\[
 \int_Z\alpha^q\wedge[\omega]^{d-q}>0.
\]
Set
\[
 L_Z(\gamma):=
 \int_Z\gamma\wedge\alpha^q\wedge[\omega]^{d-q-1}.
\]
Then
\[
 L_Z(\alpha)=0,\qquad L_Z([\omega])>0,\qquad
 -L_Z(K_X+B)=\tau L_Z([\omega])>0.
\]

Put $U:=X_{\mathrm{reg}}\setminus\operatorname{Supp}B$.
Choose $y_j\in Z\cap U$ with $y_j\to x$. Applying
Theorem~\ref{thm:normal-kahler-rational-curves} to the class
$[\omega]+j\alpha$ and the point $y_j$, we obtain a rational
curve $R_j\subset X$ through $y_j$ such that
\begin{equation}\label{eq:null-locus-approximating-curves}
 [\omega]\cdot R_j+j\alpha\cdot R_j
 \leq(n+1)\frac{L_Z([\omega]+j\alpha)}{-L_Z(K_X+B)}
 =\frac{n+1}{\tau}.
\end{equation}

Choose a projective log resolution $\pi:M\to X$ which is an
isomorphism over $U$. As in the proof of
Theorem~\ref{thm:normal-kahler-rational-curves}, we may assume there is an effective $\pi$-exceptional
$\mathbb R$-divisor $F$ and a K\"ahler form $\kappa$ on $M$
such that
\[
 [\kappa]=\pi^*[\omega]-[F].
\]
Let $C_j\subset M$ be the strict transform of $R_j$. Since
$y_j\in R_j\cap U$, the map $C_j\to R_j$ is bimeromorphic, it follows that
$C_j$ is rational and
$C_j\not\subset\operatorname{Exc}(\pi)$. Hence
$F\cdot C_j\geq0$. Using
\eqref{eq:null-locus-approximating-curves} and the nefness of
$\alpha$, we obtain
\begin{align*}
 \int_{C_j}\kappa
 =[\omega]\cdot R_j-F\cdot C_j
 \leq[\omega]\cdot R_j
 &\leq\frac{n+1}{\tau},\\
 \pi^*\alpha\cdot C_j
 =\alpha\cdot R_j&\leq\frac{n+1}{j\tau}.
\end{align*}
Lemma~\ref{lem:rational-specialization} gives, after passing
to a subsequence, an effective rational $1$-cycle
\[
 \Gamma=\sum_{\nu=1}^N m_\nu C_\nu,
 \qquad m_\nu\in\mathbb Z_{>0},
\]
such that $\lim_{j\to \infty}[C_j]=[\Gamma]$ weakly and some component
$C=C_\nu$ of $\Gamma$ maps onto a rational curve $R:=\pi(C)$
through $x$.
Since $\pi^*\alpha$ is nef, weak convergence gives
\[
 0\leq\pi^*\alpha\cdot C
 \leq\pi^*\alpha\cdot\Gamma
 =\lim_{j\to\infty}\pi^*\alpha\cdot C_j
 =\lim_{j\to\infty}\alpha\cdot R_j=0.
\]
Then we get $\alpha\cdot R=\deg(C\to R)^{-1}\pi^*\alpha\cdot C=0$.
Similarly, since $\pi^*[\omega]$ is nef, we have
\[
 \int_R\omega
 \leq\pi^*[\omega]\cdot C
 \leq\pi^*[\omega]\cdot\Gamma
 =\lim_{j\to\infty}\pi^*[\omega]\cdot C_j
 =\lim_{j\to\infty}[\omega]\cdot R_j
 \leq\frac{n+1}{\tau}.
\]
Finally, $\alpha=[K_X+B+\tau\omega]$ and
$\alpha\cdot R=0$ imply
\[
 0<-(K_X+B)\cdot R
 =\tau\int_R\omega\leq n+1.
\]
This completes the proof.
\end{proof}

\section{A proof of Kobayashi canonical ampleness conjecture}
\label{sec:hyperbolic-canonical-bigness}

In this section, we give an immediate proof of Kobayashi's canonical ampleness conjecture based on our understanding.
\begin{theorem}[= Theorem~\ref{thm:intro-hyperbolic-canonical-ampleness}]\label{thm:hyperbolic-canonical-bigness}
Let $X$ be a compact K\"ahler manifold.
If every holomorphic map $\mathbb C\to X$ is constant, then $K_X$ is ample.
\end{theorem}
OpenAI \cite{OAI26} recently claimed a proved of this result. Their main idea is to use the second variation of disc functional of the following form
\[
|f'(0)|^2_{\omega}+\frac{1}{2\pi}\int_\D\log\frac{1}{|z|}f^*\rho-\delta\int_\D f^*\omega,\qquad f\in \mathcal{A}_x.
\]
In this article, the logarithmic disc functional
$
Q_\rho(f)=\frac{1}{2\pi}\int_\D \log\frac{1}{|z|}\,f^*\rho
$
in Theorem~\ref{thm:pointwise-criterion} is used to construct rational curves when $K_X$ is not nef. 
If $X$ contains no entire curves, it contains no rational
curves, therefore we get to know $K_X$ is nef. 
It is well-known that if $K_X$ is big and $X$ contains no rational curve, then $K_X$ is ample.
Consequently, Theorem~\ref{thm:hyperbolic-canonical-bigness} is reduced to showing the bigness of $K_X$. 
It is then natural to dispense with the logarithmic term and use a simpler disc functional, which is defined in \eqref{eq:hyp-disc-functional}.

 
Since $K_X$ is nef, to prove the bigness of $K_X$, it suffices to show that $\int_Xc_1(K_X)^n>0$ by \cite{DP04}.
Fix a K\"ahler form $\omega_0$ on $X$.
Since $K_X$ is nef, the class $2\pi c_1(K_X)+t[\omega_0]$ is K\"ahler for $0<t\leq1$.
By Yau's theorem, there is a K\"ahler form
$\omega_t$ in this class satisfying
\begin{equation}\label{eq:hyp-twisted-ricci}
 \Ric(\omega_t)=-\omega_t+t\omega_0.
\end{equation}
We only need to prove that $\int_X\omega_t^n$ has a positive lower bound independent of $t$. If so, we can conclude that
\[
 (2\pi)^n\int_Xc_1(K_X)^n
 =\lim_{t\downarrow0}\int_X\omega_t^n\geq c.
\]

\begin{proof}
Fix $x\in X$. Choose holomorphic coordinates $(x^1,\ldots,x^n)$ on a
coordinate ball $U$ centered at $x$, and put
$u:=\sum_{k=1}^n|x^k|^2$ on $U$. Since there is no entire curve on $X$, Brody's lemma \cite{Bro78} gives a constant $C_0$ such that
$|df(0)|_{\omega_0}\leq C_0$ for every holomorphic disc $f:\D\to X$.
For $z\in\D$, the automorphism $\phi_z(\zeta):=\frac{z+\zeta}{1+\bar z\zeta}$
satisfies $\phi_z(0)=z$ and $\phi_z'(0)=1-|z|^2$. The above derivative bound gives us
\begin{equation}\label{eq:hyp-brody}
 (1-|z|^2)|df(z)|_{\omega_0}
 =|d(f\circ\phi_z)(0)|_{\omega_0}\leq C_0,
 \qquad z\in\D.
\end{equation}
It follows that there is some $0<r_0<1/2$ and $V\Subset U$ so that every holomorphic disc with $f(0)=x$ satisfies $f(\overline D_{2r_0})\subset V$. Write $M:=\sup_V u$.

\vspace{0.2cm}
For a holomorphic disc $f$ with $f(0)=x$, we introduce the disc functional
\begin{equation}\label{eq:hyp-disc-functional}
 \mathcal J_t(f):=A_t(f)-P(f),\qquad A_t(f):=\int_\D f^*\omega_t,\qquad
 P(f):=\frac{1}{2\pi}\int_0^{2\pi}
   u(f(r_0e^{i\theta}))\,d\theta .
\end{equation}
It is clear that $0\leq P(f)\leq M$, and $\mathcal J_t(\text{constant disc})=0$. 
Hence a minimizing sequence of $\mathcal J_t$ has uniformly bounded $\omega_t$-area.
By \eqref{eq:hyp-brody} and compactness of $X$, a subsequence denoted by $f_j$ converges
locally uniformly on $\D$ to a holomorphic disc $f$ with $f(0)=x$.
Fatou's lemma implies $\int_\D f^*\omega_t\leq\liminf_{j\to \infty}\int_\D f_j^*\omega_t$,
the uniform convergence on $\{|z|=r_0\}$ implies $\lim_{j\to \infty}P(f_j)=P(f)$. Thus $f$
is a minimizer of $\mathcal J_t$ and
\begin{equation}\label{eq:hyp-area-bound}
 A_t(f)\leq \frac{1}{2\pi}\int_0^{2\pi}u(f(r_0e^{i\theta}))\,d\theta\leq M.
\end{equation}

To apply the variational argument, we check that $f\in W^{1,p}(\D,X)$ for $2<p<3$. For $0<r<1$ we set
$f_r(z):=f(rz)$. By minimality of $f$, we get that there is a constant $C$ depending on $C_B,r_0,u$ such that
\[
 \int_{\D\setminus D_r}f^*\omega_t
 =A_t(f)-A_t(f_r)
 \leq \frac{1}{2\pi}\int_0^{2\pi}
 \bigl(u(f(r_0e^{i\theta}))-u(f(rr_0e^{i\theta}))\bigr)\,d\theta
 \leq C(1-r),
\]
where the last inequality follows from \eqref{eq:hyp-brody} on $D_{2r_0}$.
For fixed $t$, we have $\omega_0\leq C_t\omega_t$. On the annulus
$E_{\varepsilon}:=\{\varepsilon<1-|z|<2\varepsilon\}$, the last estimate and \eqref{eq:hyp-brody} give
\[
 \int_{E_\varepsilon}|df|_{\omega_0}^p\,\dA
 \leq \bigl(\sup_{E_\varepsilon}|df|_{\omega_0}\bigr)^{p-2}
       \int_{E_\varepsilon}|df|_{\omega_0}^2\,\dA
 \leq (C/\varepsilon)^{p-2}C_t(2\varepsilon)
 = C_t' \varepsilon^{3-p}.
\]
Summing over $\varepsilon=2^{-j}$, we get $f\in W^{1,p}(\D,X)$ for $p<3$.

\vspace{0.2cm}
By Theorem~\ref{lem:baseline-frame}, there is a $W^{1,p}$ holomorphic frame $e=(e_1,\ldots,e_n)$ of $f^*T^{1,0}X$, which is unitary on $\partial\mathbb D$
with respect to $\omega_t$. Set $\xi_j(z):=ze_j(z)$,
Lemma~\ref{lem:actual} gives families of holomorphic discs $f_\tau^j$ with
$f_0^j=f$, $f_\tau^j(0)=x$ and $\partial_\tau f_\tau^j|_{\tau=0}=\xi_j$.
The variational formula of $A(f)$ in Proposition~\ref{prop:variation-sums} yields that
\[
 \sum_{j=1}^n\Delta_\tau A_t(f_\tau^j)\Big|_{\tau=0}
 =4\pi n-2\int_\D f^*\Ric(\omega_t)\leq4\pi n+2M
\]
by \eqref{eq:hyp-twisted-ricci} and \eqref{eq:hyp-area-bound}.

We compute the second variation of $P(f^j_\tau)$.
Let $F:D_\varepsilon\times\D\to X$ be a holomorphic map so that $F(0,\cdot)=f$, $\xi(z):=\partial_\tau F(0,z)\in T_{f(z)}X$. Then
$(\Delta_\tau|x^k\circ F|^2)(0,z)=4|dx^k(\xi)|^2(f(z))$ for $z\in\D$.
The family $f^j_\tau$ satisfies $\partial_{\tau}f^j_\tau|_{\tau=0}=ze_j$, we therefore obtain that
\[
\sum_{j=1}^n\Delta_\tau P(f^j_\tau)\Big|_{\tau=0}
=\Delta_\tau\Big|_{\tau=0} \int_0^{2\pi}\sum_{j,k=1}^n|x^k\circ f^j_\tau|^2(r_0e^{i\theta})\,d\theta
=4r_0^2\int_0^{2\pi}\sum_{j,k=1}^n|dx^k(e_j)|^2(f(r_0e^{i\theta}))\,d\theta
\]
If we set $\omega_{\rm std}:=\frac{i}{2}\sum_{k=1}^n dx^k\wedge d\bar x^k$ on $U$,
and denote by $G$ the Gram matrix of $e$ with respect to $\omega_{\rm std}$, then $\sum_{j,k=1}^n|dx^k(e_j)|^2(f(r_0e^{i\theta}))=\tr G(r_0e^{i\theta})$.

Since $f$ minimizes $\mathcal J_t$, we get $\Delta_\tau\mathcal J_t(f_\tau^j)|_{\tau=0}\geq0$.
The inequality $\Delta_\tau\sum_{j=1}^n\mathcal{J}_t(f_\tau^j)|_{\tau=0}\geq0$ then gives us
\[
 \frac{1}{2\pi}\int_0^{2\pi}\tr G(r_0e^{i\theta})\,d\theta\leq C_1,\qquad C_1=(4\pi n+2M)/(4r_0^2).
\]
By the arithmetic-geometric mean value inequality and Jensen's inequality, we can deduce that
\begin{equation}\label{eq:hyp-euclidean-determinant}
 \frac{1}{2\pi}\int_0^{2\pi}\log\det G(r_0e^{i\theta})\,d\theta
 \leq n\log(C_1/n).
\end{equation}

Finally, we prove $\int_X\omega_t^n$ is uniformly bounded from below.
Let $G_t$ be the Gram matrix of $e$ with respect to $\omega_t$ and put
$g_t:=\log\det G_t$. As in the proof of Proposition~\ref{prop:variation-sums}, using $e=(e_1,\cdots,e_n)$ is boundary-unitary, we obtain that
\[
 \Delta_z g_t\,\dA=-2 f^*\Ric(\omega_t),\qquad g_t|_{S^1}=0.
\]
For $0<r<1$, Stokes' formula gives
\[
 \frac{d}{dr}\frac{1}{2\pi}\int_0^{2\pi}g_t(re^{i\theta})\,d\theta
 =-\frac{1}{\pi r}\int_{D_r}f^*\Ric(\omega_t).
\]
Using $\Ric(\omega_t)\geq-\omega_t$
and \eqref{eq:hyp-area-bound}, we obtain
\begin{align*}
 \frac{1}{2\pi}\int_0^{2\pi}\log\det G_t(r_0e^{i\theta})\,d\theta
 =-\int_{r_0}^1\frac{d}{dr}\Big(\frac{1}{2\pi}\int_0^{2\pi}g_t(re^{i\theta})\,d\theta\Big)dr
 =&\frac1\pi\int_{r_0}^1\frac{1}{r}
    \left(\int_{D_r}f^*\Ric(\omega_t)\right)dr\\
 \geq&-\frac M\pi\log\frac1{r_0}.
\end{align*}
Since $G_t$ and $G$ are computed in the same frame $e$,
$\det G_t/\det G=(\omega_t^n/\omega_{\rm std}^n)\circ f$ on $\D$.
This lower bound and the upper bound in \eqref{eq:hyp-euclidean-determinant} are independent of $t$, hence at some point
$x_t\in f(r_0S^1)\subset V$ it holds that
\begin{equation}\label{eq:hyp-density}
 \frac{\omega_t^n}{\omega_{\rm std}^n}(x_t)\geq c_0>0
\end{equation}
for some $c_0$ independent of $t$.

Choose a smooth local potential $\varphi$ of $\omega_0$ on $U$, so that
$i\partial\bar\partial\varphi=\omega_0$. By \eqref{eq:hyp-twisted-ricci},
$\Ric(\omega_t)\leq\omega_0$. Hence
\[
 v_t:=\log\frac{\omega_t^n}{\omega_{\rm std}^n}+\varphi
\]
is plurisubharmonic on $U$, and so is $e^{v_t}$. Since $V\Subset U$,
there is $a>0$ independent of $t$ so that $B(x_t,a)\Subset U$. The submean
inequality of $e^{v_t}=e^\varphi\frac{\omega_t^n}{\omega_{\rm std}^n}$ at $x_t$ together with \eqref{eq:hyp-density} give
\[
\int_X\omega_t^n\geq\int_{B(y_t,a)}\omega_t^n=\int_{B(y_t,a)}e^{v_t-\varphi}\omega_{\rm std}^n\geq  c
\]
for some $c=c(\varphi,a)>0$ independent of $t$. The proof is complete.
\end{proof}

\appendix
\section{Proof of the compactness theorem}
\label{app:marked-compactness}

In this appendix, we prove Proposition~\ref{prop:interior-marked-compactness}.
We introduce several necessary notions in \cite{ISCompactness}. 
\begin{definition}
A nodal curve with boundary is a complex curve $C$ with only nodal singularities and $\overline{C}:=C\cup\partial C$ is compact and near $\partial C$, $\overline C$ is a Riemann surface with boundary.
\end{definition}

\begin{definition}\label{def:para curve}
    A nodal curve with boundary $C$ is parameterized by $\D$ if there is a continuous map $\sigma:\overline\D\to \overline C$ such that
    \begin{itemize}
        \item if $a\in C$ is a node, then $\gamma_a=\sigma^{-1}(a)$ is a smooth embedded circle in $\D$,
        \item $\sigma:\overline\D\backslash \cup_{i=1}^N\gamma_{a_i}\to \overline C\backslash\{a_1,\cdots,a_N\}$ is a diffeomorphism, where $a_1,\cdots,a_N$ are nodes of $C$.
    \end{itemize}
\end{definition}
We can easily see that such nodal curve has one principal component $C_0\simeq\D$ and some $\mathbb P^1$ trees attached to $C_0$.
In fact, 
let $\gamma_1,\cdots,\gamma_{N_1}$ be the "biggest" circles, namely not contained in the interior of any other circles.
Remaining circles are contained in $\cup_{i_1=1}^{N_1}\operatorname{Int}(\gamma_{i_i})$.
Pick out the "biggest" circles $\gamma_{i_11},\cdots,\gamma_{i_1N_{i_1}}$ in $\operatorname{Int}(\gamma_{i_1})$.
Then pick out the "biggest" circles $\gamma_{i_1i_21},\cdots,\gamma_{i_1i_2N_{i_1i_2}}$ in $\operatorname{Int}(\gamma_{i_1i_2}), \ldots$. See the figure below.

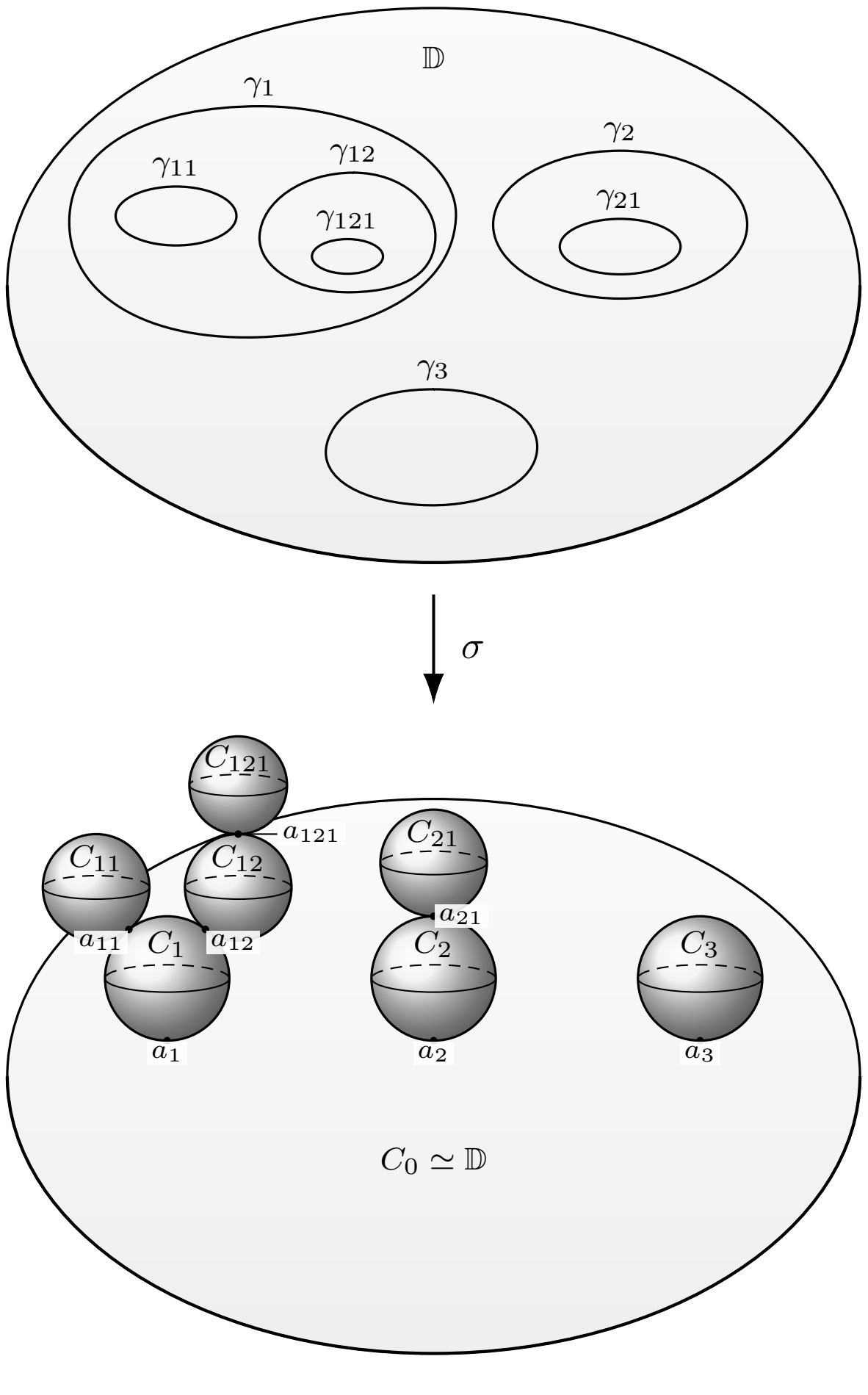
\begin{figure}[htbp]
\centering
\begin{tikzpicture}[
    scale=0.80,
    line width=0.6pt,
    every node/.style={font=\scriptsize},
    curve label/.style={
        above=1pt,
        inner sep=1pt
    },
    node label/.style={
        font=\tiny,
        inner sep=1pt,
        fill=white,
        fill opacity=0.9,
        text opacity=1
    },
    disc/.style={
        top color=gray!3,
        bottom color=gray!13,
        draw=black
    },
    sphere/.style={
        ball color=gray!12,
        draw=black
    }
]

\begin{scope}[shift={(0,8.5)}]

    \shade[disc] (0,0) ellipse (4.8 and 3.12);

    \draw[line width=0.8pt]
        (-4.8,0)
        arc[
            start angle=180,
            end angle=360,
            x radius=4.8cm,
            y radius=3.12cm
        ];

    \node at (0,2.56) {$\mathbb{D}$};

    \begin{scope}[yscale=0.65]

        \draw (-1.95,3.10)
            node[curve label] {$\gamma_1$}
            .. controls (-0.75,3.10) and (0.30,2.20)
                .. (0.25,1.15)
            .. controls (0.20,-0.15) and (-0.80,-0.85)
                .. (-2.00,-0.90)
            .. controls (-3.35,-0.95) and (-4.15,-0.05)
                .. (-4.10,1.20)
            .. controls (-4.05,2.45) and (-3.15,3.10)
                .. cycle;

        \draw (2.10,1.05) ellipse (1.43 and 1.28);
        \node[curve label] at (2.10,2.33)
            {$\gamma_2$};

        \draw (0.00,-1.80)
            node[curve label] {$\gamma_3$}
            .. controls (0.72,-1.80) and (1.23,-2.30)
                .. (1.16,-2.91)
            .. controls (1.09,-3.53) and (0.46,-3.87)
                .. (-0.18,-3.80)
            .. controls (-0.92,-3.72) and (-1.34,-3.25)
                .. (-1.18,-2.64)
            .. controls (-1.03,-2.06) and (-0.57,-1.80)
                .. cycle;

        \draw (-2.90,1.20) ellipse (0.68 and 0.51);
        \node[curve label] at (-2.90,1.71)
            {$\gamma_{11}$};

        \draw (-0.90,1.95)
            node[curve label] {$\gamma_{12}$}
            .. controls (-0.35,1.95) and (0.03,1.40)
                .. (0.02,0.83)
            .. controls (0.01,0.25) and (-0.28,-0.10)
                .. (-0.92,-0.12)
            .. controls (-1.50,-0.14) and (-1.98,0.30)
                .. (-1.96,0.85)
            .. controls (-1.94,1.45) and (-1.48,1.95)
                .. cycle;

        \draw (-0.97,0.50) ellipse (0.40 and 0.30);
        \node[curve label] at (-0.97,0.80)
            {$\gamma_{121}$};

        \draw (2.10,0.67) ellipse (0.68 and 0.48);
        \node[curve label] at (2.10,1.15)
            {$\gamma_{21}$};

    \end{scope}
\end{scope}

\draw[
    -{Latex[length=2.8mm,width=1.9mm]},
    line width=0.75pt
]
    (0,5.02)
    -- node[right=3pt,font=\small] {$\sigma$}
    (0,3.78);

\begin{scope}

    \shade[disc] (0,-0.4) ellipse (4.8 and 3.12);

    \draw[line width=0.8pt]
        (-4.8,-0.4)
        arc[
            start angle=180,
            end angle=360,
            x radius=4.8cm,
            y radius=3.12cm
        ];

    \node at (0,-1.35) {$C_0\simeq\mathbb{D}$};

    \foreach \x/\y/\r/\idx in {
        -3/0.7/0.7/1,
         0/0.7/0.7/2,
         3/0.7/0.7/3,
        -3.8/1.724695/0.6/11,
        -2.2/1.724695/0.6/12,
        -2.2/2.874695/0.55/121,
         0/2.0/0.6/21
    }{
        \begin{scope}[shift={(\x,\y)}]

            \pgfmathsetmacro{\eqradius}{0.22*\r}

            \shade[sphere] (0,0) circle (\r);

            \draw[thin,densely dashed]
                (\r,0)
                arc[
                    start angle=0,
                    end angle=180,
                    x radius=\r cm,
                    y radius=\eqradius cm
                ];

            \draw[thin]
                (-\r,0)
                arc[
                    start angle=180,
                    end angle=360,
                    x radius=\r cm,
                    y radius=\eqradius cm
                ];

            \node at (0,{0.55*\r}) {$C_{\idx}$};

        \end{scope}
    }

    \foreach \x/\idx in {-3/1,0/2,3/3}{
        \fill (\x,0) circle (1.2pt);
        \node[node label,below] at (\x,0)
            {$a_{\idx}$};
    }

    \fill (-3.430769,1.251759) circle (1.2pt);
    \node[node label,below left]
        at (-3.430769,1.251759)
        {$a_{11}$};

    \fill (-2.569231,1.251759) circle (1.2pt);
    \node[node label,below right]
        at (-2.569231,1.251759)
        {$a_{12}$};

    \fill (-2.2,2.324695) circle (1.2pt);
    \draw[thin]
        (-2.14,2.324695) -- (-1.76,2.324695);
    \node[node label,right]
        at (-1.76,2.324695)
        {$a_{121}$};

    \fill (0,1.4) circle (1.2pt);
    \node[node label,right] at (0,1.4)
        {$a_{21}$};

\end{scope}

\end{tikzpicture}
\end{figure}

The irreducible component $C_0$ of $C$ containing $\sigma(\D\backslash \cup_{i_1=1}^{N_1}\text{closed domain enclosed by }\gamma_{i_1})$ is homeomorphic to $\D$, this is because, topologically, $C_0$ is obtained by pinching $N_1$ disks in $\D$ to $N_1$ points, hence is $\D$. 
Furthermore, since $(C_0,\partial C_0)$ is a Riemann surface with boundary, one can easily prove $(C_0,\partial C_0)\simeq(\D,\partial\D)$ as Riemann surface with boundary. $C_0$ is called the principal component and $a_{1},\cdots,a_{N_1}$ are called roots.

Since $\sigma(\gamma_{i_1})=a_{i_1}$, we see that
$\sigma(\text{closed domain enclosed by }\gamma_{i_1})$ is a $\mathbb{P}^1$ tree attached to $a_{i_1}$.
In fact, the component $C_{i_1}$ of $C$ containing $\sigma(\text{domain enclosed by }\gamma_{i_1}\backslash \cup_{i_2=1}^{N_{i_1}}\text{closed domains enclosed by }\gamma_{i_1i_2})$ is homeomorphic to $S^2$, this is because topologically $C_{i_1}$ is obtained by pinching the boundary of $\overline{\D}$ to a point and pinching $N_{i_1}$ balls in $\D$ to points, hence is $S^2$.
Compact Riemann surface homeomorphic to $S^2$ can only be $\mathbb P^1$. The nodes of $C$ on $C_{i_1}$ is $a_{i_1}=\sigma(\gamma_{i_1}),a_{i_11}=\sigma(\gamma_{i_11}),\cdots,a_{i_1N_{i_1}}=\sigma(\gamma_{i_1N_{i_1}})$. The same analysis shows that the components $C_{i_11},\cdots,C_{i_1N{i_1}}$ are $\mathbb{P}^1$ attached to these nodes.

Following \cite[Definition 2.4]{ISCompactness} we define the stability of holomorphic map.
\begin{definition}
    Let $C$ be a nodal curve with boundary parametrized by $\D$ and $Y$ be a complex manifold. A holomorphic map $F:C\to Y$ is called stable if $F$ is non-constant on every sphere component $C_i$ so that at most two nodes are attached to $C_i$.
\end{definition}

A special case of \cite[Theorem~1]{ISCompactness} is the following.
\begin{theorem}\label{thm:disc compact}
    Let $X$ be a compact complex manifold,
    $M>0$ and $f_j:\D\to X$ are holomorphic maps so that $A(f_j):=\int_\D f_j^*\omega\leq M$,. 
    Then up to a subsequence, $(\iota,f_j):\D\to \mathbb{C}\times X$, where $\iota:\mathbb\D\to \mathbb{C}$ is the canonical inclusion, converge to a stable curve $(g,F):C\to \mathbb{C}\times X$ in the following sense: there are parameterizations $\sigma_j:\overline\D\to \overline\D$ and $\sigma_\infty:\overline\D\to \overline C$ in the sense of Definition~\ref{def:para curve} so that
    \begin{itemize}
        \item[(1)] $(\iota,f_j)\circ\sigma_j$ converges to $(g,F)\circ \sigma_\infty$ in $C^0_\mathrm{loc}(\D,\mathbb{C}\times X)$;
        \item[(2)] let $a_i$ be the nodes of $C$ and $\gamma_i=\sigma_\infty^{-1}(a_i)$ be the circles in $\D$, then for every compact subset $K\subset\D\backslash\cup_{i}\gamma_i$, $(\iota,f_j)\circ\sigma_j\to (g,F)\circ\sigma_\infty$ in $W^{1,p}(K,\mathbb{C}\times X)$ for any $p<\infty$;
        \item[(3)] $\sigma_j$ is identity on $\{1-\varepsilon<|z|\leq1\}$ for some $\varepsilon>0$ and for every compact subset $K\subset\overline\D\backslash\cup_j\gamma_j$, the complex structures $\sigma_j^*J_{\D}$ converge smoothly to $\sigma_\infty^*J_C$ on $K$.
    \end{itemize}
\end{theorem}
\begin{remark}\label{rem:appendix sm convergence}
\cite{ISCompactness} allows the almost complex structure of target manifold $Y$ to deform. In our situation $Y=\mathbb{C}\times X$ is equipped with the fixed complex structure, therefore $(2),(3)$ implies the $C^{\infty}$ convergence of $(\iota,f_j)\circ\sigma_j$ on compact subsets of $\D\backslash\cup_ia_i$. 
For the sake of the completeness, we give a proof of this statement.

Fix a Riemannian metric $g$ on $\D$ and a small ball $\overline{B}\subset \D\backslash\cup\gamma_i$.
Let $J_j:=\sigma_j^*J_{\D}$ be the pull-back complex structures on $\D$.
\[g_j(\cdot,\cdot):=\frac{1}{2}\Big(g(\cdot,\cdot)+g(J_j(\cdot),J_{j}(\cdot))\Big)\]
is a Riemannian metric compatible with $J_j$, converging smoothly on $\overline{B}$ by $(3)$.
Since $(\iota,f_j)\circ\sigma_j$ converges locally uniformly on $\D$ by $(1)$, after shrinking $B$ we get $(\iota,f_j)\circ\sigma_j$ maps $B$ to a common coordinate chart of $\mathbb{C}\times X$ for $j$ large.
Since $J_j=\sigma_j^*J_{\D}$ and
$\sigma_j:(\D,J_j)\to(\D,J_{\D})$ is biholomorphic, hence $u_j:=(\iota,f_j)\circ\sigma_j$ is holomorphic on $(\D,J_j)$.
In the chosen coordinates, its components are therefore
$g_j$-harmonic:
\[\Delta_{g_j}u_{j0}=\cdots=\Delta_{g_j}u_{jn}=0.\]
Then for every integer $k$, $\|u_j\|_{C^k(\frac{1}{2}B)}\leq C(k,\|g_j\|_{C^k(B)})\|u_{j}\|_{C^0(B)}$ by interior estimates of the Laplacian operator.
This proves that the maps $(\iota,f_j)\circ\sigma_j$ are locally uniformly bounded in $C^k$ on $\D\backslash \cup_i\gamma_i$ since $B$ is arbitrary. 
Finally, since $k$ is arbitrary, Arzela--Ascoli theorem implies $(\iota,f_j)\circ\sigma_j\to (g,F)\circ\sigma_\infty$ locally smoothly on $\D\backslash \cup_i\gamma_i$. 
\end{remark}
\begin{remark}\label{rem:change base D}
We prove $g|_{C_0}:C_0\to \D$ is biholomorphic and on compact subset $K\subset \D\backslash\, g(\text{roots of }C)$, $f_j$ converge smoothly to $f:=F\circ (g|_{C_0})^{-1}$.

Let $r>1-\varepsilon$, then $\sigma_j$ is identity on $\{r\leq|z|\leq1\}$. We get $g\circ\sigma_\infty=\operatorname{id}$ on this collar by $(1)$.
Let $D_{r}'\subset C_0\simeq\D$ be the domain enclosed by $\sigma_\infty(\{|z|=r\})$, maximum principle implies $g(D_r')\subset D_r$.
We prove that $D_{r}'\xrightarrow{g}D_{r}$ is biholomorphic. For any point $a\in D_{r}$, the number of zeros of $g-a$ in $D_{r}'$ is
\[
\frac{1}{2\pi i}\int_{\sigma_\infty(\{|z|=r\})}\frac{dg}{g-a}
=\frac{1}{2\pi i}\int_{|z|=r}\frac{d(g\circ \sigma_\infty)}{g\circ\sigma_\infty-a}
=\frac{1}{2\pi i}\int_{|z|=r}\frac{dz}{z-a}=1.
\]
This proves that $D_r'\to D_r$ is bijective hence biholomorphic.
Since $\sigma_\infty$ is a diffeomorphism in a neighborhood of $\{r\leq|z|\leq1\}$ and $g\circ\sigma_\infty=\operatorname{id}$ in a neighborhood of $\{r\leq|z|\leq1\}$, we get $g$ is biholomorphic on $C_0\backslash D_r'$. This proves the first statement.

On a compact subset $K\subset \D\backslash\, g(\text{roots of }C)$, we have the smooth convergence
\[
F\circ g^{-1}
=F\circ\sigma_\infty\circ(g\circ \sigma_\infty)^{-1}
=\lim_{j\to\infty} (f_j\circ \sigma_j)\circ(\sigma_j)^{-1}
=\lim_{j\to\infty}f_j.
\]
Here we use Remark~\ref{rem:appendix sm convergence} to get $\sigma_j=\iota\circ\sigma_j$ converges locally smoothly on $\D\backslash \cup_i\overline{\text{Int}\gamma_i}$ to 
\[
g\circ\sigma_\infty\big|_{\D\backslash\,\cup_i\overline{\text{Int}\,\gamma_i}}:\D\,\backslash\,\cup_i\overline{\text{Int}\,\gamma_i}\longrightarrow \D\backslash\,g(\{\text{roots of }C\}),
\]
hence its inverse converges locally smoothly to $(g\circ\sigma_\infty)^{-1}$ on $\D\backslash\,g(\{\text{roots of }C\})$. This proves the second statement.
\end{remark}

\begin{proposition}[=Proposition~\ref{prop:interior-marked-compactness}]
Let $(X,\omega)$ be a compact K\"ahler manifold.
Let $M>0$ and $f_j:\D\to X$ be holomorphic maps so that $A(f_j)=\int_\D f_j^*\omega\leq M$ and $f_j(0)=x$.
Then up to a subsequence, there is a stable curve $(g,F):C\to \mathbb{C}\times X$ with the following properties:
\begin{enumerate}
\item $g|_{C_0}:C_0\to \D$ is biholomorphic, $f_j$ converge smoothly to $f:=F\circ (g|_{C_0})^{-1}$ on every compact subset of $\D\backslash\,g(\{\text{roots of $C$}\})$;
\item there is $p\in C$ so that $g(p)=0$, $F(p)=x$;
\item For every smooth closed real $(1,1)$-form $\rho$ on $X$,
\begin{equation}\label{eq:appendix-stable-local-measure-convergence}
 \int_\D\chi\,f^*_j\rho\longrightarrow
 \int_C(g^*\chi)\, F^*\rho,
                  \qquad \forall\, \chi\in C_c(\D).
\end{equation}
\item $A_*(C)\le\liminf_j A(f_j)$.
\end{enumerate}
\end{proposition}
\begin{proof}
The statement is an almost direct corollary of Theorem~\ref{thm:disc compact} except a neck analysis for $(3)$, which is done by an integration by parts.

Apply Gromov compactness Theorem~\ref{thm:disc compact} and keep the notations in this theorem. $(1)$ is proved in Remark~\ref{rem:change base D}.
We prove $(2)$. Set $q_j:=\sigma^{-1}_j(0)$. This sequence lies in $\{|z|<1-\varepsilon\}$ because $\sigma_j=\operatorname{id}$ on $\{1-\varepsilon\leq|z|\leq1\}$. After a subsequence, $q_j\to q$. Set $p=\sigma_\infty(q)$. The local uniform convergence of $(\iota,f_j)\circ\sigma_j\to (g,F)\circ\sigma_\infty$ gives $(g,F)(p)=(g,F)\circ\sigma_\infty(q)=\lim(\iota,f_j)\circ\sigma_j(q_j)=(0,x)$.

We prove $(3)$. 
The fact that $\rho$ is $d$-closed makes the neck analysis quite easy. The convergence \eqref{eq:appendix-stable-local-measure-convergence} is equivalent to
\[
\lim_{j\to \infty}\int_\D \sigma_j^*\chi\, (f_j\circ\sigma_j)^*\rho=\int_\D (g\circ\sigma_\infty)^*\chi\, (F\circ\sigma_\infty)^*\rho.
\]
Remark~\ref{rem:appendix sm convergence} already provides the smooth convergence $\sigma_j\to g\circ\sigma_\infty$, $f_j\circ\sigma_j\to F\circ\sigma_\infty$ away from $\cup_i\gamma_i$, thus we only need to do the neck-analysis near each circle $\gamma_i$. More concretely, it suffices to prove that
\begin{equation}\label{eq:neck-analysis}
\lim_{r\downarrow0} \Big(\limsup_{j\to\infty}\int_{A_{i,r}}\big|\chi\circ\sigma_j\,(f_j\circ \sigma_j)^*\rho\big|\Big)=0,
\end{equation}
where the necks $A_{i,r}\simeq\gamma_i\times[-r,r]$ are annulus neighborhoods of the circle $\gamma_i$, contracting to $\gamma_i$ as $r\downarrow0$.
By Theorem~\ref{thm:disc compact}$(1)$, $\sigma_j\to g\circ\sigma_{\infty}$, $f_j\circ\sigma_j\to F\circ\sigma_\infty$ locally uniformly on $\D$.
So write $a_i:=\sigma_\infty(\gamma_i)\in C$, 
$x_i:=F(a_i)=F\circ\sigma_\infty(\gamma_i)\in X$,
we can choose a ball $B_i\subset X$ centered at $x_i$ so that for some $r_0>0$, $f_j\circ\sigma_j(A_{i,r})\subset\subset B_i$, when $r<r_0$ and $j> j(r_0)$.
Since $\rho$ is a real closed $(1,1)$-form on $X$,
we can write $\rho=\omega_1-\omega_2$ on $X$ as the difference of two K\"ahler forms, and $\omega_k=d\alpha_k$ on $B_i$ for some $1$-form $\alpha_k$ on $B_i$.
It is clear that $\int_{A_{i,r}}|\chi\circ\sigma_j\,(f_j\circ\sigma_j)^*\rho|\leq\|\chi\|_{C^0}\int_{A_{i,r}}(f_j\circ\sigma_j)^*(\omega_1+\omega_2)$ and we have
\[
\int_{A_{i,r}}(f_j\circ\sigma_j)^*\omega_k=\int_{\partial A_{i,r}}(f_j\circ\sigma_j)^*\alpha_k,\qquad k=1,2.
\]
But $f_j\circ\sigma_j$ converge smoothly outside $\cup_i\gamma_i$, so on $\partial A_{i,r}$, which is disjoint from $\cup_i\gamma_i$, we get
\[
\lim_{j\to\infty}\int_{\partial{A_{i,r}}}(f_j\circ\sigma_j)^*\alpha_k=\int_{\partial A_{i,r}}(F\circ\sigma_\infty)^*\alpha_k.
\]
Then using the fact that $\sigma$ is diffeomorphic outside $\cup_i\gamma_i$, we get
\[\int_{\partial A_{i,r}}(F\circ\sigma_\infty)^*\alpha_k
=\int_{\partial (\sigma_\infty(A_{i,r}))}F^*\alpha_k=\int_{\sigma_\infty(A_{i,r})}F^*\omega_k.\]
We explain the Stokes formula used in the second equality. $\sigma_\infty(A_{i,r})$ has two components $B_{i,r}^+,\, B_{i,r}^-$, which are both biholomorphic to the disc. So we apply Stokes formula on each $B_{i,r}^{\pm}$:  $\int_{\partial B_{i,r}^\pm}F^*\alpha_k=\int_{B_{i,r}^\pm}F^*\omega_k$. This proves the second equality.
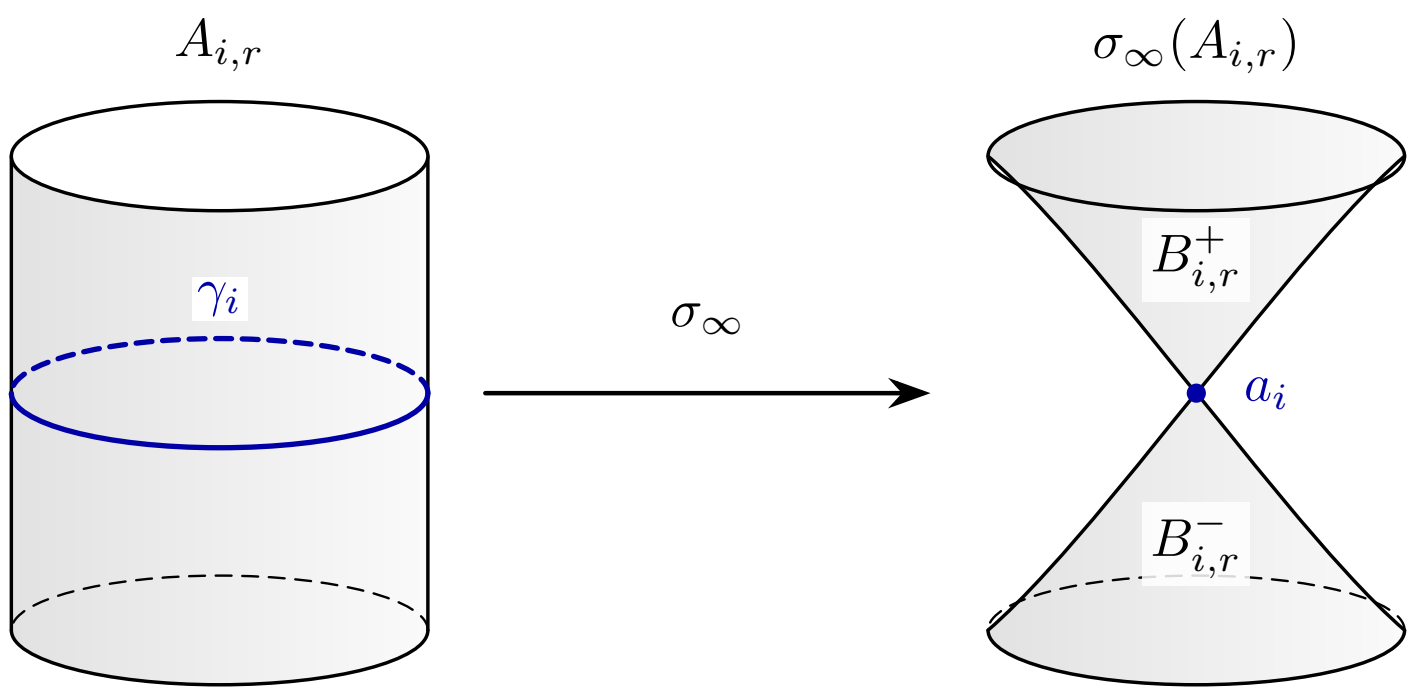
\begin{figure}[htbp]
\centering
\begin{tikzpicture}[
    scale=0.95,
    line width=0.65pt,
    line cap=round,
    line join=round,
    every node/.style={font=\small},
    core/.style={
        draw=blue!65!black,
        line width=1pt
    },
    hidden/.style={
        densely dashed,
        line width=0.45pt
    },
    branch label/.style={
        fill=white,
        fill opacity=0.85,
        text opacity=1,
        inner sep=1.5pt
    }
]

\begin{scope}[xshift=-3.4cm]

    \shade[left color=gray!23,right color=gray!4]
        (-1.45,1.65)
        arc[start angle=180,end angle=0,
            x radius=1.45cm,y radius=0.38cm]
        -- (1.45,-1.65)
        arc[start angle=0,end angle=-180,
            x radius=1.45cm,y radius=0.38cm]
        -- cycle;

    \draw (-1.45,-1.65) -- (-1.45,1.65);
    \draw ( 1.45,-1.65) -- ( 1.45,1.65);

    \filldraw[fill=white]
        (0,1.65) ellipse (1.45 and 0.38);

    \draw[hidden] (1.45,-1.65)
        arc[start angle=0,end angle=180,
            x radius=1.45cm,y radius=0.38cm];
    \draw (-1.45,-1.65)
        arc[start angle=180,end angle=360,
            x radius=1.45cm,y radius=0.38cm];

    \draw[core,densely dashed] (1.45,0)
        arc[start angle=0,end angle=180,
            x radius=1.45cm,y radius=0.38cm];
    \draw[core] (-1.45,0)
        arc[start angle=180,end angle=360,
            x radius=1.45cm,y radius=0.38cm];

    \node[
        text=blue!65!black,
        above=3pt,
        fill=white,
        inner sep=1pt
    ] at (0,0.38) {$\gamma_i$};

    \node at (0,2.42) {$A_{i,r}$};

\end{scope}

\draw[
    -{Stealth[length=2.8mm,width=2mm]},
    line width=0.85pt
]
    (-1.55,0) -- (1.55,0)
    node[midway,above=7pt] {$\sigma_\infty$};

\begin{scope}[xshift=3.4cm]

    \foreach \sgn in {-1,1}{
        \begin{scope}[yscale=\sgn]

            \shade[left color=gray!23,right color=gray!4]
                (-1.45,1.65)
                .. controls (-0.95,1.22) and (-0.20,0.23)
                    .. (0,0)
                .. controls (0.20,0.23) and (0.95,1.22)
                    .. (1.45,1.65)
                arc[start angle=0,end angle=180,
                    x radius=1.45cm,y radius=0.38cm]
                -- cycle;

            \draw (-1.45,1.65)
                .. controls (-0.95,1.22) and (-0.20,0.23)
                    .. (0,0)
                .. controls (0.20,0.23) and (0.95,1.22)
                    .. (1.45,1.65);

        \end{scope}
    }

    \draw (0,1.65) ellipse (1.45 and 0.38);

    \draw[hidden] (1.45,-1.65)
    arc[start angle=0,end angle=180,
        x radius=1.45cm,y radius=0.38cm];

    \draw (-1.45,-1.65)
    arc[start angle=180,end angle=360,
        x radius=1.45cm,y radius=0.38cm];

    \node[branch label] at (0, 0.93) {$B_{i,r}^{+}$};
    \node[branch label] at (0,-1.05) {$B_{i,r}^{-}$};

    \fill[blue!65!black] (0,0) circle (1.9pt);
    \node[text=blue!65!black,right=5pt]
        at (0,0) {$a_i$};

    \node at (0,2.42) {$\sigma_\infty(A_{i,r})$};

\end{scope}

\end{tikzpicture}
\caption{$\sigma_\infty$ collapses
the circle $\gamma_i\subset A_{i,r}$ to the nodal point $a_i$.}
\label{fig:local-pinching}
\end{figure}

When $r\downarrow0$, $\sigma_\infty(A_{i,r})$ contracts to a single point $\sigma_\infty(\gamma_i)=a_i$, hence $\lim_{r\downarrow0}\int_{\sigma_\infty(A_{i,r})}F^*\omega_k=0$. This proves \eqref{eq:neck-analysis}.

Finally, we prove $(4)$. Statement~$(3)$ gives
$\int_C g^*\chi F^*\omega=\lim_{j}\int_\D \chi f_j^*\omega\leq\liminf_j\int_\D f_j^*\omega=\liminf_jA(f_j)$ for every cut-off $0\leq\chi\leq1$ compact supported in $\D$.
Since $\chi$ is arbitrary, $A_*(C)\leq\liminf_jA(f_j)$ follows as we want.
The proof is therefore complete.
\end{proof}

\vspace{0.2cm}
\section{Disc functionals and the pseudoeffective threshold}
\label{app:psef threshold}
Let $(X,\omega)$ be a compact K\"ahler manifold, we define the pseudoeffective threshold $t_*$ of adjoint class by
\[
t_*:=\inf\{\,t>0:[K_X+t\omega]\text{ is pseudoeffective}\}.
\]
In this appendix, we construct compact K\"ahler manifold $(X,\omega)$ so that
\[
\sup_{f\in\mathcal A(X,x)}
Q_{\theta-\delta\omega}(f)
=
\begin{cases}
+\infty,&0<\delta<2\pi t_*,\\
0,&\delta\ge 2\pi t_*.
\end{cases}
\qquad\forall\,x\in X.
\]

A holomorphic vector bundle $V$ on a complex manifold $Y$ is called unitary flat if the transition matrices are locally constant with values in unitary matrices, see e.g. \cite[Section~2.2]{Beauv13}.
On complex torus $T=\mathbb{C}^m/\Lambda$ we can construct unitary flat bundles in the following manner. Let $\lambda_1,\cdots,\lambda_{2m}$ be a basis of lattice $\Lambda$. If $\alpha$ is a real number, we define
\[
\chi:\Lambda\longrightarrow U(1),\qquad
\chi\!\left(\sum_{j=1}^{2m}k_j\lambda_j\right)
=e^{2\pi i\alpha k_1}.
\]
Let
\[
L_\chi:=(\mathbb C^m\times\mathbb C)/{\sim},
\qquad
(z,v)\sim(z+\lambda,\chi(\lambda)v).
\]
Then $L_\chi$ is a unitary flat line bundle on $T$.
\begin{proposition}
Let $T$ be a complex torus with flat K\"ahler form $\eta$ and $V$ be a unitary flat holomorphic vector bundle on $T$ of rank $r+1\geq2$. 
Let
\[
\pi:X=\mathbb P(V)\longrightarrow T
\]
be the projectivized bundle.
Fiberwise Fubini--Study forms define a global semipositive $(1,1)$-form $\beta$
on $X$. Put
\[
\omega=a\beta+b\pi^*\eta,
\qquad a,b>0.
\]
Then $X$ is a uniruled compact K\"ahler manifold,
\[
\theta:=\operatorname{Ric}(\omega)=(r+1)\beta.
\]
The pseudoeffective threshold $t_*=\frac{r+1}{2\pi a}$. 
Moreover, every $x\in X$ satisfies
\[
\sup_{f\in\mathcal A(X,x)}
Q_{\theta-\delta\omega}(f)
=
\begin{cases}
+\infty,&0<\delta<2\pi t_*,\\
0,&\delta\ge 2\pi t_*.
\end{cases}
\]
\end{proposition}
\begin{proof}
The Fubini--Study metric on $\mathbb{P}^{r}$ is $i\partial\bar\partial\log(1+|z_1|^2+\cdots+|z_r|^2)$ hence is invariant under unitary transformation.
Therefore by gluing Fubini--Study metrics on each trivialization $U_\alpha\times\mathbb{P}^r$ of $\mathbb{P}(V)$, we obtain a global semipositive $(1,1)$-form $\beta$.
It is then clear that $\omega:=a\beta+b\pi^*\eta$ is K\"ahler for $a,b>0$ and $\operatorname{Ric}(\omega)=(r+1)\beta$.

We prove the threshold $t_*=\frac{r+1}{2\pi a}$.
Using $-2\pi[K_X]=[\operatorname{Ric}(\omega)]=(r+1)[\beta]$,
we get $[K_X+t\omega]=(at-\frac{r+1}{2\pi})[\beta]+bt[\pi^*\eta]$.
If $t\geq\frac{r+1}{2\pi a}$ clearly $[K_X+t\omega]$ is pseudoeffective.
Conversely if a closed positive $(1,1)$-current $T$ represents $[K_X+t\omega]$, then using $\pi^*\eta^{\dim T+1}=0$ we get
\[
0\leq\int_XT\wedge\beta^{r-1}\wedge (\pi^*\eta)^{\dim T}
=(at-\frac{r+1}{2\pi})\int_X\beta^r\wedge(\pi^*\eta)^{\dim T}.
\]
Since $\int_X\beta^r\wedge(\pi^*\eta)^{\dim T}=(\int_{\mathbb{P}^r}\omega_{FS}^r)\cdot(\int_T\eta^{\dim T})>0$, we get $t\geq\frac{r+1}{2\pi a}$. This proves the threshold $t_*=\frac{r+1}{2\pi a}$.

Finally, for every $x\in X$, we pick an embedding $\mathbb{P}^1\to\mathbb{P}^r\simeq \pi^{-1}(\pi(x)),\,[z_1:z_2]\mapsto[z_1:z_2:0:\cdots:0]$.
Consider discs $f_\lambda(z):=[\lambda z:1:\cdots:0]$, $\lambda>0$, a direct computation shows
\[
Q_{\theta-\delta\omega}(f_\lambda)=\frac{r+1-a\delta}{2\pi}\int_\D\log(1/|z|)i\partial\bar\partial\log(1+\lambda^2|z|^2)=\frac{r+1-a\delta}{2}\log(1+\lambda^2).
\]
If $\delta<2\pi t_*=\frac{r+1}{a}$, then $\lim_{\lambda\to+\infty}Q_{\theta-\delta\omega}(f_\lambda)=+\infty$.
If $\delta\geq 2\pi t_*=\frac{r+1}{a}$, then $\theta-\delta\omega\leq 0$, $Q_{\theta-\delta\omega}(f)\leq0$ for any $f\in\mathcal{A}(X,x)$.
Moreover $\lim_{\lambda\downarrow0}Q_{\theta-\delta\omega}(f_\lambda)=0$ shows $\sup_{f\in\mathcal{A}(X,x)}Q_{\theta-\delta\omega}(f)=0$.
The proof is therefore complete.
\end{proof}

\end{document}